\documentclass[11pt, reqno]{article}

\usepackage{titlesec}

\titlespacing\section{0pt}{8pt plus 4pt minus 2pt}{0pt plus 2pt minus 2pt}
\titlespacing\subsection{0pt}{4pt plus 4pt minus 2pt}{0pt plus 2pt minus 2pt}
\titlespacing\subsubsection{0pt}{2pt plus 4pt minus 2pt}{0pt plus 2pt minus 2pt}

\usepackage{hyperref}

\usepackage{enumitem}
\usepackage{tikz}
\usepackage{array} 
\usepackage{amsfonts}
\usepackage{amsmath}
\usepackage{amssymb}
\usepackage{graphicx}
\usepackage{color}
\usepackage{esint}    
\usepackage{cancel} 
\usepackage{enumitem} 

\newtheorem{thm}{Theorem}[section]
\newtheorem{cor}[thm]{Corollary}
\newtheorem{lem}[thm]{Lemma}
\newtheorem{prop}[thm]{Proposition}
\newtheorem{defn}[thm]{Definition}
\newtheorem{rem}[thm]{Remark}
\numberwithin{equation}{section}

\newcommand{\dx}{\,{\rm d}x}
\newcommand{\dy}{\,{\rm d}y}

\newcommand{\dt}{\,{\rm d}t}

\newcommand{\rd}{{\rm d}}

\def\LL{\mathrm{L}} 

\newcommand{\AM}{\mathcal{L}^{\frac{1}{2}}}
\newcommand{\AI}{\mathcal{L}^{-1}}
\newcommand{\AIM}{\mathcal{L}^{-\frac{1}{2}}}

\newcommand{\F}{\mathcal{F}}
\newcommand{\J}{\mathcal{J}}

\newcommand{\ka}{\overline{\kappa}}
\newcommand{\kb}{\underline{\kappa}}

\newcommand{\RR}{\mathbb{R}}
\newcommand{\NN}{\mathbb{N}}
\newcommand{\G}{\mathbb{G}_\Omega}

\renewcommand{\c}{\mathtt{c}}
\newcommand{\E}{\mathcal{E}}

\newcommand{\A}{\mathcal{A}}

\newcommand{\EL}{\mathsf{E}}
\newcommand{\IL}{\mathsf{I}}
\newcommand{\R}{\mathsf{R}}
\newcommand{\Q}{{\mathcal{Q}}}
\newcommand{\QL}{\mathsf{Q}}

\newcommand{\V}{\mathtt{V}}

\newcommand{\pkj}{\phi_{k,j}}

\def\ee{\mathrm{e}} 
\def\dist{\mathrm{dist}} 
\def\dom{\mathrm{dom}} 

\def\qed{\,\unskip\kern 6pt \penalty 500
\raise -2pt\hbox{\vrule \vbox to8pt{\hrule width 6pt
\vfill\hrule}\vrule}\par}
\definecolor{darkblue}{rgb}{0.05, .05, .9}
\definecolor{darkgreen}{rgb}{0.1, .65, .1}
\definecolor{darkred}{rgb}{0.8,0,0}

\begin{document}
\title{\vspace{-2.5cm}\textbf{The Cauchy-Dirichlet Problem\\ for the Fast Diffusion Equation\\ on Bounded Domains\\[7mm]}}

\author{Matteo Bonforte$^{\,a}$ and Alessio Figalli$^{\,b}$\\[3mm]
    \textit{In memoria di Emmanuele DiBenedetto}\\[3mm]
    }
\date{} 

\maketitle

\thispagestyle{empty}

\begin{abstract}
The Fast Diffusion Equation (FDE) $u_t= \Delta u^m$, with $m\in (0,1)$, is an important model for singular nonlinear (density dependent) diffusive phenomena. Here, we focus on the Cauchy-Dirichlet problem posed on smooth bounded Euclidean domains. In addition to its physical relevance, there are many aspects that make this equation particularly interesting from the pure mathematical perspective. For instance: mass is lost and solutions may extinguish in finite time, merely integrable data can produce unbounded solutions, classical forms of Harnack inequalities (and other regularity estimates) fail to be true, etc.

In this paper, we first provide a survey (enriched with an extensive bibliography) focussing on the more recent results about existence, uniqueness, boundedness and positivity (i.e., Harnack inequalities, both local and global), and higher regularity estimates (also up to the boundary and possibly up to the extinction time).   We then prove new global (in space and time) Harnack estimates in the subcritical regime.  In the last section, we devote   a special attention to the asymptotic behaviour, from the first pioneering results to the latest sharp results, and we present some new asymptotic results in the subcritical case.
\end{abstract}
\vspace{.5cm}

\noindent {\bf Keywords. } Singular nonlinear diffusion, Regularity, Harnack Inequalities, extinction profiles, asymptotic behaviour, convergence rates, Poincar\'e inequalities, entropy methods. \\[2mm]
{\sc Mathematics Subject Classification}.   {\sc 35K55, 35K67, 35B40, 35P30, 35J20.}

\vfill

\begin{itemize}[leftmargin=*]\itemsep1pt \parskip1pt \parsep0pt
\item[(a)] Departamento de Matem\'{a}ticas, Universidad Aut\'{o}noma de Madrid,\\
		ICMAT - Instituto de Ciencias Matem\'{a}ticas, CSIC-UAM-UC3M-UCM, \\
		Campus de Cantoblanco, 28049 Madrid, Spain.\\
		E-mail:\texttt{~matteo.bonforte@uam.es }\;\;\\
		Web-page: \texttt{http://verso.mat.uam.es/\~{}matteo.bonforte}
\item[(b)] ETH Z\"urich, Department of Mathematics,\\
R\"amistrasse 101, 8092 Z\"urich, Switzerland.\\
E-mail:\texttt{~alessio.figalli@math.ethz.ch}\\
Web-page:\texttt{~https://people.math.ethz.ch/\~{}afigalli}
\end{itemize}
\smallskip\noindent {\sl\small\copyright~2022 by the authors. This paper may be reproduced, in its entirety, for non-commercial purposes.}
\normalsize

\newpage

{\footnotesize \itemsep0pt \parskip0pt \parsep0pt
\tableofcontents}
\normalsize

\newpage
\section{Introduction}
The Fast Diffusion Equation
\begin{equation}\label{FDE}\tag{FDE}
u_t= \Delta u^m,\qquad m \in (0,1),
\end{equation}
is an important model for singular nonlinear (density dependent) diffusive phenomena. It was used by Barenblatt-Zel’dovich-Kompaneets to model gas-kinetics, then for thin liquid film dynamics by Van der Waals forces, and also for plasma (in nuclear reactors) when the high temperatures change the standard physical laws (Fourier's Law); see for instance the monographs by Barenblatt \cite{B1979} and V\'azquez \cite{VazBook}.

In their pioneering work \cite{BH}, Berryman-Holland introduced the Cauchy-Dirichlet problem for the \eqref{FDE} on bounded Euclidean domains and with homogeneous lateral boundary conditions, and analyzed its asymptotic behaviour. This was motivated by an experiment by Drake-Greenwood-Navratil-Post \cite{DGNP1977} on anomalous diffusion of hydrogen plasma across an octupole magnetic field, where the density of
the plasma satisfies the \eqref{FDE} with $m = 1/2$, also studied by Okuda-Dawson \cite{OD1973}. These experiments revealed that, after a few milliseconds, the solution evolves into a fixed shape which then decays in time and eventually extinguish.  We refer to the first chapters of \cite{VazBook} for a more detailed and historically rich description of this model. Since the 70s this model has been thoroughly studied  and we can say that, although today the theory is quite complete at least for the Cauchy problem on the whole space, see
\cite{VazLN, BBDGV,  BDGV, BGV2010, DKM2015,BDNS2022}, still many interesting questions remain open. In particular, for the Cauchy-Dirichlet problem on bounded domains,
\begin{equation}\label{CDP}\tag{CDP}
		\left\{\begin{array}{lll}
			u_t(t,x)=\Delta u^m(t,x) & \qquad\mbox{in }(0,+\infty)\times\Omega\,,\\
			u(t,x)=0 & \qquad\mbox{on }(0,+\infty)\times\partial\Omega\,,\\
			u(0,\cdot)=u_0  & \qquad \mbox{in } \Omega\,,
		\end{array}\right.
\end{equation}
where $0<m<1$ and $\Omega\subset \RR^N$ is a bounded smooth domain (at least of class $C^{2,\alpha}$), many basic questions remained open until very recently.  
We shall discuss below the many aspects that make this equation particularly interesting also from the pure mathematical perspective: for instance, mass is lost and solutions may extinguish in finite time, merely integrable data can produce unbounded solutions, classical forms of Harnack inequalities (and other regularity estimates) fail to be true, etc.

\noindent\textbf{Plan of the paper. }This paper is partly a survey but it contains also some new results, see in particular Sections \ref{SSec.new.proof.BC}, \ref{Sec.GHP.new}, and \ref{Ssec.Asymp.Subcrit.New}. We begin with a short survey about the theory of the Cauchy-Dirichlet problem \eqref{CDP} for the FDE, focussing on the more recent results about existence, uniqueness (Section \ref{Sec.exist.uniq}), the time monotonicity estimates of Benilan-Crandall (Section \ref{Sec.BC.ES}, with a new proof in Section \ref{SSec.new.proof.BC}), estimate on the extinction time and norm decay (Section \ref{Sec.FET} and new results in Section \ref{Sec.Rayleigh}), smoothing effects (Section \ref{Sec.Smoothings}), boundedness and positivity (local and global Harnack inequalities, see Section \ref{Sec.Harnack}, and new results in Section \ref{Sec.GHP.new}), and higher regularity estimates (up to the boundary and possibly up to the extinction time, see Section \ref{Ssec.Regularity}). We also prove some new results that we put in context below. In the last part of the paper, Section \ref{Sec.Asypt.beh},  we devote a special attention to the asymptotic behaviour, from the first pioneering results up to the most recent and up-to-date. We explain in some detail the proof of the sharp convergence rates towards equilibrium in generic domains, given by us in \cite{BF2021}, and we also comment on alternative proofs that appeared recently \cite{JX2019,JX2021,Ak2021,CMCS2022}. We also discuss the critical case, related to the Yamabe flow in Section \ref{Ssec.Yamabe}. Finally, we prove some new asymptotic results in the subcritical case $m<m_s$, which essentially represent the only available results in this range, see Section \ref{Ssec.Asymp.Subcrit.New}.


\section{Existence, uniqueness, different concepts of solutions. }\label{Sec.exist.uniq}Existence and uniqueness is related to the concept of solution that we want to consider.

\noindent\textbf{Smooth solutions. }Starting from sufficiently smooth data it is possible to establish existence of smooth solutions: the first result is done by Sabinina \cite{S1962} with $u_0\in C^1_0(\Omega)$. A frequently used approach to produce classical solutions (i.e., $C^{1+\alpha/2,2+\alpha}_{t,x}$ in the interior) is through smooth approximations that are solutions to the so-called ``lifted problem'', see e.g. \cite{VazBook, BFR2017} (this approach eliminates the degeneracy/singularity in the equation). Another approach is by fixed point arguments, where it is possible to obtain even smoother solutions, see \cite[Theorem 3.2]{JX2019} and \cite{JX2021}. Indeed, bounded positive weak  solutions are $C^\infty_{t,x}$ in the interior by standard parabolic theory, as we shall discuss later. Smooth solutions are unique as a consequence of the maximum/comparison principle.

\noindent\textbf{Weak solutions. }Starting from less regular data, there are several concepts of weak solutions that appear in the literature: weak, weak energy, very weak, distributional, mild-$\LL^1$, mild-$H^{-1}$, etc., and we refer to \cite{VazBook,VazLN,BSV2015,BV-NA,BII2022,AnSm05} for more details. Recently, Weak Dual Solutions (WDS) have proven to be a convenient setup. They are a class of ``limit solutions'': they can be constructed as non-decreasing limit in the strong $\LL^1_{\Phi_1}$ topology\footnote{Here, $\LL^1_{\Phi_1}$ denotes the weighted $\LL^1$ space in $\Omega$, where the weight $\Phi_1 \asymp \dist(\,\cdot\,,\partial\Omega)$ is just any smooth extension of the distance to the boundary. However, it is often convenient to take $\Phi_1$ as the ground state of the operator -- the classical Laplacian $-\Delta$ in this case -- because it simplifies technical estimates, although the same results hold with any smooth extension of the distance to the boundary, since we consider only smooth domains,  see for instance \cite{BSV2015,BV-NA,BFV2018, BII2022} for a more detailed discussion.  } of non-negative approximations by any other concept of solution mentioned above (usually by mild $\LL^1$ or $H^{-1}$). The minimal one turns out to be unique and strong in $\LL^1_{\Phi_1}$, when the initial datum is nonnegative and in $\LL^1_{\Phi_1}$.  This is the concept of solution that allow for biggest class of non-negative initial data known so far.
\begin{defn}[Weak Dual Solutions]\label{Def.WDS}
Let $T>0$. We say that $u\in C((0,T):L^1_{\Phi_1}(\Omega))$ is a \emph{weak dual solution} of \eqref{CDP} if $u^m\in L^1\big([0,T]:L^1_{\Phi_1}(\Omega)\big)$ and
	\begin{align}\tag{WDS}\label{WDS}
	\int_{0}^T\int_{\Omega}  (-\Delta)^{-1}u\;\partial_t\psi\dx\dt=\int_{0}^{T}\int_{\Omega}u^m\psi\dx\dt\qquad\forall\,\psi\text{ s.t. } \frac{\psi}{\Phi_1}\in C^1_c((0,T):L^{\infty}(\Omega))\,.
	\end{align}
We say that $u$ is a WDS of the Cauchy-Dirichlet Problem \eqref{CDP}, corresponding to the initial datum $u_0\in L^1_{\Phi_1}(\Omega)$, if moreover $u\in C([0,T):L^1_{\Phi_1}(\Omega))$, and
$
\lim\limits_{t\to 0^+}\|u(t)-u_0\|_{L^1_{\Phi_1}(\Omega)}=0.
$
A WDS is called \emph{strong} if in addition $t\,u_t \in L^\infty((0,T) : L^1_{\Phi_1}(\Omega))$.
A WDS is called \emph{Minimal Weak Dual Solution} (MWDS) if it is obtained as the non-decreasing limit of a sequence of semigroup solutions.
\end{defn}\vspace{-2mm}
Note that this concept of solutions involves $(-\Delta)^{-1}$, the inverse of the Dirichlet Laplacian, whose kernel is the classical Green function, see Section \ref{Ssec2.1} for more details. Roughly speaking, WDS can be seen as distributional solution to the dual equation: $(-\Delta)^{-1}u_t=-u^m$, used probably for the first time by Pierre \cite{P1982} to prove uniqueness with measure data.   WDS were introduced by V\'azquez and one of the authors for the first time in \cite{BV-ARMA, BV-NA}, inspired by an elliptic analogous concept, introduced by Brezis (unpublished notes). In a recent paper \cite{BII2022}, Ibarrondo, Ispizua and one of the authors established existence, uniqueness and boundedness of strong MWDS in a more general framework (when  $\Delta$ is replaced by a possibly nonlocal operator $-\mathcal{L}$). The result of \cite[Theorems 2.4 and 2.5]{BII2022} read as follows in the present setting:  \it For every $0\leq u_0\in L^1_{\Phi_1}(\Omega)$, there exist a unique strong MWDS $u$ of \eqref{CDP} with
	\[
	\lim\limits_{t\rightarrow 0^+}\|u(t)-u_0\|_{L^1_{\Phi_1}\!(\Omega)}=0\quad\mbox{and}\quad \|u_t(t)\|_{L^1_{\Phi_1}\!(\Omega)}\leq \frac{2\|u_0\|_{L^1_{\Phi_1}\!(\Omega)}}{(1-m)\,t}\,.
	\]
Moreover, the T-contraction estimates hold: $\|(u(t)-v(t))_{\pm}\|_{L^1_{\Phi_1}(\Omega)}\leq \|(u_0-v_0)_{\pm}\|_{L^1_{\Phi_1}(\Omega)}$  for any $t\ge 0$ and any MWDS $u(t),v(t)$ corresponding to  $0\le u_0,v_0\in L^1_{\Phi_1}(\Omega)$. \rm

As far as nonnegative solutions are concerned, it can be proven that all the previous mentioned concepts of solutions (weak, mild or semigroup, $H^{-1}$) are indeed WDS, see \cite{BSV2015, BV-NA, BII2022}. When we deal with signed solutions, the biggest class known so far is obtained by gradient flow techniques via the Brezis-Komura Theorem \cite{Brezis1,BrezisBook,Komura} see also the excellent notes \cite{Ambrosio-Notes}: indeed, the FDE is a gradient flow in $H^{-1}$ of the $L^{1+m}$-norm. In \cite{BII2022} nonnegative WDS are obtained as increasing limits of nonnegative $H^{-1}$ solutions. Indeed, $H^{-1}$-solutions are strong in the $H^{-1}$ sense, i.e., $u_t\in H^{-1}$, they are unique, and they can be shown to be weak energy solutions, i.e., $u^m\in H^1_0$.

\noindent\textit{Nonlinear contractive semigroups. }Notice that the FDE generates a nonlinear contractive semigroup only in a few spaces: in $\LL^1$ by the theory of Benilan-Crandall-Pazy-Pierre \cite{BCPbook, CP-JFA, VazBook}, in $H^{-1}$ by the (gradient flow) theory of Brezis and Komura \cite{Ambrosio-Notes, BrezisBook, Brezis1,Komura, VazBook}, and with respect to the 2-Wasserstein distance as first observed by Otto \cite{Otto}, see also \cite{AGSbook,IPS2019}. In some $\LL^p$ spaces there can be still contractivity, but the relation between $p,m$ and $N$ is quite involved, as observed by Chmaycem-Jazar-Monneau \cite{CJM2016}. Moreover, it is well known that these semigroups have the stronger $T$-contractive property in $\LL^1$ and $H^{-1}$, namely $\|(u(t)-v(t))_{\pm}\|\le \|(u(0)-v(0))_{\pm}\|$,  
which implies uniqueness and comparison: ordered initial data produce unique and ordered solutions. More recently, T-contraction in $\LL^1_{\Phi_1}$, hence uniqueness and comparison, has been shown to hold also for MWDS, see \cite{BII2022}.

\noindent\textbf{Initial traces. }The problem of existence and uniqueness of initial traces (i.e., the weak limit of solutions as $t\to 0^+$) has been solved by Dahlberg-Kenig, see \cite{Da-Ke,DaskaBook}.

\noindent\textbf{Optimal class of data. }We have seen that there are many different concepts of solutions.  
It shall be reminded that for the Cauchy problem on the whole space, Herrero-Pierre \cite{HP1985} showed existence and uniqueness of very weak solution with data $u_0\in \LL^1_{\rm loc}(\RR^N)$. We shall try to use here the most general concept of solution, depending on the initial datum: to the best of our knowledge, $\LL^1_{\Phi_1}(\Omega)$ and $H^{-1}(\Omega)$ are the biggest spaces where existence, uniqueness and comparison hold, if one considers nonnegative or signed solutions respectively. Solutions corresponding to data in $\LL^1_{\Phi_1}(\Omega)$ and $H^{-1}(\Omega)$ generate unique smooth solutions in most of the cases: only when $m$ is in the very fast diffusion range, we shall require further integrability of the initial datum to obtain bounded solutions. Since we will mainly deal with nonnegative data, from now on we will consider MWDS corresponding to $0\le u_0\in \LL^1_{\Phi_1}(\Omega)$. We refer to the book of V\'azquez, Chapters 6.6 and 6.7 , for more details about the setup in $H^{-1}$ (originally due to Brezis \cite{Brezis1}) or in $\LL^1_{\dist(\,\cdot\,,\partial\Omega)}$. For the nonlinear semigroup theory we refer to the book of Benilan-Crandall-Pazy \cite{BCPbook} and the original papers of Crandall-Ligget \cite{CL1971}, Benilan-Crandall \cite{BCr, BCr2}, and Crandall-Pierre \cite{CP-JFA}, among many others.

\noindent\textbf{Related problems. }\textit{The limit $m\to0^+$. }In this case, there are two possible limiting equations: the logarithmic diffusion $u_t=\Delta\log(u)$ and the signed fast diffusion $u_t=\Delta {\rm sign}(u)$, as shown by the authors in \cite[Section 3]{BF2012}. In the logarithmic case, nonnegative solutions of the \eqref{CDP} fail to exist, as first proven by V\'azquez \cite{V1}, see also \cite{BoSeVa}. In the other case, the authors \cite{BF2012} provide the explicit dynamic of solutions to the Total Variation Flow when $N=1$, and of their distributional derivatives, which solve  the signed fast diffusion $u_t=\Delta {\rm sign}(u)$\,.

\noindent\textit{Neumann boundary conditions. }These have been considered by Iacobelli \cite{I2019} and by Iacobelli-Patacchini-Santambrogio \cite{IPS2019}, for ultrafast diffusion (i.e., when $m\le 0$) possibly with weights. This model comes from the quantization problem \cite{I2019}, and a quite complete theory of existence, uniqueness boundedness and asymptotic behaviour has been developed in \cite{IPS2019}.

\noindent\textit{Dynamical boundary conditions. } These are considered by Schimperna-Segatti-Zelik \cite{SSZ2016}, where they prove existence, uniqueness, and a number of estimates aimed to clarify the asymptotic behaviour.


\section{The time monotonicity estimates of Benilan and Crandall}\label{Sec.BC.ES}
The celebrated Benilan-Crandall estimates \cite{BCr} for solution to the FDE $u_t=\Delta u^m$ with $m\in (0,1)$ hold in the distributional sense and read
\begin{equation}\label{Be-Cr}
u_t\le\frac{u}{(1-m)t}\,.
\end{equation}
Indeed this is just the weak formulation of the following time monotonicity:
\begin{equation*}
t\mapsto t^{-\frac{1}{1-m}}u(t,x) \qquad\mbox{is monotone non-increasing for a.e. }x\in \Omega\,.
\end{equation*}
This monotonicity plays a key role in regularity estimates. See \cite{BCr,VazBook} and also \cite{BII2022} for more details.

\noindent\textbf{Proof by Scaling and Comparison. }We shall first present what is probably the simplest proof of the Benilan-Crandall estimates. To the best of our knowledge, this has been first shown in \cite{BCr}.

Take $\lambda\ge 1$ and consider the rescaled (in time) solution to \eqref{CDP}
\[
    u_{\lambda}(t, x)=\lambda^{-\frac{1}{1-m}}u(\lambda t, x),
        \qquad\mbox{with}\qquad u_{\lambda}(0, x)=\lambda^{-\frac{1}{1-m}}u_0(x)\le u_0(x)\,,
\]
By comparison, since $\lambda\ge 1$, it follows that $u_{\lambda}(t, x)\le u(t,x)$ for almost all $t\ge 0$ and $x\in \Omega$. On the other hand, letting now  $\lambda=\frac{t+h}{t}\ge 1$, we obtain
\[\begin{split}			
\frac{u(t+h, x)-u(t, x)}{h}
			&= \frac{1}{h}\left[\left(\frac{t+h}{t}\right)^{\frac{1}{1-m}}-1\right]u_{\lambda}(t, x) +\frac{u_{\lambda}(t, x)-u(t, x)}{h}
			 \leq \frac{(t+h)^{\frac{1}{1-m}}-t^{\frac{1}{1-m}}}{h}  \frac{u_{\lambda}(t, x)}{t^{\frac{1}{1-m}}}.
\end{split}
\]
Taking limits as $h\to 0^+$ (in the distributional sense), we obtain \eqref{Be-Cr}.

\noindent\textbf{Proof by maximum principle. }When the solution  is classical there is another simple proof purely based on the maximum principle: Define
\[
w(t,x):=u(t,x) - (1-m)t u_t(t,x)\qquad\mbox{that satisfies the equation}\qquad w_t= m \Delta\left(u^{m-1} w\right)\,.
\]
Since $w(0,x)=u_0\ge 0$ and $w(t,x)=0$ on $(0,\infty)\times \partial\Omega$\,, and the equation satisfies the maximum principle, we conclude that $w(t,x)\ge 0$ for a.e. $t>0$ and $x\in \Omega$, which is exactly the pointwise version of \eqref{Be-Cr}. One can easily extend this proof to more general solutions, for instance with merely integrable data, by standard approximation techniques and uniqueness.

We refer to the works of Benilan-Crandall-Pierre \cite{BCr, CP-JFA} and to  V\'azquez's book \cite[Chapter 8]{VazBook}, where different (rigorous) proofs are collected in different frameworks and for different classes of solutions.

\noindent\textbf{More general nonlinearities. }The above time monotonicity estimates \eqref{Be-Cr} can be generalized to solutions of the so-called \textit{filtration equation }$u_t=\Delta\varphi(u)$, where $\varphi$ is allowed to be a non-homogeneous nonlinearity satisfying suitable conditions. This was done for the first time by Crandall-Pierre in \cite{CP-JFA}, where a third different proof of \eqref{Be-Cr} can be found.

The Benilan-Crandall estimates \eqref{Be-Cr} are a key ingredient in many regularity theories for the FDE and in general for Porous Medium type equations, see \cite{BCPbook,BII2022,BV-NA,DaskaBook,FM2017,VazBook}. Also, they are the key ingredient in the ``almost representation formula'' \eqref{Prop.PE.ineq}, see Section \ref{Ssec2.1}.

\subsection{A proof of the Benilan-Crandall estimates by absorption}\label{SSec.new.proof.BC}
In a recent paper by Jin-Xiong \cite{JX2019}, it appeared another proof of that Benilan-Crandall estimates, which relies on the equation satisfied by the function $w=u_t/u$, which is a curvature-like quantity. In particular $w$ satisfies a ``nice'' equation with weights that is used there to study optimal boundary regularity, as we discuss in Section \ref{Ssec.Regularity}. Inspired by \cite{JX2019}, we have found a simplified proof, that exploits the absorption term in the equation for $w$. By Kato's inequality, the positive and negative parts of solutions are subsolutions, in particular, we have that
\[
v(t,x):=(u_t(t,x))_+\qquad\mbox{satisfies}\qquad v_t \le m\Delta\big(u^{m-1}v\big)\,.
\]
Define now the positive part of the ``curvature''
\[
w(t,x):=\frac{v(t,x)}{u(t,x)}=\frac{(u_t(t,x))_+}{u(t,x)}\qquad\mbox{that satisfies}\qquad
w_t \le \frac{m}{u} \Delta\left(u^m w\right)-w^2\,.
\]
Hence $w$ is a nonnegative subsolution to the same equation as $v$, but with an extra absorption term, that we are going to exploit.
We consider the following ``weighted norms'' for $q>1$:
\[
N_q[w](t) = \int_\Omega w^q(t,x) u(t,x) \dx\,.
\]
The time derivative along the flow has the expression:
\[\begin{split}
\frac{\rd}{\dt}N_q[w]&(t) = q\int_\Omega w^{q-1}w_t u \dx + \int_\Omega w^q u_t \dx
 \le -qm\int_\Omega \nabla\left(w^{q-1}\right)\cdot\nabla\left(u^m   w \right)\dx
-(q-1)\int_\Omega w^{q+1} u\dx\\
&=-q(q-1)m\int_\Omega w^{q-2}  |\nabla w|^2 u^m\dx
+m(q-1)   \int_\Omega  w^q \Delta u^m \dx -(q-1)\,\int_\Omega w^{q+1} u\dx\\
&= - \frac{4 m (q-1)}{q}\int_\Omega   \left|\nabla w^{\frac{q}{2}}\right|^2 u^m\dx
-(1-m)(q-1)\int_\Omega    w^{q+1} u \dx\le -(1-m)(q-1)\int_\Omega    w^{q+1} u \dx\,,
\end{split}\]
where we have integrated by parts, used that $u\ge 0$, and that $\Delta u^m =  u_t \le (u_t)_+= v=  w u$. Using H\"older inequality and the monotonicity of the $\LL^1$ norm (namely $\|u(t)\|_{\LL^1(\Omega)}\le \|u_0\|_{1}$) as follows,
\[
\int_\Omega    w^{q+1} u \dx \ge \left(\int_\Omega u \dx\right)^{-\frac{1}{q}} \left(\int_\Omega    w^q u \dx\right)^{\frac{q+1}{q}}
\ge \frac{N_q[w]^{\frac{q+1}{q}}}{\|u_0\|_{\LL^1(\Omega)}^{\frac{1}{q}}}\,,
\]
we deduce the following differential inequality
\begin{equation}\label{Nq.diff.ineq}
\frac{\rd}{\dt}N_q[w](t)\le -\frac{(1-m)(q-1)}{\|u_0\|_{\LL^1(\Omega)}^{\frac{1}{q}}}N_q[w]^{\frac{q+1}{q}}\,.
\end{equation}
Integrating the above differential inequality \eqref{Nq.diff.ineq} on $[0,t]$, we obtain
\begin{equation}\label{smoothing.q.1.w}
  N_q^{\frac{1}{q}}[w](t)\le  \frac{q}{(q-1)(1-m)} \frac{\|u_0\|_{\LL^1(\Omega)}^{\frac{1}{q}}}{t}\;.
\end{equation}
This inequality can be seen as a generalized Benilan-Crandall inequality. Indeed, when $q \to \infty$ the left-hand side converges to $ \|w(t)\|_{\infty}$, and
 \eqref{smoothing.q.1.w} implies
\[
\frac{u_t(t,x)}{u(t,x)}\le  \frac{(u_t(t,x))_+}{u(t,x)}\le \|w(t)\|_{\infty}\le \frac{1}{(1-m)t}\,,
\]
which is exactly the Benilan-Crandall inequality \eqref{Be-Cr}. \newpage

 This proof holds for sufficiently regular solutions, then by standard approximation techniques we can extend the inequality (in the distributional sense) up to Weak Dual Solutions. One possibility  is to  approximate WDS by means of $C^{2,3}_{t,x}$-smooth solutions for which the above proof is rigorous and inequality \eqref{Be-Cr} holds for all $x\in \overline{\Omega}$, and all $t>0$. Such solutions exist and are unique when the initial datum belong to a special class (dense in the cone of nonnegative functions of $L^1_{\Phi_1}(\Omega)$), see Theorem 3.4 of \cite{JX2019} for a proof based on a fixed point theorem. A (simpler) rigorous proof the Benilan-Crandall inequality \eqref{Be-Cr} that holds for WDS (also the context of nonlocal diffusions) can be obtained by scaling and comparison, as explained above, see also Section 5.2 of \cite{BII2022}.


\section{Local VS global smoothing effects}\label{Sec.Smoothings}
Local upper bounds for local solutions to the FDE take the form
\begin{equation}\label{Smoothing.Local}
 \sup_{x\in B_{R/2}}u(t,x)
 \le\frac{c_{1}}{t^{ N\vartheta_p}}\,\left[\int_{B_{  R}}|u_0(x)|^p\dx\right]^{ 2\vartheta_p}
 +c_{2}\left[\frac{t}{R^2}\right]^{\frac{1}{1-m}}\,,
\end{equation}
where $\vartheta_p=1/(2p-N(1-m))=1/2(p-p_c)$, and the constants $c_{i}$ depend on $m,N$
and $p$. The above bounds hold for all $t,R>0$ under some restrictions on the integrability of the initial datum:
\begin{equation}\label{Hyp.u0.Lp}
\mbox{$u_0\in L^p(\Omega)$ with $p\ge 1$ if $m\in (m_c,1)$ or $p>p_c$ if $m\in (0,m_c]$.}
\end{equation}
Such estimates have been proven for different concepts of weak solutions by several authors:  by DiBenedetto, Gianazza, Vespri \cite{DiBook, DGVbook} using nonlinear extensions of the celebrated De Giorgi method, or through Moser iteration by Dahlberg, Daskalopoulos, Herrero, Kenig, Pierre, Simonov, V\'azquez and the first author \cite{HP1985, Da-Ke,DaskaBook,BV-ADV,BS2019}.
Two critical exponents appear naturally:
\[
m_c:= \frac{N-2}{N}\qquad\mbox{and}\qquad p_c:=\frac{N(1-m)}{2}\,.
\]
It is well known since the celebrated counterexample by Brezis-Friedman \cite{BF1983}, that \textit{in the very fast diffusion range, i.e., when $m<m_c$, locally integrable data may not generate bounded solutions. }The integrability condition $p>p_c$ is necessary to avoid concentration and   blow-up. In  the case of the Cauchy problem on the whole space, the fundamental solution does not exist anymore. More precisely, the solution corresponding to a Dirac delta at a point remains a measure until it extinguishes. Moreover, there exist explicit very singular solutions that behave like $(T-t)^{\frac{1}{1-m}}|x-x_0|^{-\frac{2}{1-m}}$, so that they remain unbounded until they extinguish, see V\'azquez's monograph \cite{VazLN}. These examples confirm that the condition $p>p_c$ is s necessary to obtain bounded solutions.

We notice that a second pair of exponents appear in the game: if we consider energy solutions, or even $H^{-1}$ solutions, it is natural to consider initial data with finite energy, i.e., $u_0\in \LL^{1+m}(\Omega)$. The question is now when such data produce bounded solutions: this happens when
\[
1+m>p_c \qquad\mbox{that is}\qquad m>m_s:=\frac{N-2}{N+2}\,.
\]
The exponent \textit{$m_s$ is commonly called Sobolev or Yamabe exponent: }it is the inverse of the critical exponent $p_s:=\tfrac{N+2}{N-2}=2^\ast-1$ in semilinear elliptic equations, and also the exponent of the Yamabe flow, see Sections \ref{Ssec.LEF} and \ref{Ssec.Yamabe}. As we shall see, this exponent plays a crucial role in the Cauchy Dirichlet problem, as it distinguishes the ``good behaviour'' from the ``bad behaviour'' both in terms of uniform regularity and asymptotic behaviour.

\color{black}
\begin{figure}[h]
	\centering
	\begin{tikzpicture}[xscale=0.9]
	\node[scale=1.2] (00) at (0,-0.5) {$\mathbf{0}$};
	\node (0) at (-0.25,0){};
	\node (0-1) at (0,0.3){};
	\node (0-2) at (0,-0.3){};
	
	\node[scale=1.2] (11) at (10,-0.5) {$\mathbf{1}$} ;
	\node (1) at (12,0){};
	\node (1-1) at (10,0.3){};
	\node (1-2) at (10,-0.3){};
	\node [scale=1.1](m) at (11.75,-0.5){${m}$};
	
	\path[line width=0.3mm,->] (0) edge (1);
	\path[line width=0.6mm,-] (0-1) edge (0-2);
	\path[line width=0.6mm,-] (1-1) edge (1-2);
	
	\node (ms) at (3.33,-0.5){$\mathbf{m_s}$};
	\node (ms-1) at (3.33,0.3){};
	\node (ms-2) at (3.33,-0.3){};
	
	\node (mc) at (6.67,-0.5){$\mathbf{m_c}$};
	\node (mc-1) at (6.67,0.3){};
	\node (mc-2) at (6.67,-0.3){};

	\path[line width=0.4mm,-] (ms-1) edge (ms-2);
	\path[line width=0.4mm,-] (mc-1) edge (mc-2);
	
	\node[fill=red!40,draw,text width=2cm] (Lp) at (1.6,1.3) {$\;\;L^p\to L^\infty$};
	\node[fill=cyan!30,draw,text width=2.5cm] (L1+m) at (5,1.3) {$\;\;L^{1+m}\to L^\infty$};
	\node[fill=green!30,draw,text width =2.5cm] (L1) at (8.3,1.3) { $\;\;\;\;L^1\to L^\infty$};
	\node[fill=yellow!40,draw,scale=1.2,text width=4.52cm] (H*) at (6.65,2.1) {\hskip 1.5cm\small$H^*\to L^\infty$};
	
	\node[scale=0.75] (p1) at (1.6,0.6) {$p>p_c>1+m$};
	\node[scale=0.75] (p2) at (5,0.6) {$1+m>p_c>1$};
	\node[scale=0.75] (p3) at (8.4,0.6) {$1>p_c>0$};
	
	\node[scale=0.9,draw,rounded corners,text width=2.75cm] (pc) at (-2.75,0.5) {$\begin{aligned}
		\mathbf{p_c}&=\frac{N(1-m)}{2}\\
		\mathbf{m_c}&=\frac{N-2}{N}\\
		\mathbf{m_s}&=\frac{N-2}{N+2}
		\end{aligned}$};
	\node (blank) at (14,0){};
	
	\node[scale=0.9,text width=5cm,align=center](Very FDE) at (3.33,-1.5){Very Fast Diffusion};
	\node[scale=0.9,text width=2.4cm,align=center](GoodFDE) at (8.35,-1.5){Good FDE};
	\node[scale=0.9,text width=2cm,align=center](HE) at (10,-2.2){Heat Eq.};
	\node[scale=0.9,text width=2cm,align=center](Very FDE) at (11.45,-1.5){PME};
	
	\node[font=\fontsize{10}{10}\selectfont,color=red!70](VFDE brace) at (3.33,-1) {$\underbrace{\hskip 6cm}$};
	\node[color=darkgreen!70](GoodFDE brace) at (8.35,-1) {$\underbrace{\hskip 2.8cm}$};
	\draw[line width=0.3mm,<-,color=orange!70] (HE) -- (11);
	\node[color=violet] (PME brace) at (11.37,-1) {$\underbrace{\hskip 2.2cm}$};


    \node[font=\fontsize{10}{10}\selectfont,color=red!85!black](VFDE overbrace) at (1.6,2.75) {$\overbrace{\hskip 2.75cm}$};
    \node[scale=0.9,text width=2.4cm,align=center](VFDE overbrace text) at (1.6,3.1) {Subcritical};

    \node[font=\fontsize{10}{10}\selectfont,color=green!70!black](GFDE overbrace) at (6.65,2.75) {$\overbrace{\hskip 5.75cm}$};
    \node[scale=0.9,text width=2.4cm,align=center](GFDE overbrace text) at (6.65,3.1) {Supercritical};

    \node[scale=0.9,text width=2.5cm,align=center](YE) at (3.33,3.8){Sobolev/Yamabe\\ exponent};
    \node[scale=0.9,text width=2cm,align=center](YE0) at (3.33,0.4){};
    \draw[line width=0.3mm,<-,color=orange!70] (YE) -- (YE0);

	\end{tikzpicture}

  \vspace{-4mm}\caption{In this figure we can appreciate the validity of the $\LL^p-\LL^\infty$ and $H^*-\LL^\infty$ smoothing effects in relation with the critical exponents, in the different fast diffusion regimes.}

\end{figure}
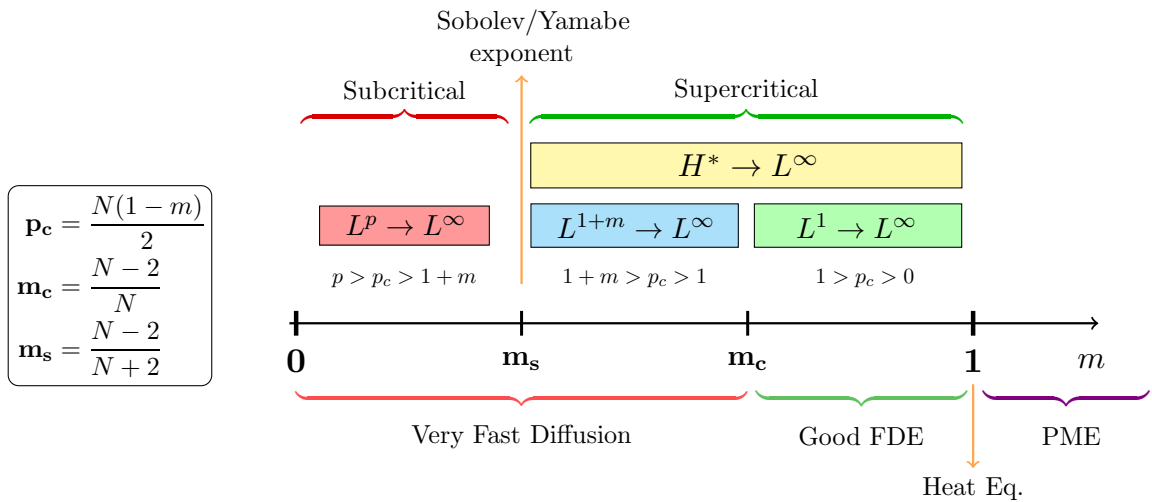

Let us also remark that many common concepts of weak solutions require $u^m\in H^1_0(\Omega)$.   By the standard Sobolev inequality
$\|f\|_{\LL^\frac{2N}{N-2}(\Omega)}\le \mathcal{S}_2 \|\nabla f\|_{\LL^2(\Omega)}$,    this implies $u\in \LL^{1+m}(\Omega)$ only when $m>m_s$.

 The Sobolev inequality is equivalent to the Hardy-Littlewood-Sobolev (HLS) inequality:
\begin{equation}\label{HLS}
\|f\|_{H^{-1}(\Omega)}\le \mathcal{S}_2\|f\|_{\LL^q} \qquad\mbox{for all $q\ge \frac{2N}{N+2}$}
\end{equation}
guarantees that when $u_0\in \LL^{p}(\Omega)$ with $p\ge p_c\wedge 1$ then $u_0\in H^{-1}(\Omega)$, hence when $m>m_s$, we have that $1+m>p_c$ and $u\in \LL^{1+m}(\Omega)$ implies $u\in H^{-1}(\Omega)$.

However, when dealing with a Cauchy-Dirichlet problem, global smoothing effects turn out to be true:  there exists a constant $\ka >0$ that depends only on $m,N,p$, such that for all $t>0$
\begin{equation}\label{Smoothing.Global}
 \|u(t)\|_{\LL^\infty(\Omega)}\le \ka\frac{\|u_0\|_{\LL^p(\Omega)}^{2p\vartheta_p}}{t^{ N\vartheta_p}}\qquad\mbox{with}\qquad\vartheta_p=\frac{1}{2p-N(1-m)}=\frac{1}{2(p-p_c)}\,.
\end{equation}
The above smoothing can be deduced from the local one \eqref{Smoothing.Local}, just extending the solution $u$ by zero outside $\Omega$ and letting $R\to\infty$.  This can also be proven directly as consequence of Gagliardo-Nirenberg-Sobolev inequalities (GNS) and Moser iteration \cite{BV-ADV,BS2019,BDNS2022,CH2021, BII2022, DaskaBook, IPS2019}, or through De Giorgi method \cite{DiBook,DGVbook}, or also by means of Logarithmic Sobolev inequalities via a nonlinear extension of Gross' method \cite{BG2006,BG2007} also in the context of Riemannian Manifolds \cite{BGV2008}.
A more recent technique, that we call Green function method\footnote{The so-called Green function method relies on duality, Green function estimates, and time monotonicity, and it was introduced for the Porous Medium case by  V\'azquez and the first author in \cite{BV-ARMA,BV-NA}.}, provides an alternative proof (without using GNS) of the above smoothing also in the more general framework of WDS, and it also allows to prove weighted smoothing effects:  there exists a constant $\ka >0$ such that for all $t>0$
\begin{equation}\label{Smoothing.Global.Weighted}
 \|u(t)\|_{\LL^\infty(\Omega)}\le \ka\frac{\|u_0\|_{\LL^p_{\Phi_1}(\Omega)}^{p\vartheta_{p,1}}}{t^{ N\vartheta_{p,1}}}
 \qquad\mbox{with}\qquad\vartheta_{p,1}=\frac{1}{p-N(1-m)}=\frac{1}{p-p_{c,1}}\,,
\end{equation}
 where $\ka$ only depends on $m,N,p,|\Omega|$.   A new pair of critical exponent appears naturally
\begin{equation}
\label{eq:mc1}
m_{c,1}:= \frac{N-1}{N}\qquad\mbox{and}\qquad p_{c,1}:= N(1-m) \,.
\end{equation}
We refer to \cite[Section 3]{BII2022} for complete proofs of the smoothing effects, both through Moser Iteration and through the Green function method (both  methods work in the present case),   together with a thorough explanation of all the critical exponents.

\noindent\textit{The ``Sobolev exponent on $\RR^N$''. }It is interesting to note that $m_{c,1}$ corresponds to the ``Sobolev exponent on $\RR^N$'' for the FDE. Indeed, fast diffusion equations on the whole space and Gagliardo-Nirenberg-Sobolev inequalities are deeply related, and $m_{c,1}$ correponds to the ``critical'' Sobolev case. We refer to the memoir by Dolbeault-Nazaret-Simonov and the first author \cite{BDNS2022} for a thorough discussion about this intriguing relation. There, constructive and quantitative stability estimates for GNS inequalities is obtained by means of a new flow method, based on entropy techniques and regularity estimates. Constructive proofs (with explicit constants) of local and global smoothing effects can be found there, many of which apply also to the Dirichlet case, directly or with minor modifications.\\
On a different point, Otto showed that solutions to the Cauchy problem for the FDE can be obtained as gradient flow solutions with respect to the 2-Wasserstein distance, see \cite{Otto, AGSbook}. Here, the exponent $m_{c,1}$ represents the threshold below which displacement convexity in 2-Wasserstein distance fails.

\noindent\textit{Boundedness when $u_0\in H^{-1}$. }When $m>m_s$, we can take $p=1+m>p_c$ in the smoothing effects \eqref{Smoothing.Global} and then bound the $L^{1+m}$ norm with the $H^{-1}$ norm as in  Proposition \ref{Prop.Q*-estimates}, to show that 
 there exists a constant $\ka >0$ that depends only on $m,N,|\Omega|$, such that for all $t> t_0\ge 0$
\begin{equation*}
\|u(t)\|_{\LL^\infty(\Omega)}\leq \ka\;\frac{\|u(t_0)\|^{4\,\vartheta_{1+m}}_{H^{-1}(\Omega)}}{(t-t_0)^{(N+2)\vartheta_{1+m}}}\qquad\qquad\mbox{with}\quad \vartheta_{1+m}=\frac{1}{(N+2)m - (N-2)}\,.
\end{equation*}
We conclude that initial data (without sign restriction) in $H^{-1}(\Omega)$ produce bounded solutions.
\subsection{Upper boundary estimates }
In \cite[Theorem 2.9]{BII2022} also boundary estimates are provided; more precisely, if $u_0$ satisfies \eqref{Hyp.u0.Lp} we have
\[
\left\|\frac{u^m(t)}{\Phi_1}\right\|_{\LL^\infty(\Omega)}
\le  \ka
\frac{\|u(t_0)\|_{L^p(\Omega)}^{2p\vartheta_p}}{(t-t_0)^{mN\vartheta_p+1}}.
\]
See Theorem \ref{thm.GHP.up} for the precise statement, also for more general data $u_0\in L^p_{\Phi_1}(\Omega)$. In this paper we extend the above estimates globally in time, see Section \ref{Ssec.Upper.Bdry.T} for the precise statements. 

The above upper boundary estimate is sharp when $m>m_s$; indeed, it is possible to prove lower bounds of the same form. When $m\le m_s$, we can only match the spatial behaviour, which remains sharp for all $m\in (0,1)$, see Section \ref{Ssec.Lower.Bdry} for more details. The above estimates are the upper part of the so-called Global Harnack Principle, and of the Boundary Harnack inequalities discussed in Sections \ref{SSec.GHP} and \ref{SSec.BHI} respectively; see also Sections \ref{sec.GHP} and \ref{Ssec.BHI-T} for new results and alternative proofs.


\section{Mass conservation VS extinction in finite time}\label{Sec.FET}
The Dirichlet boundary condition does not allow for mass conservation for any $m>0$, in contrast with what happens for the Cauchy problem in the whole space. In the latter case, mass is preserved whenever $m\ge m_c:=\tfrac{N}{N-2}$, while it is not in the very fast diffusion range $m<m_c$ where solution can extinguish in finite time. To the best of our knowledge, the first estimates of the extinction time were obtained by Benilan-Crandall \cite{BCr2}, see also Vazquez's books \cite{VazBook,VazLN} for a more detailed exposition.

For the CDP, in the fast diffusion regime $m<1$, solutions always extinguish in finite time $T$; this, somehow, shows the super-diffusive character of the equation with respect to the Heat equation ($m=1$) or the Porous Medium Equation ($m>1$, for which finite speed of propagation holds). We define the Finite Extinction Time (FET) of a solution $u$ with initial datum $u_0$ as follows:
\[
T=T(u_0):=\inf\left\{t_1>0\;:\; u(t,x)=0\qquad\mbox{for a.e. $x\in \Omega$ and a.e. $t\ge t_1$}\right\}\,.
\]
Estimates on the extinction time, both from above and below in terms of $L^p$ norms of $u_0$, can be found in many papers, among which \cite{AnSm05, AnSm14, BV-ADV,BS2019,DKV1991, K2, JX2021} and in the monograph \cite{VazLN}. To the best of our knowledge, global lower bound (i.e., in terms of a norm on the whole $\Omega$) were firstly proven in \cite{BII2022}. Summing up, it is known that if $u_0$ satisfies \eqref{Hyp.u0.Lp} then there exists $\c_1,\c_p>0$,  depending only on $p,m,N, |\Omega|$,  such that
\[
\c_1 \|u_0\|_{L^1_{\Phi_1}(\Omega)}^{1-m}\le T(u_0) \le \c_p \|u_0\|_{L^p(\Omega)}^{1-m}.
\]
These estimates are strictly related to the issue of sharp extinction rates for $\LL^p$ norms, that we discuss below; we refer to Section \ref{Sec.Rayleigh}, in particular Subsection \ref{Sssec.Extinction.Norms},  for more details.

\noindent\textbf{Related problems. }\textit{The limit $m\to0^+$. }As we have seen, there are two possible limiting equations, and solutions exist only for the signed fast diffusion $u_t=\Delta {\rm sign}(u)$, see \cite[Section 3]{BF2012} for more details on the case $N=1$, where also the extinction time is explicitly calculated. On the other hand, if we impose \textit{Neumann boundary conditions, } then the mass is preserved and solutions do not extinguish, see \cite{I2019,IPS2019} for the ultrafast diffusion case, i.e., when $m\le 0$.

\subsection{Extinction rates for $\LL^p$ norms} In the supercritical range $m>m_s$, the first rates of extinction for $\LL^p$ norms were shown by Berryman-Holland \cite{BH}, then extended to stronger norms by Kwong \cite{K2} and DiBenedetto-Kwong-Vespri \cite{DKV1991}. Indeed it can be shown that
\[
 \|u(t)\|_{\LL^{p}(\Omega)} \asymp \|u(t)\|_{\LL^{1+m}(\Omega)} \asymp (T-t)^{\frac{1}{1-m}}\qquad\mbox{for all $p\in [1+m,\infty]$.}
\]
The estimate for the $\LL^{1+m}$ norm is based on a differential inequality satisfied by a ``nonlinear Rayleigh quotients'' along the FDE flow, see Section \ref{Sec.Rayleigh}. The rates for the $\LL^{1+m}$ can be then extended to  other $\LL^p$ norm using the smoothing effect \eqref{Smoothing.Global}, see Section \ref{SSec.Linfty.ext.T}.

\noindent\textit{Some new results. }In Section \ref{Sssec.Extinction.Norms} we show several results (some of them new) on extinction rates for various norms, including $H^*$ and weighted $\LL^p$ norms. In particular, Lemma \ref{Lem.Lp.norms} deals with extinction rates of $\LL^p$ norms for all $m\in(0,1)$, and their optimality is discussed in Remark \ref{Sharp.Lp.extinction}.  As a consequence, in Sections \ref{Ssec.Upper.Bdry.T} and \ref{Ssec.Lower.Bdry} we show upper and lower boundary estimates up to the extinction time, that fairly combine into Global Harnack Inequalities in different forms, see Sections \ref{Ssec.BHI-T} and \ref{sec.GHP}.


\section{Local VS global Harnack inequalities}\label{Sec.Harnack}
The celebrated results of Moser \cite{Moser1964,Moser1967,Moser1971} showed that nonnegative local  weak solution to linear uniformly parabolic equations  $u_t=\nabla\cdot(A(t,x) \nabla u)$ with bounded measurable coefficients (i.e., with $0<\lambda_0|\xi|^2\le \sum_{i,j=1}^d A_{i,j}\xi_i\xi_j\le \lambda_1|\xi|^2$ for all $\xi\in \RR^N$) satisfy the following Harnack Inequalities (HI)
\begin{equation}\label{HI.linear}
\sup_{D^{-}_R(t_0,x_0)} u\le\overline{\mathsf h}\,\inf_{D^{+}_R(t_0,x_0)}  u\,.
\end{equation}
 Here, the standard parabolic cylinders have the form
\[
\begin{split}
& D_R^+(t_0,x_0):=(t_0+\tfrac34\,R^2,t_0+R^2)\times B_{R/2}(x_0)\,,\\
& D_R^-(t_0,x_0):=\left(t_0-\tfrac34\,R^2,t_0-\tfrac14\,R^2\right)\times B_{R/2}(x_0)\,,
\end{split}
\]
Since the infimum is taken at later times, they are usually called Forward HI.
The constant $\overline{\mathsf h}$ can be explicitly expressed in the form $\overline{\mathsf h}:=\mathsf h^{\lambda_1+1/\lambda_0}$ where $\mathsf h$ only depends on $N$. Notice that in the current case (i.e., for the Laplacian) we have  $\lambda_0=\lambda_1=1$. See also \cite[Chapter 3]{BDNS2022} where a constructive proof of the above HI is given and the constant $\mathsf h$ is explicitly calculated.

When dealing with solutions to the homogeneous Dirichlet problem on a bounded domain $\Omega\subset \RR^N$, these results have been improved by Fabes-Garofalo-Salsa in \cite{FGS1986} to HI of elliptic type (i.e., supremum and infimum can be taken at the same time) and of backward type (i.e., infimum is taken at a previous time), namely, for all $|h|\le \delta R^2$ and all $B_R(x_0)\subset \Omega$
\[
\sup_{x\in B_R(x_0)} u(t,x)\le \overline{\mathsf h }\,\inf_{x\in B_R(x_0)} u(t\pm h,x)\,.
\]
where $\overline{\mathsf h}$ depends on $N, \lambda_0,\lambda_1$ and on the distance $\delta = \dist(B_R(x_0),\partial\Omega)$, and blows up when $\delta\to 0$. In \cite{FGS1986} there is also a version of backward/elliptic/forward boundary Harnack inequalities, i.e., the corresponding forms of HIs valid close to the boundary $\partial\Omega$, see also \cite{FSY1999}.

Note that HI of backward/elliptic type are not true for nonnegative solutions to the Cauchy problem on the whole space, the couterexample being the Gaussian (fundamental solution), and it was quite surprising to see that Dirichlet boundary conditions allow for such stronger inequalities.

\subsection{Local Harnack inequalities for FDE}\label{Ssec.Harnack.loc}
In general, the above HI \eqref{HI.linear} do not hold when $m\neq 1$, at least not on standard parabolic cylinders of the form $D_R^{\pm}(t_0,x_0)$. This is due to the singular character of the fast diffusion equation, which may possibly cause extinction in finite time. Note that local solutions need not to extinguish (as they can be ``pieces'' of solution to a Neumann problem, for instance), and another important factor enters the game: the scaling properties of the equation, roughly speaking $R^2\sim u^{m-1}T$. This means that the singular (in the present case $m\in (0,1)$) or degenerate (when $m>1$) character    of the equation induces a natural change of geometry in the parabolic cylinders, whose size needs to be intrinsically related to the size of the solution itself. For this reason they have been called ``intrinsic cylinders'' by DiBenedetto. Indeed, it can be shown that the size of intrinsic cylinders depends implicitly on an averaging value of the solution on the same cylinder, see the monographs \cite{DiBook,DGVbook, UrBook08}. Note that the intrinsic cylinders a priori may collapse or blow up when $u$ is not bounded or bounded away from zero. It turns out that these \textit{intrinsic cylinders are the natural domains of Harnack inequalities for degenerate and singular parabolic equations}:
\[
I_R(t_0,x_0)= \left(t_0-c\,u(t_0,x_0)^{1-m}R^2,
t_0+c\,u(t_0,x_0)^{1-m}R^2\right)\times B_R(x_0).
\]
In such domains, Intrinsic Harnack Inequalities (IHI) take a simpler form, closer to the  linear case:\\
\textit{There exist positive constants $\overline{c}$ and $\overline{\delta}$ depending only on $m,N$, such that for all $(t_0,x_0)\in Q=(0,T)\times\Omega$ and all cylinders of the type $I_{8R}\subset Q$, we have
\begin{equation}\label{IHI.DGV}
\overline{c}\,u(t_0,x_0)\le \inf_{x\in B_R(x_0)}u(t,x)
\end{equation}
for all times $t_0-\overline{\delta}\,u(t_0,x_0)^{1-m}\,R^2<t<t_0+\overline{\delta}\,u(t_0,x_0)^{1-m}\,R^2$. The constants
$\overline{\delta}$ and $\overline{c}$ tend to zero as $m\to 1$ or as $m\to m_c$\,.}

Notice that the above IHI are of backward/elliptic/forward type, and they have been first proven by DiBenedetto-Gianazza-Vespri in \cite{DGV2010}\footnote{The result by DiBenedetto-Gianazza-Vespri in \cite{DGV2010} is for local solutions, while for solutions to the \ref{CDP} forward IHI were first proven by DiBenedetto-Kwong \cite{DK1992}.}, where they also provide a counterexample to the validity of the above IHI in the subcritical range $m\le m_c$. 
The authors left open the intriguing question of which form, if any, Harnack inequalities would have taken in the very fast diffusion range. The answer was given by V\'azquez and the first author in \cite{BV-ADV}, where they showed the following:\\
\it Let $ u $ be a nonnegative local strong solution to $u_t=\Delta u^m$ on $(0,T)\times\Omega$ and let $0\le u_0\in \LL^p_{\rm loc}(\Omega)$, with $ p \geq 1 $ if $ m \in (m_c, 1) $ and $ p > p_c $ if $ m \in (0, m_c] $.  Let $B_{8R}(x_0)\subseteq\Omega$, $t_0\in [0,T)$, and
\begin{equation}\label{def.t*}
t_*(t_0)=\kappa_*\,R^{2-N(1-m)}\|u(t_0)\|^{1-m}_{\LL^1(B_R(x_0))}.
\end{equation}
Then, for any $\varepsilon\in (0, 1)$, there exists $\ka_3>0$ such that for any $t,t\pm\theta\in [t_0+\varepsilon t_*(t_0), t_0+t_*(t_0)]\cap (0,T)$
\begin{equation}\label{harnack.easy.thm}
\sup_{x \in B_R(x_0)} u(t, x) \leq \ka_3   \inf_{x \in B_R(x_0)}u(t \pm \theta, x).
\end{equation}
The constants $\kappa_*,\ka_3>0$ always depend on $N,m$ and have an explicit form; $\ka_3$ may also depend on $R, x_0$ and $\varepsilon$, and, when $0<m\le m_c$, it depends on the quotient $H_p(u_0, x_0, 2R)$ defined as
\begin{equation}\label{def.Hp}
    H_p\left( f, x_0, R \right):=
    \left[\frac{|B_R(x_0)|\left(\int_{B_R(x_0)}{f^p\dx}\right)^\frac{1}{p}}
            {|B_R(x_0)|^{\frac{1}{p}}\int_{B_R(x_0)}{f  \dx}} \right]^{2p \vartheta_p}\,.
\end{equation}
\rm Note that when we can take  $p=1$ (i.e., when $m>m_c$), $H_p$ simplifies to a constant, so that $\ka_3$ only depends on $m,N$ and we recover \eqref{IHI.DGV} with an explicit constant and with the size of intrinsic cylinders explicitly measured by $t_*$, i.e., measured in terms of the local mass of the solution at a previous time. Notice that $t_*$ can also be interpreted as (an estimate of) the \textit{minimal life time}, i.e., the time for which a local nonnegative solution will be non-trivial (when starting with nonzero mass). Indeed, $t_*$ also provides a lower bound for the extinction time for solutions to the   Cauchy-Dirichlet problem,    but it is still a local information. Actually, global bounds on the extinction time are known for solutions to the   Cauchy-Dirichlet problem,   as mentioned in Section \ref{Sec.FET}, see also Section \ref{Sec.Rayleigh}. The above result has been fist proven in \cite[Theorem 3.1]{BV-ADV}, using a Moser iteration (for the upper bounds) and an Aleksandrov moving plane argument combined with a ``Flux Lemma'' (for the lower bounds). After that, DiBenedetto-Gianazza-Vespri extended the result to the ``measurable coefficient case'' through a De Giorgi technique, see the monograph \cite{DGVbook} and other generalizations appear also in Fornaro-Henriques-Vespri \cite{FHV21}. The above IHI \eqref{harnack.easy.thm} was extended to solutions of the FDE with Caffarelli-Kohn-Nirenberg weights by Simonov and the first author in \cite[Theorem 1.6]{BS2019}, with a new constructive proof for the lower bounds. In \cite{BV-ADV,BS2019} the form of $\ka_3$ is explicitly given:
\[
\ka_3
\asymp \varepsilon^{-\frac{2 p \vartheta_p}{1-m}} H_p(u_0, x_0, 2R) \, \widetilde{H}_p^{\frac{c\,\widetilde{H}_p^{1/2}}{m(1-m)}}\,,
\]
where $\widetilde{H}_p(f, x_0, R):= 1+H_p(f, x_0, R)^{1-m}$
with $H_p$ as in \eqref{def.Hp}.
This expression allows to quantify the H\"older continuity exponent, as we shall see below. We refer also to the memoir \cite[Chapters 3 and 4]{BDNS2022} by Dolbeault-Nazaret-Simonov and the first author, where constructive proofs and explicit constants can be found.

\subsection{Global Harnack Principle. }\label{SSec.GHP}Local Harnack inequalities are ``interior estimates'' and hold independently of the boundary data, hence we can expect better results for solutions to a Cauchy-Dirichlet problem, as it happens in the linear case. Also it is possible to have estimates that capture the sharp boundary behaviour, possibly up to the extinction time. The first ``global Harnack estimate'', or Global Harnack Principle (GHP), was proven in the range $m>m_s$ by Kwong \cite{K1,K3} and  DiBenedetto-Kwong-Vespri \cite{DKV1991}: \it for all $t_0>0$ there exist constants $\ka[u_0]\ge \kb[u_0]>0$, depending on $N,m,\Omega,t_0$ and $u_0$, such that
\begin{equation}\label{GHP.DKV}
\kb[u_0](T-t)^{\frac{m}{1-m}} \le \frac{u^m(t,x)}{\Phi_1(x) }\le
\ka[u_0] \,
(T-t)^{\frac{m}{1-m}} \,, \qquad\mbox{ for all $x\in \overline{\Omega}$ and $t\in[t_0,T]$\,. }
\end{equation}\rm
This is one of the key points where the ``Sobolev exponent'' $m_s$ plays a role. The constants $\ka[u_0],\kb[u_0]$ depend on $u_0$ through $\Q[u_0]$, i.e., through $\|\nabla u_0^m\|_{\LL^2(\Omega)}$ and $\|u_0\|_{1+m}$, and blow up when $t_0\to 0^+$, or when $m\to 1$ or $m\to m_s$.

In this paper we provide a different proof of the GHP, which is not based on barriers as in \cite{DKV1991}, and which allows to eliminate the dependence on $\|\nabla u_0^m\|_{\LL^2(\Omega)}$ in the constants. Moreover, our technique allows to extend the GHP to all $m\in (0,1)$ as follows:
\begin{equation}\label{GHP.BF}
 \kb[u_0,t]   (T-t)^{\frac{m}{1-m} +   \frac{p-2m+1}{1-m}2\vartheta_p[p-(1+m)]_+}   
\le \frac{u^m(t,x)}{\Phi_1(x) }\le
\ka[u_0,t] \,
(T-t)^{\frac{m}{1-m}-  \frac{2\vartheta_p }{1-m}[p-(1+m)]_+},
\end{equation}
  where $\kb,\ka>0$ depend on $m,p,N,\Omega$ and $u_0$, and have explicit expressions, see Theorem \ref{thm.GHP.I}. When $t\ge \frac{2}{3}$, $\ka[u_0,t]$ and $\kb[u_0,t]$ do  not depend on $t$.

Notice that, when $m>m_s$, one can choose $p=1+m$ in the above estimate and recover an improved version of the GHP \eqref{GHP.DKV} of \cite{DKV1991}: our proof allows us to deal with a larger class of initial data and to obtain an explicit constant that depends only on $\|u_0\|_{H^{-1}(\Omega)}$ and $\|u_0\|_{1+m}$. We refer to Section \ref{sec.GHP} for more details, in particular to Theorem \ref{thm.GHP.I} and Remark \ref{Rem.GHP}.

When $m\in (m_{c,1},1)$ we can even eliminate the dependence on $u_0$ in the constants of the GHP \eqref{GHP.DKV}
\begin{equation*}
\kb (T-t)^{\frac{m}{1-m}} \le \frac{u^m(t,x)}{\Phi_1(x) }\le
\ka  \,
(T-t)^{\frac{m}{1-m}} \,,
\end{equation*}
where $\ka,\kb>0$  only depend on $m,N,\Omega$, see Theorem \ref{thm.GHP.II}, one of our main results.

We complete the panorama with Theorem \ref{thm.NoGHP}, which shows that  when $m\in (0,m_s]$ it is  not possible to have a GHP with matching time powers, of the form \eqref{GHP.DKV}. As a consequence, our GHP \eqref{GHP.BF} is the only known global Harnack estimate valid in the whole range $m\in (0,1)$. We would like to point out that when $m\le m_s$ the power of time in the upper bound can be negative, hence it estimates a maximal blow-up rate for solutions, see Sections \ref{Ssec.Yamabe} and \ref{Ssec.Asymp.Subcrit.New} for a thorough discussion.

\subsection{Harnack Inequalities up to the Boundary (BHI) }\label{SSec.BHI} Boundary Harnack Inequalities (BHI) for solutions to singular parabolic equations are not straightforward to state, and are particularly relevant when the domain is non-smooth, i.e., when it has a Lipschitz or H\"older continuous boundary, a delicate issue that we have chosen not to discuss here. For more details, we refer to the papers of Kuusi-Mingione-Nystrom \cite{KMN2014}, Avelin-Gianazza-Salsa \cite{AGS2016} and references therein. In particular, in \cite{AGS2016} Carleson estimates and BHI are proven for solutions to \eqref{CDP} in the good fast diffusion range, namely when $m\in (m_c,1)$.

In this paper we show uniform estimates valid for all $m\in (0,1)$, which we still call BHI, because of their formal similarities and consequences:
\[
\sup_{x\in \Omega}\frac{u^m(t,x)}{\Phi_1(x) }
\le \frac{\overline{\mathsf h}[u_0,t]}{(T-t)^{\frac{2\vartheta_p }{1-m}[p+2(1-m)][p-(1+m)]_+}}
\inf_{x\in \Omega}\frac{u^m(t,x)}{\Phi_1(x) } \,.
\]
This is precisely proven in Theorem \ref{thm.BHI.1}.
The above   estimate   easily implies the more classical form of BHI:   for all $m\in (0,1)$,   given two bounded nonnegative solutions $u,v$ to the \eqref{CDP} with the same extinction time $T$, we have that
\[
\frac{u^m(t,x)}{v^m(t,x)}
\le  \frac{\overline{\mathsf h}[u_0,t ]^2}{(T-t)^{\frac{4\vartheta_p }{1-m}[p+2(1-m)][p-(1+m)]_+}}
 \frac{u^m(t,y)}{v^m(t,y)}\qquad\mbox{for all }x,y\in \overline{\Omega}\,. \,.
\]
The GHP or BHI are the key properties to prove higher regularity estimates, as we shall discuss next. Indeed, once solutions are positive and bounded they are smooth (in the interior). The regularity at the boundary is more delicate.


\section{Interior VS boundary regularity}\label{Ssec.Regularity}

\subsection{Interior regularity. }
\noindent\textbf{H\"older regularity. }Local bounded solutions have been proven to be continuous up to the boundary by Sacks \cite{Sa1983} and DiBenedetto \cite{D1983}; indeed, in the latter paper, a logarithmic modulus of continuity was obtained for the first time. Later on, bounded solutions were proven to be H\"older continuous in intrinsic cylinders by Chen-DiBenedetto \cite{CD1988,CD1992}, using a nonlinear adaptation of De Giorgi method. Analogously to what happens in the linear case, also in the nonlinear setting IHI imply H\"older continuity (in intrinsic cylinders), see DiBenedetto-Kwong \cite{DK1992} and the monographs by DiBenedetto \cite{DiBook}, DiBenedetto-Gianazza-Vespri \cite{DGVbook} and Urbano \cite{UrBook08}. More recently, Simonov and the first author \cite{BS2019, BDNS2022} gave a constructive proof of the $\LL^\infty-C^\alpha$ estimates through a nonlinear adaptation of Moser's method. Theorem 1.8 of \cite{BS2019} adapted to our setting (avoiding the intrinsic cylinders notation) reads:\\ \it Let $Q_0=(t_0,t_0+t_*)\times B_{4R_0}(x_0)\subset [0,T]\times\Omega$, and $t_*$ as in \eqref{def.t*}. There exist  $\alpha\in (0,1)$ and  $\ka_\alpha'>0$ such that
\begin{equation*}
|u(t,x)-u(\tau,y)|
\le\frac{\ka_\alpha'\, \|u\|_{\LL^\infty(Q_0)}}{t_*^{\alpha/2}  R_0^\alpha} \left(|x-y| + \|u\|_{\LL^\infty(Q_0)}^{\frac{m-1}{2}}|t-\tau|^{\frac{1}{2}} \right)^\alpha\,
\end{equation*}
for all $t,\tau \in \left[t_0+\frac{5}{8}t_*, t_0+\frac{7}{8}t_*\right]\cap(0,T)$ and for all $x, y \in B_{R_0}(x_0)$.\rm

Note that both $\alpha $ and  $\ka_\alpha' $ are explicit and depend only on $m,N$ and $H_p=H_p(u(t_0),x_0,4R_0)$ defined in \eqref{def.Hp}. Indeed, it is shown in \cite{BS2019} that the exponent $\alpha$ depends on $H_p$ and $t_*$ in a quantitative way:
\[
   \alpha\sim {\rm exp}\left(- \frac{c_6}{t_*}H_p^{\frac{c_7(1-m)}{m}H_p^{(1-m)/2} }\right)\,,
\]
where   $c_6,c,7>0$   only depend on $N,m,p$. Having at disposal a IHI with a uniform constant, i.e., independent of the solution or on the initial datum, would allow to prove quantitative uniform H\"older regularity.

\noindent\textbf{$C^\infty$ regularity and more. }The unwritten principle of parabolic regularity states that when solutions are bounded and positive then they are smooth (as much as the operator allows). This happens in our case.
Indeed, when $m\in (m_s,1)$, DiBenedetto-Kwong-Vespri \cite{DKV1991} show that, as a consequence of the GHP, nonnegative   bounded weak solutions to the \eqref{CDP} even analytic in the interior, and up to the extinction time: we recall Theorem 1.1 of \cite{DKV1991}, adapted to our notation: \\\it Let $m\in (m_s,1)$ and let $u$ be a bounded solution to \eqref{CDP}. Then, for every $t\ge\varepsilon>0$, $k=0,1,2,\dots$ there exist $\ka_{k,\varepsilon}>0$ such that
\[
\left|D_x^k u(t,x)\right|\le \ka_{k,\varepsilon} \dist(x,\partial\Omega)^{1-k}(T-t)^{\frac{1}{1-m}},
\]
and
\[
\left|\partial_t^k u(t,x)\right|\le \ka_{k,\varepsilon} \dist(x,\partial\Omega)^{1-k\frac{1+m}{m}}(T-t)^{\frac{1}{1-m}-k},
\]
where $\ka_{k,\varepsilon}\sim k!$ depends on $m,N, \Omega$, $\|u_0\|_{\LL^{1+m}(\Omega)}, \|\nabla u_0^m\|_{\LL^2 (\Omega)}$ and $\varepsilon$. Note that $\ka_{k,\varepsilon}\to \infty$ when $\varepsilon\to 0^+$ and also when $m\to m_s$.\rm

Recently, Jin-Xiong \cite{JX2019,JX2022} have shown that solutions are smooth up to the boundary, again assuming the validity of the GHP. More precisely they have proven that \it  solutions to \eqref{CDP} that satisfy the GHP are classical and smooth up to the boundary, namely $u^m(\,\cdot\,,x)\in C^\infty(0,T)$ for all $x\in \overline{\Omega}$ and $\partial_t^ku^m(t,\,\cdot\,)\in C^{2+\frac{1}{m}}(\overline{\Omega})$ for all $t\in (0,T)$. Moreover, $u^m\in C^{\infty}((0,T)\times \overline{\Omega})$ when $1/m$ is an integer.    \\
\rm See Theorem 1.1 of  \cite{JX2019} for the case $m\in [m_s,1)$, then extended to all $m\in (0,1)$ in Theorem 1.1 of  \cite{JX2022}. As the authors observe, the $C^{2+\frac{1}{m}}(\overline{\Omega})$ regularity is optimal, in view of solutions of separation of variables whose elliptic part is known to possess such optimal regularity.

\subsection{Boundary regularity. }The GHP tells us that $u^m\asymp \dist(\,\cdot\,,\partial\Omega)$. The question of the sharp boundary regularity consists in establishing whether or not we have
\[
\frac{u^m(t,\,\cdot\,)}{\dist(\,\cdot\,,\partial\Omega)}\in C^{\alpha}(\overline{\Omega})\quad\mbox{for some }\alpha\in (0,1).
\]
The answer has been recently given by Jin-Xiong \cite{JX2019,JX2021,JX2022}: in the first paper \cite{JX2019} they establish how the GHP implies sharp boundary regularity estimates in the case $m\in [m_s,1)$ locally in time. Later, in \cite{JX2021}, they extend the results up to the extinction time $T$. Finally, they prove how the GHP implies the optimal boundary regularity for all $m\in (0,1)$ and all times $t<T$. We summarize here the main results of \cite{JX2019,JX2021,JX2022}, adapted to our notation: \\
\it Let $m\in (0,1)$ and let $u$ be a bounded solution to \eqref{CDP} that satisfies the GHP. Then, for every $\varepsilon>0$   and $j,k=0,1,2,\dots$ there exist $\ka_{k,\varepsilon}>0$ such that
for all $t\in (\varepsilon, T-\varepsilon)$
\[
\left\|\frac{ \partial_t^j u^m(t,\,\cdot\,) }{\dist(\,\cdot\,,\partial\Omega)}\right\|_{\LL^\infty(\overline{\Omega})}
+ \left\|  \partial_t^j u^m(t,\,\cdot\,) \right\|_{C^{2+\frac{1}{m}}(\overline{\Omega})}
\le \ka_{k,\varepsilon}\,.
\]
Moreover, when $1/m$ is an integer, we have that for all $k=0,1,2,\dots$ there exist $\ka_{j,k,\varepsilon}>0$ such that
\[
\left\|\frac{ \partial_t^j u^m(t,\,\cdot\,) }{\dist(\,\cdot\,,\partial\Omega)}\right\|_{\LL^\infty(\overline{\Omega})}
+\left\|D_x^k \partial_t^j u^m(t,\,\cdot\,)\right\|_{\LL^\infty(\overline{\Omega})}\le \ka_{j,k,\varepsilon} ,
\]
where $\ka_{k,\varepsilon}, \ka_{j,k,\varepsilon}>0$ depend on $k,m,N, \Omega$, $\|u_0\|_{\LL^{1+m}(\Omega)}, \|\nabla u_0^m\|_{\LL^2 (\Omega)}$ and $\varepsilon$. Also, $\ka_{j,k,\varepsilon}\xrightarrow[]{\varepsilon\to 0^+} \infty$. \\
Finally, when $m\in (m_s,1)$ the above estimates can be extended up to $t=T$, in which case the above estimates have an extra term on the right-hand side, namely a multiplying factor $(T-t)^{\frac{m}{1-m}-j}$.
\rm

\noindent The proof is contained in \cite{JX2019,JX2022} and uses many tools (in a suitable weighted setting) like De Giorgi-Nash-Moser iterations, Campanato spaces, Schauder estimates, and other ingenious bootstrap arguments. The main novelty consists in a careful  analysis of the ``curvature term'' $w=u_t/u$. Since $w$ satisfies a ``nice'' weighted parabolic equation (see also Section \ref{SSec.new.proof.BC}), it is possible to deduce energy estimates involving weighted $\LL^q$ norms (with $u$ as a time dependent weight) that iterated provide the basic regularity estimates, that finally can be bootstrapped up to $C^{2+\frac{1}{m}}$ (or $C^\infty$).

As for interior regularity, the key assumption is again the validity of a GHP of the form \eqref{GHP.BF}. Indeed, we are not aware of a proof of the GHP in the subcritical range $m\le m_s$ except the one contained in this paper. Hence, the results of this paper ensure the sharp boundary regularity for all $m\in (0,1)$, and also for a larger class of nonnegative solutions: we are able to weaken the known assumptions on the initial datum such that GHP holds for the corresponding solution. Recall that we allow $0\le u_0\in \LL^p_{\Phi_1}(\Omega)$, which does not automatically extend to a $\LL^p_{\rm loc}(\RR^N)$ function.


\section{Global boundary estimates via duality and Green functions}\label{Sec.GHP.new}
This section contains new result about upper and lower estimates for (weak dual) solutions to the    \eqref{CDP} for the FDE     that are global both in space (up to $\partial\Omega$) and in time (up to the extinction time). The GHP of Theorems \ref{thm.GHP.I} and \ref{thm.GHP.II} in Section \ref{sec.GHP} is new when $m\le m_s$, while when $m\in (m_s,1)$ it extends the GHP of \cite{DKV1991} to a bigger class of solutions, as discussed in Section \ref{SSec.GHP}, cf. also Remark \ref{Rem.GHP}. We complete the panorama with Theorem \ref{thm.NoGHP.ineq}, showing that \textit{in the subcritical $m<m_s$ range a GHP uniform up to the extinction time (as the one valid for $m>m_s$) is not possible in general}, at least not in star-shaped domains. Hence our form of GHP  seems to be optimal.
\subsection{Dual equation and an almost representation formula}\label{Ssec2.1}
It is well-known\footnote{Sharp estimates for the Green functions are often derived from optimal heat kernel bounds. The sharp lower bounds were proven for the first time by Zhang \cite{Z2002},  while the sharp upper bounds were known before, see Davies' book \cite{DaviesBook}\,.} that the Green function of the Dirichlet Laplacian on a smooth domain $\Omega\subset \RR^N$ satisfies the following two-sided (sharp) estimates
\begin{equation}\label{G4}\tag{G4}
 \G(x,y)\asymp \frac{c_1}{|x-y|^{N-2}}\left( \frac{\dist(x,\partial\Omega)}{|x-y|}\wedge 1\right)\left( \frac{\dist(y,\partial\Omega)}{|x-y|}\wedge 1\right)
\end{equation}
with upper and lower constants $0<c_0<c_1$ which depend only on $N,\Omega$. The first eigenfunction of the Laplacian is $C^\infty$ in the interior and satisfies $\Phi_1\asymp \dist(\cdot,\partial\Omega)$, hence it is convenient to use it as a smooth extension of the distance function. We also assume that $\|\Phi_1\|_{\LL^2(\Omega)}=1$. Recall that $\G$ is the kernel of the inverse of the Dirichlet Laplacian, namely
\[
\AI f(x)=\int_\Omega f (y) \G(x,y)\dy
\]
Applying $\AI=(-\Delta)^{-1}$ to both sides of the FDE $u_t=\Delta u^m$ we obtain the dual equation
\[
\AI u_t = - u^m
\]
This clarifies the advantage of using weak dual solutions of Definition \ref{Def.WDS}, which are nothing but the weak formulation of the dual problem. Formally integrating the above inequality in time, over $[t_0,t_1]$ we obtain that WDS satisfy 
\[
\int_{t_0}^{t_1}u^m(t,x) \dt = \int_\Omega \big[u(t_0,y)-u(t_1,y)\big]\G(x,y)\dy\,.
\]
The same can be obtained by using the (a priori not admissible\footnote{  We recall that the test functions $\psi$ in the definition of WDS need to be in the space $\psi/\Phi_1\in C^1_c((0,T):L^{\infty}(\Omega))$, and our candidate $\chi_{[t_0,t_1]}(t)\G(x_0,x)$ does not belong to that space.})  test function $\chi_{[t_0,t_1]}(t)\G(x_0,x)$ in the weak formulation \eqref{WDS}. Next we can estimate the time integral using the time monotonicity provided by the Benilan-Crandall estimates \eqref{Be-Cr}: for all $0<t_0\le t\le t_1$ and a.e. $x\in \Omega$ we have
\[
\left(\frac{t}{t_1}\right)^{\frac{1}{1-m}}u(t_1,x)\le u(t)\le \left(\frac{t}{t_0}\right)^{\frac{1}{1-m}} u(t_0,x).
\]
This easily implies that
\[
 \,\frac{t_1^{\frac{1}{1-m}}-t_0^{\frac{1}{1-m}}}{t_1^{\frac{m}{1-m}}}  u^{m}(t_1,x)\le \frac{1}{1-m}\int_{t_0}^{t_1}u^{m}(t,x) \dt \le \,\frac{t_1^{\frac{1}{1-m}}-t_0^{\frac{1}{1-m}}}{t_0^{\frac{m}{1-m}}}u^{m}(t_0,x)\,.
\]
This provides a formal proof of what we call ``almost representation formula'', given in the following:
\begin{lem}[Fundamental pointwise estimates]\label{Prop.PE}
  Let $m\in (0,1)$ and let $u$ be a  WDS to \eqref{CDP}   in the sense of Definition \ref{Def.WDS}. Then, for all   $0\le  t_0 < t_1$  and a.e. $x\in \Omega$ we have
\begin{equation}\label{Prop.PE.ineq}
\frac{u^m(t_1,x)}{t_1^{\frac{m}{1-m}}} \le \frac{1}{1-m}\int_\Omega  \frac{u(t_0,y)-u(t_1,y)}{ t_1^{\frac{1}{1-m}}-t_0^{\frac{1}{1-m}} } \G(x,y)\dy \le \frac{u^m(t_0,x)}{t_0^{\frac{m}{1-m}}}\,.
\end{equation}
\end{lem}
The above Lemma is a particular case of Lemma 3.4 of \cite{BII2022} when $s=1$. The proof is obtained by approximating in the weak dual formulation \eqref{WDS} the test function $\chi_{[t_0,t_1]}(t)\G(x_0,x)$ by means of   admissible ones,   thus making  rigorous the heuristic proof given above.   Note that when $t_0=0$, only the first inequality of \eqref{Prop.PE.ineq} is really meaningful,  since the latter term becomes $+\infty$.

\subsection{Upper Boundary estimates}
We start by recalling a recent result of \cite{BII2022}, proven by means of the so-called Green function method used in the case $m>1$ to obtain the sharp results by Ros-Oton, V\'azquez, and the two authors \cite{BV-ARMA,BV-NA, BFR2017, BFV2018}. Indeed, boundary estimates easily follow by the standard smoothing effect combined with the almost representation formula, as we shall explain next.
\begin{thm}[Upper boundary estimates]\label{thm.GHP.up}
Let $u$ be a   WDS to \eqref{CDP}    corresponding to $0\le u_0\in \LL^1_{\Phi_1}(\Omega)$.  \\
(i) Let $u_0\in \LL^p(\Omega)$ with $p\ge 1$ if $m\in (m_c,1)$ and $p>p_c$ if $m\in (0,m_c]$, where $m_c=\frac{N-2}{N}$ and $p_c=\tfrac{N(1-m)}{2}$. There exists $\ka>0$ depending only on $m,p,N$ and $\Omega$ such that for all $0\le t_0\le t_1$\vspace{-2mm}
\begin{equation}\label{Smoothing.p.bdry}\begin{split}
\left\|\frac{u^m(t)}{\Phi_1}\right\|_{\LL^\infty(\Omega)}
& \le  \ka
\frac{\|u(t_0)\|_{L^p(\Omega)}^{2p\vartheta_p}}{(t-t_0)^{N\vartheta_p+1}}
\qquad\mbox{with}\qquad \vartheta_p=\frac{1}{2p-N(1-m)}\,.\\
\end{split}
\end{equation}
(ii) Let $u_0\in \LL^p_{\Phi_1}(\Omega)$ with $p\ge 1$ if $m\in (m_{c,1},1)$ and $p>p_{c,1}$ if $m\in (0,m_{c,1}]$, where $m_{c,1}=\frac{N-1}{N}$ and $p_{c,1}=  N(1-m)$. There exists $\ka>0$ depending only on $m,p,N$ and $\Omega$ such that for all $0\le t_0\le t_1$\vspace{-2mm}
\begin{equation}\label{Smoothing.p.bdry.Phi1}\begin{split}
\left\|\frac{u^m(t)}{\Phi_1}\right\|_{\LL^\infty(\Omega)}
& \le  \ka
\frac{\|u(t_0)\|_{L^p_{\Phi_1}(\Omega)}^{p\vartheta_{p,1}}}{(t-t_0)^{N\vartheta_{p,1}+1}}
\qquad\mbox{with}\qquad \vartheta_{p,1}=\frac{1}{p-N(1-m)}\,.\\
\end{split}
\end{equation}
\end{thm}
\noindent {\bf Proof.~} Recall the fundamental upper bounds \eqref{Prop.PE.ineq}: for all $0\le t_0\le t_1$
\begin{align} \label{Prop.PE.ineq1}
u^m(t_1,x)&\le \frac{t_1^{\frac{m}{1-m}}}{1-m}\int_\Omega  \frac{u(t_0,y)-u(t_1,y)}{ t_1^{\frac{1}{1-m}}-t_0^{\frac{1}{1-m}} } \G(x,y)\dy
\le \frac{1}{1-m}\frac{t_1^{\frac{m}{1-m}}}{ t_1^{\frac{1}{1-m}}-t_0^{\frac{1}{1-m}} } \|u(t_0)\|_{\LL^\infty(\Omega)} \|\G(x,\cdot)\|_{\LL^1(\Omega)}\nonumber\\
&\le c_1 \frac{\Phi_1(x)}{1-m}\frac{t_1^{\frac{m}{1-m}}}{ t_1^{\frac{1}{1-m}}-t_0^{\frac{1}{1-m}} }\|u(t_0)\|_{\LL^\infty(\Omega)}\,,
\end{align}
where we have used the Green function estimate \eqref{G4}, which implies that $\|\G(x,\cdot)\|_{\LL^1(\Omega)}\dy\le c_1\Phi_1(x)$, where $c_1>0$ depends only on $N$ and $\Omega$, see for instance Lemma 4.2 of \cite{BFV2018} (taking $s=\gamma=1$ there).

Inequality \eqref{Smoothing.p.bdry} for $t_0=0$ follows by combining \eqref{Prop.PE.ineq1} with the smoothing effects \eqref{Smoothing.Global} applies with $t_0=t_1/2$. Then we can extend the smoothing to all $t_0>0$ by time-shift invariance.
The proof of inequality \eqref{Smoothing.p.bdry.Phi1} is similar, just by using the smoothing effect \eqref{Smoothing.Global.Weighted}\,.\qed

\subsection{Energies, nonlinear Rayleigh quotients, and extinction rates}\label{Sec.Rayleigh}
We recall a monotonicity property of the nonlinear Rayleigh quotient, firstly shown by Berryman-Holland in their pioneering paper \cite{BH}. We shall introduce a \textit{dual nonlinear Rayleigh quotient} which also decreases along the flow. Finally, we shall draw consequences in terms of extinction rates of $\LL^p$ norms. For the well-known results we will just explain the main ideas, referring to the papers of Dibenedetto-Kwong-Vespri \cite{DKV1991} and Jin-Xiong \cite{JX2019} for a rigorous justification of the classical energy (in)equalities. The dual quotient has been introduced  and its monotonicity along the flow rigorously justified by Ibarrondo-Ispizua and the first author
\cite[Section 4.2]{BII2022}, in the framework of nonlocal fast diffusion equations. Consider the nonlinear Rayleigh quotient for all $m\in (0,1)$
\[
\Q[f]:=\frac{\|\nabla f^m\|_{\LL^2(\Omega)}^2}{\|f\|_{\LL^{1+m}(\Omega)}^{2m}}\,.
\]
It can been shown that $t\mapsto \Q[u(t)]$ it is non-increasing along the flow: let $u$ be a WDS of \eqref{CDP}, then
\[
\frac{\rd}{\dt}\Q[u(t)]\le 0 \qquad\mbox{hence}\qquad \Q[u(t)]\le \Q[u_0]\,.
\]
This has consequences for the behaviour of the $\LL^{1+m}$ energy:
\begin{equation}\label{deriv.1+m}
\frac{1}{1+m}\frac{\rd}{\dt}\|u(t)\|_{\LL^{1+m}(\Omega)}^{1+m}=- \|\nabla u^m(t)\|_{\LL^2(\Omega)}^2
= - \Q[u(t)] \|u(t)\|_{\LL^{1+m}(\Omega)}^{2m} \ge  - \Q[u_0] \|u(t)\|_{\LL^{1+m}(\Omega)}^{2m}\,.
\end{equation}
Integrating over $[t,T]$,   $T$ being the extinction time,   we get that for some $\c_m>0$ (depending only on $m\in (0,1)$)
\[
 \|u(t)\|_{\LL^{1+m}(\Omega)} \le \c_m \Q[u_0]^{\frac{1}{1-m}} (T-t)^{\frac{1}{1-m}}\,.
\]
This estimate automatically provides extinction rates for all $\LL^p$ norms with $p\le 1+m$. In the case $m>m_s$ this will be sufficient to have sharp extinction rate in $\LL^\infty$ norm, while when $m<m_s$ this will not be the case, see the new results of Section \ref{Ssec.Asymp.Subcrit.New}.
The monotonicity of the quotient follows from \eqref{deriv.1+m} and from the energy inequality:
\begin{equation}\label{deriv.H.norm}
\frac{1}{2}\frac{\rd}{\dt}\|\nabla u^m(t)\|_{\LL^2(\Omega)}^2
\le -m\int_\Omega \frac{(u^m\Delta u^m)^2}{u^{1+m}}\dx
\le   -m \frac{\|u^m\Delta u^m\|_{\LL^1(\Omega)}^2}{\|u\|_{\LL^{1+m}(\Omega)}^{1+m}}
=  -m \frac{\|\nabla u^m(t)\|_{\LL^2(\Omega)}^4 }{\|u\|_{\LL^{1+m}(\Omega)}^{1+m}}\,.
\end{equation}
The above proof,   that uses the Cauchy-Schwarz inequality in the form $\int_\Omega \tfrac{f^2}{g}\dx\ge\tfrac{\|f\|_{\LL^1(\Omega)}^2}{\|g\|_{\LL^1(\Omega)}}$,     is due to Berryman and Holland \cite{BH}.
The first inequality of \eqref{deriv.H.norm} is indeed an identity when $m>m_s$, as recently shown by Jin-Xiong \cite{JX2019}.

On the other hand, we observe that in many cases we have that $\Q[u_0]=+\infty$, since $u_0$ may be just in $\LL^{1}_{\Phi_1}$ or $H^{-1}$. To overcome this issue, one possibility is to use the following estimates: for all    $0\le t_0< t_1 < t_2$
\begin{equation}\label{energy.H1.Hstar}
\|\nabla u^m(t_2)\|_{\LL^2(\Omega)}^2 \le \frac{1}{2m} \frac{\|u(t_1)\|_{\LL^{1+m}(\Omega)}^{1+m}}{t_2-t_1}
\le \frac{1}{2m(1+m)} \frac{\|u(t_0)\|_{H^{-1}(\Omega)}^{2}}{(t_2-t_1)(t_1-t_0)}\,,
\end{equation}
see   Lemma 3.6    of \cite{BII2022} for a rigorous proof.

Let $\mathcal{L}=-\Delta$ and recall that the kernel of $\AI=(-\Delta)^{-1}$ is given by the Green function $\G(\,\cdot\,,\,\cdot\,)$ so that
\[
\|f\|_{H^{-1}(\Omega)}^2=\int_\Omega f(x)\AI f(x)\dx = \int_\Omega \left|\AIM f(x)\right|^2 \dx\,.
\]
To be able to deal with less regular data, it is convenient to consider the ``dual Rayleigh quotients''
\[
\Q^*[f]=\frac{\|f\|_{1+m}^{1+m}}{\|f\|_{H^{-1}(\Omega)}^{1+m}}\,.
\]
This quotient perfectly fit in the $H^{-1}$ setting introduced by Brezis\footnote{Recall that the FDE is a gradient flow in $H^{-1}$ of the energy given by the $\LL^{1+m}$-norm, as first noticed by Brezis in \cite{Brezis1}, see also the excellent notes by Ambrosio-Brue-Semola \cite{Ambrosio-Notes}. We also refer to \cite{BII2022} where $\Q^*$ has been introduced for the first time. There, the Brezis-Komura gradient flow approach to nonlinear diffusions in Hilbert spaces is applied to nonlocal fast diffusion equations, including the case at hand. Also WDS are constructed and complete proofs can be found.} and also turns out to be monotone along the flow, namely
\[
\frac{\rd}{\dt}\Q^*[u(t)]\le 0 \qquad\mbox{hence}\qquad \Q^*[u(t)]\le \Q^*[u_0]\,.
\]
This easily follows from the fact that
\[
\frac{\rd}{\dt}\|u(t)\|_{H^{-1}(\Omega)}^2=-2\|u(t)\|_{\LL^{1+m}(\Omega)}^{1+m}
\]
and
\begin{equation}\label{deriv.1+m.Hstar}
\begin{split}
\frac{\rd}{\dt}\frac{\|u(t)\|_{\LL^{1+m}(\Omega)}^{1+m}}{1+m}=-\int_\Omega \frac{(\AM u^m)^2(\AIM u)^2}{(\AIM u)^2}
\le - \frac{\|(\AM u^m)(\AIM u)\|_{\LL^1(\Omega)}^2}{\|(\AIM u)^2\|_{\LL^1(\Omega)}}
\le - \frac{\|u(t)\|_{\LL^{1+m}(\Omega)}^{2(1+m)}}{\|u(t)\|_{H^{-1}(\Omega)}^2}\,,
\end{split}
\end{equation}
where we have used again   the Cauchy-Schwarz inequality as above.   Note that integrating the above inequality and using that $\|u(t)\|_{H^{-1}(\Omega)}\le \|u(t_0)\|_{H^{-1}(\Omega)}$ proves the last inequality of \eqref{energy.H1.Hstar}, while the first one can be proven analogously to \eqref{deriv.H.norm}.

As a first consequence of the monotonicity of $\Q^*$, we can deduce an extinction rate for the $H^{-1}$ norm:
\begin{equation*}
\frac{\rd}{\dt}\|u(t)\|_{H^{-1}(\Omega)}^2=-2\|u(t)\|_{\LL^{1+m}(\Omega)}^{1+m}=-2\Q^*[u(t)] \|u(t)\|_{H^{-1}(\Omega)}^{1+m}
\ge -2\Q^*[u_0] \|u(t)\|_{H^{-1}(\Omega)}^{1+m}\,,
\end{equation*}
which integrated over $[t,T]$,   $T$ being the extinction time,    gives for all $t\in [0,T]$
\[
\|u(t)\|_{H^{-1}(\Omega)}^{1-m}\le 2\Q^*[u_0] (T-t)\,.
\]
We immediately deduce an extinction rate for the $\LL^{1+m}$ norm which does not depend on $\|\nabla u_0^m\|_{\LL^2(\Omega)}$. Indeed, inequality \eqref{energy.H1.Hstar} gives for all $t_1\ge t_0\ge 0$:
\begin{equation*}
\begin{split}
(1+m) \|u(t_1 )\|_{\LL^{1+m}(\Omega)}^{1+m}
\le \frac{\|u(t_0)\|_{H^{-1}(\Omega)}^{2}}{t_1-t_0}
\le \left(2\Q^*[u_0]\right)^{\frac{2}{1-m}}\frac{(T-t_0)^{\frac{2}{1-m}}}{t_1-t_0}\,.
\end{split}
\end{equation*}
Choose $t_1\ge T/3$ and $t_0=t_1-\tfrac{T-t_1}{2}$ to get (notice that $t_0\ge 0$ since $t_1\ge T/3$)
\begin{equation}\label{deriv.1+m.Hstar}
\begin{split}
\|u(t_1 )\|_{\LL^{1+m}(\Omega)}^{1+m}
\le \frac{1}{1+m}\left(2\Q^*[u_0]\right)^{\frac{2}{1-m}}\frac{(T-t_0)^{\frac{2}{1-m}}}{t_1-t_0}
= \c_*^{1+m} \Q^*[u_0]^{\frac{2}{1-m}}(T-t_1)^{\frac{1+m}{1-m}}
\end{split}
\end{equation}
where $\c_*>0$ only depends on $m$, recalling that $T-t_0= 3\tfrac{T-t_1}{2}$ and $t_1-t_0=\tfrac{T-t_1}{2}$\,.

A similar argument, exploiting the first inequality in \eqref{energy.H1.Hstar} allows to obtain analogous extinction rates for $\|\nabla u^m(t)\|_{\LL^2(\Omega)}$. Summing up, we have proven the following:
\begin{prop}\label{Prop.Q*-estimates}
Let $m\in (0,1)$, let $u$ be a   WDS to \eqref{CDP}    corresponding to the initial datum $0\le u_0\in H^{-1}(\Omega)\cap \LL^{1+m}(\Omega)$ and let $T=T(u_0)$ be its extinction time. We have that for all $t\in [0,T]$
\[
\|u(t)\|_{H^{-1}(\Omega)} \le (2\Q^*[u_0])^{\frac{1}{1-m}} (T-t)^{\frac{1}{1-m}}\,.
\]
Moreover, there exist $\c_*>0$ depending only on $m$ such that for all $t\in [0,T]$
\begin{equation}\label{ext.1+m}
\|u(t)\|_{\LL^{1+m}(\Omega)} \le \c_* \Q^*[u_0]^{\frac{2}{1-m^2}}\begin{cases}
T^{\frac{2}{1-m^2}} \,t^{- \frac{1}{1+m}} &\qquad\mbox{when }  0\le  t< \frac{T}{3}\\
(T-t)^{\frac{1}{1-m}}  &\qquad\mbox{when }\frac{T}{3}\le  t\le T\,,
\end{cases}
\end{equation}
and also
\begin{equation}\label{ext.grad}
\|\nabla u^m(t)\|_{\LL^{2}(\Omega)}^2 \le \c_* \Q^*[u_0]^{\frac{2}{1-m}}\begin{cases}
T^{\frac{2}{1-m}} \, t^{-2}  &\qquad\mbox{when }  0\le  t< \frac{T}{3}\\
(T-t)^{\frac{2m}{1-m}}  &\qquad\mbox{when }\frac{T}{3}\le t\le T\,.
\end{cases}
\end{equation}
\end{prop}
\subsubsection{Estimates of some weighted $L^p$ norms and of the extinction time}\label{Sssec.Extinction.Norms}
\begin{lem}[$L^1_{\Phi_1}$-norm estimates]
Let $u$ be a   WDS to \eqref{CDP}    corresponding to the initial datum $0\le u_0\in  \LL^1_{\Phi_1}(\Omega)$. Then we have that the extinction time $T=T(u_0)$ can be estimated from below as
\begin{equation}\label{lower.ext.time}
T\ge \c_0^{-1}\|u_0\|_{\LL^1_{\Phi_1}(\Omega)}^{1-m} \qquad\mbox{where}\qquad \c_0:= \lambda_1(1-m)\|\Phi_1\|_{1}^{1-m}\,.
\end{equation}
Moreover, for all $0\le t_0\le t  \le T$
\begin{equation}\label{L1-Phi1.ineq}
 \left(\int_\Omega u(t_0)\Phi_1\dx\right)^{1-m} - \c_0(t-t_0)\le \left(\int_\Omega u(t)\Phi_1\dx\right)^{1-m}\le\c_0(T-t)\,.
\end{equation}
\end{lem}
\noindent {\bf Proof.~}We differentiate in time and use H\"older inequality to get
\[
\frac{\rd}{\dt}\int_\Omega u\Phi_1\dx=\int_\Omega u^m \Delta \Phi_1\dx=-\lambda_1\int_\Omega u^m\Phi_1\dx \ge -\lambda_1\|\Phi_1\|_{1}^{1-m}\left(\int_\Omega u\Phi_1\dx\right)^m\,,
\]
which integrated over $[t_0,t]$ gives the first of inequalities \eqref{L1-Phi1.ineq}. The inequality on the right in \eqref{L1-Phi1.ineq} follows by letting $t= T$ in the first inequality of \eqref{L1-Phi1.ineq}. Letting $t_0=0$ and $t=T$ in the first  inequality of \eqref{L1-Phi1.ineq} proves \eqref{lower.ext.time}.\qed

\begin{lem}[$L^p$-norms estimates and extinction rates]\label{Lem.Lp.norms}
Let $m\in (0,1)$, $p> 1$ and $p\ge p_c$, and let $u$ be the   WDS to \eqref{CDP}    corresponding to $0\le u_0\in \LL^p(\Omega)$. Then $u$ extinguishes at a finite time $T=T(u_0)$ that can be estimated from above as follows:
\begin{equation}\label{upper.ext.time}
T\le \c_p\|u_0\|_{\LL^p(\Omega)}^{1-m} \qquad\mbox{where}\qquad \c_p:=\frac{(p+m-1)^2}{4m (1-m)  p(p-1)}\mathcal{S}_2^{2}\,,
\end{equation}
where $\mathcal{S}_2>0$ is the constant in the Sobolev-Poincar\'e inequality (see \eqref{Sob.Poinc} below).
Moreover, for all $t \ge t_0 \ge 0$ we have
\begin{equation}\label{Lp.est}
\c_p^{-1}(T-t)\le \|u(t)\|_{\LL^p(\Omega)}^{1-m}\le \|u(t_0)\|_{\LL^p(\Omega)}^{1-m}- \c_p^{-1}(t-t_0)\,.
\end{equation}
Furthermore  we have the following extinction rates: for all $t\in (0,T]$ and all $p> 1$ with $p> p_c$ there exists $\c_{m,p,\Omega}>0$ such that
\begin{equation}\label{Lp.est.2}
\|u(t)\|_{\LL^p(\Omega)}^p \le \c_{m,p,\Omega} \Q^*[u_0]^{\frac{2m}{1-m}}\frac{T^2}{ t^2}\left( \frac{\|u_0\|_{\LL^p(\Omega)}^{2p\vartheta_p}}{t^{N\vartheta_p}}\right)^{(p-(1+m))_+}
(T-t)^{\frac{1+m}{1-m}}\,.
\end{equation}
\end{lem}
\begin{rem}[About optimal extinction rates]\label{Sharp.Lp.extinction}\rm~\\[-6mm]
\begin{enumerate}[leftmargin=16pt, label=(\roman*)]\itemsep2pt \parskip2pt \parsep0pt
\item Note that in the limit $p\to 1^+$ the constant $\c_p^{-1}\to 0$, but $\c_{m,p,\Omega}$ remains strictly positive and finite.
\item When $m>m_s$ we can always take $p=1+m>p_c$ and obtain that
\begin{equation}\label{sharp.L1+m.bounds}
\c_{1+m}^{-\frac{1}{1-m}}(T-t)^{\frac{1}{1-m}}\le \|u(t)\|_{\LL^{1+m}(\Omega)} \le \c_{m,\Omega}^{\frac{1}{1+m}} \Q^*[u_0]^{\frac{2m}{1-m^2}}
(T-t)^{\frac{1}{1-m}}\,.
\end{equation}
This is the sharp extinction decay for the $\LL^{1+m}$ norm, and it has been shown for the first time in \cite{BH}, see also \cite{DKV1991,BGV2012-JMPA,JX2021,BII2022}. The novelty here is represented by the fact that we do not need to assume $\Q[u_0]<\infty$, i.e., $|\nabla u_0^m|\in \LL^2(\Omega)$, just that $\Q^*[u_0]<\infty$, that is $u_0\in H^{-1}(\Omega)\cap \LL^{1+m}(\Omega)$.

\item In the critical case $m=m_s$, we have that $p_c=1+m$. On the one hand, the bounds \eqref{sharp.L1+m.bounds} continue to hold, but they do not imply the $\LL^\infty$ decay with the same power. Indeed, in this case, $u_0\in L^{1+m}$ does not necessarily produce bounded solutions, see Proposition \ref{Prop.Linfty.ext.T} below.

\item In the subcritical case $m<m_s$ we have $p>p_c>1+m$ and we can only prove that
\[
\c_{p}^{-\frac{p}{1-m}}(T-t)^{\frac{p}{1-m}}\le \|u(t)\|_{\LL^{p}(\Omega)}^p\le \c_{m,p,\Omega} \Q^*[u_0]^{\frac{2m}{1-m}} \frac{T^2}{t^2} \left( \frac{\|u_0\|_{\LL^p(\Omega)}^{2p\vartheta_p}}{t^{N\vartheta_p}}\right)^{(p-(1+m))_+}
(T-t)^{\frac{1+m}{1-m}}\,.
\]
As we shall see, having upper bounds for such norms with the ``optimal power'' is impossible in general. Indeed, this would imply a GHP up to the extinction time, which does not hold in smooth and star-shaped domain, see Section \ref{sec.GHP} and in particular Theorem \ref{thm.NoGHP}. See also Section \ref{Ssec.Asymp.Subcrit.New} for more details about the behaviour of $\LL^p$ norms close to the extinction time.
\end{enumerate}
\end{rem}

\noindent {\bf Proof of Lemma~\ref{Lem.Lp.norms}.~}We split the proof into two Steps.

\noindent$\bullet~$\textsc{Step 1. }We first prove \eqref{Lp.est}, since it implies  \eqref{upper.ext.time} by letting $t=T$. Recall that
\begin{equation}\label{lem.Lp.decay.1}
\begin{split}
\frac{\rd}{\dt}\int_\Omega u(t,x)^p\dx &=-\frac{4mp(p-1)}{(p+m-1)^2}\int_\Omega \left|\nabla u^{\frac{p+m-1}{2}}\right|^2\dx\\
&\le -\frac{4mp(p-1)}{(p+m-1)^2}\mathcal{S}_2^{-2} \left(\int_\Omega u(t,x)^p\dx\right)^{1-\frac{1-m}{p}}
\end{split}
\end{equation}
where $\mathcal{S}_2^2$ is the constant in the Sobolev-Poincaré inequality:
\begin{equation}\label{Sob.Poinc}
\mathcal{S}_2^{-2}\|f\|_{\LL^{\frac{2p}{p+m-1}}(\Omega)}^2\le \|\nabla f\|_{\LL^2(\Omega)}^2\,.
\end{equation}
Note that
\[
\frac{2p}{p+m-1}\le 2^\ast \qquad\mbox{if and only if}\qquad p\ge p_c=\frac{N(1-m)}{2}\,.
\]
Integrating over $[t_0,t_1]$ we obtain
\[
\left(\int_\Omega u(t_1,x)^p\dx\right)^{\frac{1-m}{p}} \le \left(\int_\Omega u(t_0,x)^p\dx\right)^{\frac{1-m}{p}} - \frac{4m (1-m)  p(p-1)}{(p+m-1)^2}\mathcal{S}_2^{-2}(t_1-t_0)\,,
\]
from which we deduce that there exists an extinction time and that \eqref{Lp.est} follows.

\noindent$\bullet~$\textsc{Step 2. }\textit{Extinction rates with $\Q^*$. }We only have to deal with the case $p>1+m$, since otherwise the estimate follows by the decay of the $L^{1+m}$ norm \eqref{ext.1+m} simply by H\"older inequality. We next estimate the following integral for all $t\in [t_0,T]$:
\begin{align}\label{norm.decay.2.1} 
 \int_\Omega \left|\nabla u^{\frac{p+m-1}{2}}(t)\right|^2 &\dx
 = \frac{(p+m-1)^2}{4m^2}\int_\Omega u^{p-(1+m)} \left|\nabla u^m(t)\right|^2 \dx\nonumber\\
 &\le  \frac{(p+m-1)^2}{4m^2}\|u(t)\|_{\LL^\infty(\Omega)}^{p-(1+m)}\|\nabla u^m(t)\|_{\LL^2(\Omega)}^2 \nonumber\\
  &\le  \frac{(p+m-1)^2}{4m^2}\|u(t)\|_{\LL^\infty(\Omega)}^{p-(1+m)}  \frac{\c_*T^2}{t^2}  \Q^*[u_0]^{\frac{2}{1-m}}(T-t)^{\frac{2m}{1-m}}  \nonumber\\
 &\le \frac{\c_*T^2}{t_0^2} \frac{(p+m-1)^2}{4m^2}\Q^*[u_0]^{\frac{2}{1-m}}\left(\ka \frac{\|u_0\|_{\LL^p(\Omega)}^{2p\vartheta_p}}{t_0^{N\vartheta_p}}\right)^{p-(1+m)}(T-t)^{\frac{2m}{1-m}\,,}\nonumber
\end{align}
where in the second inequality we have estimated $\|\nabla u^m(t)\|_{\LL^2(\Omega)}^2 $ using Proposition~\ref{Prop.Q*-estimates}, while in the third we have used the smoothing effects \eqref{Smoothing.Global} and   that $t\ge t_0$.
Combining the above estimate with \eqref{lem.Lp.decay.1} we get
\[
\begin{split}
\frac{\rd}{\dt}\int_\Omega u(t,x)^p\dx &=-\frac{4mp(p-1)}{(p+m-1)^2}\int_\Omega \left|\nabla u^{\frac{p+m-1}{2}}\right|^2\dx\\
&\ge - \frac{\c_*T^2}{t_0^2}  \frac{p(p-1)}{m} \Q^*[u_0]^{\frac{2}{1-m}}\left(\ka \frac{\|u_0\|_{\LL^p(\Omega)}^{2p\vartheta_p}}{t_0^{N\vartheta_p}}\right)^{p-(1+m)}(T-t)^{\frac{2m}{1-m}}\,,
\end{split}
\]
which integrated over $[t_0,t_1]$ gives
\[
\|u(t_1)\|_{\LL^p(\Omega)}^p - \|u(t_0)\|_{\LL^p(\Omega)}^p \ge \c_{m,p} \frac{T^2}{t_0^2} \Q^*[u_0]^{\frac{2}{1-m}}\left(  \frac{\|u_0\|_{\LL^p(\Omega)}^{2p\vartheta_p}}{t_0^{N\vartheta_p}}\right)^{p-(1+m)}
\left[(T-t_1)^{\frac{1+m}{1-m}}-(T-t_0)^{\frac{1+m}{1-m}}\right]
\]
and inequality \eqref{Lp.est.2} follows by letting $t_1=T$ and $t_0=t$.\qed

\begin{lem}[Monotonicity of $L^p_{\Phi_1}$-norms]\label{lem.Lp.weighted}
Let $m\in (0,1)$ and let $u$ be the   WDS to \eqref{CDP}    corresponding to $0\le u_0\in \LL^p_{\Phi_1}(\Omega)$, with $p\geq 1$. Then we have that
\begin{equation}\label{lem.Lp.weighted.monot.Lp}
\|u(t)\|_{\LL^p_{\Phi_1}(\Omega)}\le \|u(t_0)\|_{\LL^p_{\Phi_1}(\Omega)}\qquad\mbox{for all }t\ge t_0\ge 0\,.
\end{equation}
\end{lem}
\noindent {\bf Proof.~}Differentiating the weighted $\LL^p$ norm and integrating by parts (a rigorous proof can be obtained by approximation by means of smooth solutions) we obtain:
\[
\begin{split}
\frac{\rd}{\dt}\int_\Omega u(t,x)^p\Phi_1 \dx &=-\frac{4mp(p-1)}{(p+m-1)^2}\int_\Omega \left|\nabla u^{\frac{p+m-1}{2}}\right|^2\Phi_1 \dx
    +\frac{m p  }{p+m-1}\int_\Omega u^{p+m-1}\Delta\Phi_1 \dx\\
& =-\frac{ p(p-1)}{m}\int_\Omega u^{p-(1+m)} \left|\nabla u^m\right|^2\Phi_1 \dx -\frac{m p \lambda_1 }{p+m-1}\int_\Omega u^{p+m-1} \Phi_1 \dx\le 0
\end{split}
\]
which immediately gives the monotonicity of the weighted norm \eqref{lem.Lp.weighted.monot.Lp}.\qed

\subsubsection{Upper $L^\infty$ bounds close to the extinction time}\label{SSec.Linfty.ext.T}

We conclude by proving upper bounds for the $\LL^\infty$ norm close to the extinction time.

\begin{prop}\label{Prop.Linfty.ext.T}Let $m\in (0,1)$ and $p>\max\{1, p_c\}$, let $u$ be the   WDS to \eqref{CDP}    corresponding to $0\le u_0\in \LL^p(\Omega)$, and let $T=T(u_0)$ be its extinction time. Then,  for all $t\in [\tfrac{T}{2},T]$ we have
\begin{equation}\label{decay.Linfty.p}
\|u(t)\|_{\LL^\infty(\Omega)}\le \overline{\c} \Q^*[u_0]^{\frac{4m}{1-m}\vartheta_p}
\left(\frac{\|u_0\|_{\LL^p(\Omega)}^{2p\vartheta_p}}{T^{N\vartheta_p}}\right)^{2\vartheta_p(p-(1+m))_+}
(T-t)^{\frac{1}{1-m} -\frac{2\vartheta_p}{1-m}[p-(1+m)]}\,.
\end{equation}
  where $\overline{\c}>0$ depends only on $m,N,p,\Omega$.
\end{prop}
\begin{rem}
  \rm When $m<m_s$ we have that $p>p_c>1+m$, so the  estimates above gives an upper bound for the explosion rate of the $\LL^\infty$ norm. This will be discussed  thoroughly in Section \ref{Ssec.Asymp.Subcrit.New}.
\end{rem}

\noindent {\bf Proof.~}For all $t>T/2$ we can choose
\[
t_0=t-\frac{T-t}{2}\quad\mbox{so that}\quad t-t_0 =\frac{T-t}{2}
\quad\mbox{and}\quad T-t_0 = 3\frac{T-t}{2} \,.
\]
Notice that we have always $t_0\le t$ and $t_0>T/4$ (since $t>T/2$). Using the smoothing effects \eqref{Smoothing.Global} together with the decay of the $\LL^p$ norm \eqref{Lp.est.2} we obtain
\[
\begin{split}
\|u(t)\|_{\LL^\infty(\Omega)}
\le \ka \frac{\|u(t_0)\|_{\LL^p(\Omega)}^{2p\vartheta_p}}{(t-t_0)^{N\vartheta_p}}
\le \ka \left[\c_{m,p,\Omega} \Q^*[u_0]^{\frac{2m}{1-m}}
\left( \frac{\|u_0\|_{\LL^p(\Omega)}^{2p\vartheta_p}}{t_0^{N\vartheta_p}}\right)^{(p-(1+m))_+}\right]^{2\vartheta_p}
\frac{(T-t_0)^{2\vartheta_p\frac{1+m}{1-m}}}{(t-t_0)^{N\vartheta_p}}\\
\le \overline{\c} \Q^*[u_0]^{\frac{4m}{1-m}\vartheta_p}
\left( \frac{\|u_0\|_{\LL^p(\Omega)}^{2p\vartheta_p}}{T^{N\vartheta_p}}\right)^{2\vartheta_p(p-(1+m))_+}
(T-t)^{\frac{1}{1-m}\frac{\vartheta_{p}}{\vartheta_{1+m}}}
\end{split}
\]
We conclude the proof of \eqref{decay.Linfty.p} by observing that
\[
\frac{\vartheta_{p}}{\vartheta_{1+m}}= 1-2[p-(1+m)]\vartheta_p \mbox{\,.\qed}
\]

When $m>m_{c,1}$ (see \eqref{eq:mc1}) we can get a better estimate.

\begin{prop}\label{prop.smoothing.L1phi}Let $m\in (m_{c,1},1)$, let $u$ be the   WDS to \eqref{CDP}    corresponding to $0\le u_0\in \LL^1_{\Phi_1}(\Omega)$, and let $T=T(u_0)$ be its extinction time. Then,  for all $t\in [\tfrac{T}{2},T]$ we have
\begin{equation}\label{decay.Linfty.11}
\|u(t)\|_{\LL^\infty(\Omega)}\le \overline{\c}
(T-t)^{\frac{1}{1-m}}\,,
\end{equation}
where $\overline{\c}>0$ only depends on $m,N,\Omega$.
\end{prop}
\noindent {\bf Proof.~}
\[
\begin{split}
\|u(t)\|_{\LL^\infty(\Omega)}
\le \ka \frac{\|u(t_0)\|_{\LL^1_{\Phi_1}(\Omega)}^{\vartheta_{1,1}}}{(t-t_0)^{N\vartheta_{1,1}}}
=\ka  \left(\frac{\|u(t_0)\|_{\LL^1_{\Phi_1}(\Omega)}}{(t-t_0)^{\frac{1}{1-m}}}\right)^{\vartheta_{1,1}}(t-t_0)^{\frac{1}{1-m} }
\le \overline{\c} (T-t)^{\frac{1}{1-m} }\,,
\end{split}
\]
where we have used the $L^1_{\Phi_1}$ upper estimates \eqref{L1-Phi1.ineq}
\[
\left(\int_\Omega u(t_0)\Phi_1\dx\right)^{1-m}\le\c_0(T-t_0)
\]
where $\c_0$ only depends on $m,N,\Omega$. Next we have chosen $t\ge T/3$ and $t_0=t-\tfrac{T-t}{2}$ (notice that $t_0\ge 0$ since $t\ge T/3$) so that $T-t_0= 3\tfrac{T-t}{2}$ and $t-t_0=\tfrac{T-t}{2}$\,.\qed

\subsection{Upper boundary estimates for all times}\label{Ssec.Upper.Bdry.T}

The results of the previous sections allow us to prove the upper part of the Global Harnack Principle.

\begin{thm}[Upper boundary estimates up to extinction time]\label{thm.GHP.up.T}
Let $m_c=\frac{N-2}{N}$ and $p_c=\tfrac{N(1-m)}{2}$, and let $u$ be a   WDS to \eqref{CDP}    corresponding to   $u_0\in \LL^p(\Omega)\cap H^{-1}(\Omega)$ with $p\ge 1$ if $m\in (m_c,1)$ and $p>p_c$ if $m\in (0,m_c]$. Then there exists a   $\ka[u_0,t]>0$    such that, for all $t\in (0,T]$,
\begin{equation}\label{thm.GHP.up.T.ineq}\begin{split}
\left\|\frac{u^m(t)}{\Phi_1}\right\|_{\LL^\infty(\Omega)}
&\le \ka[u_0,t] \,
(T-t)^{\frac{m}{1-m}-\frac{2\vartheta_p }{1-m}[p-(1+m)]_+} \,.
\end{split}
\end{equation}
Here   $\ka[u_0,t] $   depends on $m,p,N,\Omega$ and $u_0$ and has the form
\begin{equation}\label{thm.GHP.up.T.ineq.const}
 \ka[u_0,t]:= c_{m,p,N,\Omega} \left\{\begin{array}{lll}
   \ka   \|u_0\|_{L^p(\Omega)}^{2p\vartheta_p} \, t^{-\frac{2p\vartheta_p}{1-m}} & \quad\mbox{when }t\in\left(0,\frac{2}{3}T\right)\,,\\
   \Q^*[u_0]^{\frac{4m}{1-m}\vartheta_p}
    \left(\frac{\|u_0\|_{\LL^p(\Omega)}^{2p\vartheta_p}}{T^{N\vartheta_p}}\right)^{2\vartheta_p(p-(1+m))_+} &  \quad\mbox{when }t\in\left[\frac{2}{3}T, T\right]\,,
 \end{array}\right.
\end{equation}
where $c_{m,p,N,\Omega},   \ka  >0$ depend only on $m,p,N,\Omega$.
\end{thm}
\noindent {\bf Proof.~}Recall the fundamental upper bounds \eqref{Prop.PE.ineq}: for all $0\le t_0\le t_1$
\begin{align}
u^m(t_1,x)&\le \frac{t_1^{\frac{m}{1-m}}}{1-m}\int_\Omega  \frac{u(t_0,y)-u(t_1,y)}{ t_1^{\frac{1}{1-m}}-t_0^{\frac{1}{1-m}} } \G(x,y)\dy
\le \frac{1}{1-m}\frac{t_1^{\frac{m}{1-m}}}{ t_1^{\frac{1}{1-m}}-t_0^{\frac{1}{1-m}} } \|u(t_0)\|_{\LL^\infty(\Omega)} \|\G(x,\cdot)\|_{\LL^1(\Omega)}\nonumber\\
&\le c_1 \frac{\Phi_1(x)}{1-m}\frac{t_1^{\frac{m}{1-m}}}{ t_1^{\frac{1}{1-m}}-t_0^{\frac{1}{1-m}} }\|u(t_0)\|_{\LL^\infty(\Omega)}\,. \label{Prop.PE.ineq1x1}
\end{align}
Here we have used the Green function estimate \eqref{G4} which implies that $\|\G(x,\cdot)\|_{\LL^1(\Omega)} \le c_1\Phi_1(x)$, where $c_1>0$ depends only on $N$ and $\Omega$, see for instance Lemma 4.2 of \cite{BFV2018} (taking $s=\gamma=1$ there).

As in Theorem \ref{thm.GHP.up}, when $t\in (0,T/2)$ we can use the smoothing effects \eqref{Smoothing.Global} or \eqref{Smoothing.Global.Weighted} and   \eqref{Prop.PE.ineq1x1}    to obtain \eqref{Smoothing.p.bdry}, that is, for all $t>0$ (recall that $\vartheta_p=\tfrac{1}{2p-N(1-m)}$)
\begin{equation}\label{Smoothing.p.bdry.111}
\begin{split}
\left\|\frac{u^m(t)}{\Phi_1}\right\|_{\LL^\infty(\Omega)}
&\le  \ka
\frac{\|u_0\|_{L^p(\Omega)}^{2p\vartheta_p}}{t^{N\vartheta_p+1}}
=\ka  \frac{\|u_0\|_{L^p(\Omega)}^{2p\vartheta_p}}{t^{\frac{2p\vartheta_p}{1-m}}}  t^{\frac{m}{1-m}}
\le 2^{\frac{m}{1-m}}\ka  \frac{\|u_0\|_{L^p(\Omega)}^{2p\vartheta_p}}{t^{\frac{2p\vartheta_p}{1-m}}}  (T-t)^{\frac{m}{1-m}}\,,
\end{split}
\end{equation}
  since $t\le 2(T-t)$ when $t\le \frac{2}{3}T$. This is exactly \eqref{thm.GHP.up.T.ineq} when $t\le \frac{2}{3}T$.

  When $t\in [\tfrac{2}{3}T,T]$,    the upper bound \eqref{thm.GHP.up.T.ineq} follows from  \eqref{Prop.PE.ineq1x1}  using the convexity of $u\mapsto u^{\frac{1}{1-m}}$ and the upper bounds \eqref{decay.Linfty.p}. More precisely,   we start from \eqref{Prop.PE.ineq1x1} that holds for all $0\le t_0\le t_1$, and estimate $\|u(t_0)\|_{\LL^\infty(\Omega)}$ with \eqref{decay.Linfty.p} (and this requires $t_0\ge T/2$) to obtain
\[
u^m(t_1,x)\le \frac{c_1 \Phi_1(x)}{ t_1 -t_0 }  \overline{\c} \Q^*[u_0]^{\frac{4m}{1-m}\vartheta_p}
\left(\frac{\|u_0\|_{\LL^p(\Omega)}^{2p\vartheta_p}}{T^{N\vartheta_p}}\right)^{2\vartheta_p(p-(1+m))_+}
(T-t_0)^{\frac{1}{1-m}\frac{\vartheta_{p}}{\vartheta_{1+m}}}\,.
\]
Choosing $t_1=t$ and
\begin{equation}\label{choice.t0}
t_0=t-\frac{T-t}{2}\quad\mbox{so that}\quad t-t_0 =\frac{T-t}{2}
\quad\mbox{and}\quad T-t_0 = 3\frac{T-t}{2} \,,
\end{equation}
 Note that $t_0\ge T/2$ when $t\ge \frac{2}{3}T$. As a consequence, we obtain
\[
\left\|\frac{u^m(t)}{\Phi_1}\right\|_{\LL^\infty(\Omega)}\le  c_1    \overline{\c} \Q^*[u_0]^{\frac{4m}{1-m}\vartheta_p}
\left(\frac{\|u_0\|_{\LL^p(\Omega)}^{2p\vartheta_p}}{T^{N\vartheta_p}}\right)^{2\vartheta_p(p-(1+m))_+}
(T-t)^{\frac{1}{1-m}\frac{\vartheta_{p}}{\vartheta_{1+m}}-1}\,.
\]
Since $\vartheta_p/\vartheta_{1+m}=m-2\vartheta_p[p-(1+m)]$,
this concludes the proof.\qed

When $m>m_{c,1}$ we can remove the dependence on $u_0$ from the constant, close to the extinction time.
\begin{thm}[Upper boundary estimates up to extinction time when $m>m_{c,1}$]\label{thm.GHP.up.T.1}
Let $m\in (m_{c,1},1)$, let $u$ be the   WDS to \eqref{CDP}    corresponding to $0\le u_0\in \LL^1_{\Phi_1}(\Omega)$, and let $T=T(u_0)$ be its extinction time. Then there exists a constant $\ka>0$ depending only on $m,N,\Omega$ such that for all $t\in [\tfrac{T}{2},T]$
\begin{equation}\label{thm.GHP.up.T.ineq.1}\begin{split}
\left\|\frac{u^m(t)}{\Phi_1}\right\|_{\LL^\infty(\Omega)}
&\le \ka \, (T-t)^{\frac{m}{1-m}} \,.
\end{split}
\end{equation}
\end{thm}
\noindent {\bf Proof.~}Repeat the proof of the Theorem \ref{thm.GHP.up.T} until formula \eqref{Prop.PE.ineq1x1}, then use the upper bounds  \eqref{decay.Linfty.11} (whose constant does not depend on $u_0$) instead of \eqref{decay.Linfty.p}, and finally choose $t_0$ as in \eqref{choice.t0}. \qed

\subsection{Lower boundary estimates for all times}\label{Ssec.Lower.Bdry}
We shall prove here the lower part of the GHP, by splitting the result into several useful lower bounds.
\begin{thm}[Global lower bounds I]\label{thm.GHP.low.I}
  Let $m\in (0,1)$ and   let $u$ be a   WDS to \eqref{CDP}    corresponding to $0\le u_0\in \LL^1_{\Phi_1}(\Omega)$ and let $T=T(u_0)>0$ be its extinction time.  Then for all $0\le t\le T$ and all $x\in \overline{\Omega}$
\begin{equation}\label{thm.GHP.low.I.ineq}
\frac{u^m(t,x)}{\Phi_1(x) }\ge \frac{c_0 t^{\frac{m}{1-m}}\|u(t)\|_{\LL^1_{\Phi_1}(\Omega)}}{(1-m)\left(T^{\frac{1}{1-m}}-t^{\frac{1}{1-m}}\right)}\,.
\end{equation}
 where $c_0>0$ is the constant in the assumption \eqref{G4}, that depends only on $N,\Omega$.
\end{thm}
\noindent {\bf Proof.~}We use the pointwise upper bound of \eqref{Prop.PE.ineq}  together with the fact that $u\ge 0$ and $T=T(u_0)$ is the extinction time, to get that for all $0\le t \le T$ and a.e. $x\in \Omega$
\[
\begin{split}
\frac{u^m(t,x)}{t^{\frac{m}{1-m}}}
 &\ge \frac{1}{1-m}\int_\Omega  \frac{u(t,y)- \cancel{u(T,y)}}{{T^{\frac{1}{1-m}}-t^{\frac{1}{1-m}}}}  \G(x,y)\dy
\ge \frac{c_0 \Phi_1(x)}{1-m}\int_\Omega  \frac{u(t,y)  \Phi_1(y)}{ T^{\frac{1}{1-m}}-t^{\frac{1}{1-m}} }\dy
\end{split}\]
where in the last step we have used \eqref{G4} in the weaker form $\G(x,y)\ge c_0\Phi_1(x)\Phi_1(y)$\,.\qed

\begin{lem}[$L^1_{\Phi_1}$-norm lower estimates]\label{lem.lower.L1}
Let $m\in (0,1)$ and $u$ be a   WDS to \eqref{CDP}    corresponding to  $0\le u_0\in \LL^p(\Omega)$ with $p\ge 1$ if $m\in (m_c,1)$ and $p>p_c$ if $m\in (0,m_c]$, where $m_c=\frac{N-2}{N}$ and $p_c=\tfrac{N(1-m)}{2}$. Let $T=T(u_0)>0$ be the extinction time. Then there exists  $\underline{\c}[u_0,t]>0$ such that
\begin{equation}\label{lem.lower.L1.ineq}
\int_\Omega u(t)\Phi_1\dx  
\ge   \underline{\c}[u_0,t]   (T-t)^{\frac{1}{1-m}+\frac{2\vartheta_p}{1-m}(p-2m+1)\,[p-(1+m)]_+}\qquad\mbox{for all $t\in [0,T]$.}
\end{equation}
  Where
\begin{equation}\label{c0.low}
\underline{\c}[u_0,t]=  \left( \frac{t}{T}\wedge \frac{2}{3} \right)^{\frac{2p(p-m+1)\vartheta_p}{1-m}}
\begin{cases}
\underline{\c}_1[u_0,T] & \mbox{when }t\in (0,\frac{2}{3}T)\\
\underline{\c}_2[u_0,T] &\mbox{when } t\ge \frac{2}{3}T\\
\end{cases}
\end{equation}
Where
\begin{equation}\label{c0.low-1}
\underline{\c}_1[u_0,T]=\frac{\lambda_1  \,T^{1+\frac{2p(p-m+1)\vartheta_p}{1-m}}  }{ \c_p ^{\frac{p}{1-m}}2^{\frac{p-m+1}{1-m}} \ka^{p-m+1}  \|u_0\|_{\LL^p(\Omega)}^{2p(p-m+1)\vartheta_p}}
\end{equation}
and
\begin{equation}\label{c0.low-2}
\underline{\c}_2[u_0,T]=
\tfrac{\lambda_1(1-m)  }{\c_p^{\frac{p}{1-m}}}
\left\{\begin{array}{lll}
\left(\overline{\c}^{1-m}\ka \right)^{-1}&\!\!\!\!\!\!\mbox{when }m_{c,1}<m<1\\
\left(\overline{\c}^{1-m}c_{m,p,N,\Omega}\right)^{-1}  \Q^*[u_0]^{-4m\frac{2-m}{1-m}\vartheta_{1+m}} &\!\!\!\!\!\!\mbox{when }m_s<m<m_{c,1}\\
\frac{\overline{\c}^{2m-p}c_{m,p,N,\Omega}^{-1}\Q^*[u_0]^{-\frac{4m}{1-m}\vartheta_p(p-2m+1)}}{1+ 2\vartheta_p (p-2m+1)[p-(1+m)] } \left[\frac{T^{N }}{\|u_0\|_{\LL^p(\Omega)}^{2p }}\right]^{2\vartheta_p^2(p-2m+1)(p-(1+m))_+}
&\!\!\!\!\!\!\mbox{when }0<m\le m_s\\
\end{array}
\right.
\end{equation}
where $\overline{\c},\ka$ depend only on $m,N,p$ and $\Omega$, see \eqref{decay.Linfty.11} and \eqref{thm.GHP.up.T.ineq.1}. Here  $\c_p,c_{m,p,N,\Omega}>0$ depend only on $m,p,N$ and $\Omega$, see \eqref{upper.ext.time} and \eqref{thm.GHP.up.T.ineq.const}. When $m>m_{c,1}$ we have that $\underline{\c}_2[u_0,T]$ does not depend  on $u_0$ nor on $T$.
\end{lem}
\noindent {\bf Proof.~}We differentiate in time and use H\"older inequality to get  for all $t\in [0,T]$:
\begin{equation}\label{lem.lower.L1.01}
\frac{\rd}{\dt}\int_\Omega u\Phi_1\dx=-\lambda_1\int_\Omega u^m\Phi_1\dx
\le  -\frac{\lambda_1 (T-t)^{\frac{p}{1-m}}}{\c_p^{\frac{p}{1-m}}  \|u\|_{\LL^\infty(\Omega)}^{p-2m}\left\|\frac{u^m}{\Phi_1}\right\|_{\LL^\infty(\Omega)}}\,.
\end{equation}
Here the inequality follows from
\[
\int_\Omega u^m\Phi_1\dx=\int_\Omega \frac{u^p}{u^{p-2m}}\frac{\Phi_1}{u^m}\dx \ge \frac{\|u\|_{\LL^p(\Omega)}^p}{\|u\|_{\LL^\infty(\Omega)}^{p-2m}\left\|\frac{u^m}{\Phi_1}\right\|_{\LL^\infty(\Omega)}}
\ge \frac{(T-t)^{\frac{p}{1-m}}}{\c_p^{\frac{p}{1-m}}\|u\|_{\LL^\infty(\Omega)}^{p-2m}\left\|\frac{u^m}{\Phi_1}\right\|_{\LL^\infty(\Omega)}}\,,
\]
where we used the lower bound on the $\LL^p$ norm \eqref{Lp.est}, and $\c_p$ is as in \eqref{upper.ext.time}.

 We split two steps, namely the case $t\in [0,\frac{2}{3}T]$  and $t\in [ \frac{2}{3}T, T]$, since  we are going to combine the above estimates with the upper bounds \eqref{decay.Linfty.p} on the $\LL^\infty$ norm, that takes different forms in the two time ranges.

\noindent$\bullet~$\textsc{Step 1. }\textit{The case $t\in [0,\frac{2}{3}T]$. }In this case \eqref{lem.lower.L1.01} becomes,
 \begin{equation}\label{lem.lower.L1.01x}
\frac{\rd}{\dt}\int_\Omega u\Phi_1\dx
\le  -\frac{\lambda_1 (T-t)^{\frac{p}{1-m}}}{\c_p^{\frac{p}{1-m}}\|u\|_{\LL^\infty(\Omega)}^{p-2m}\left\|\frac{u^m}{\Phi_1}\right\|_{\LL^\infty(\Omega)}}
\le -\frac{\lambda_1  t^{\frac{2p(p-m+1)\vartheta_p}{1-m}}}{\c_p^{\frac{p}{1-m}}\kb_1[u_0]}
=: - \kb_2[u_0] t^{\frac{2p(p-m+1)\vartheta_p}{1-m}}\,.
\end{equation}
since $t\le \frac{2}{3}T$ implies $T-t\ge \frac{T}{3}$, and we have used the smoothing effect \eqref{Smoothing.Global}  and the smoothing \eqref{thm.GHP.up.T.ineq.const} (both can be written in the form \eqref{Smoothing.p.bdry.111}, valid when  $t\le \frac{2}{3}T$), that is
\[
 \|u(t)\|_{\LL^\infty(\Omega)}\le 2^{\frac{1}{1-m}} \ka  \frac{\|u_0\|_{L^p(\Omega)}^{2p\vartheta_p}}{t^{\frac{2p\vartheta_p}{1-m}}}  (T-t)^{\frac{1}{1-m}}\quad\mbox{and}\quad
 \left\|\frac{u^m(t)}{\Phi_1}\right\|_{\LL^\infty(\Omega)}
 \le 2^{\frac{m}{1-m}} \ka  \frac{\|u_0\|_{L^p(\Omega)}^{2p\vartheta_p}}{t^{\frac{2p\vartheta_p}{1-m}}}  (T-t)^{\frac{m}{1-m}}
\]
that can be combined as follows
\[
\|u(t)\|_{\LL^\infty(\Omega)}^{p-2m} \left\|\frac{u^m(t)}{\Phi_1}\right\|_{\LL^\infty(\Omega)}
    \le 2^{\frac{p-m+1}{1-m}} \ka^{p-m+1} \frac{\|u_0\|_{\LL^p(\Omega)}^{2p(p-m+1)\vartheta_p}}{t^{\frac{2p(p-m+1)\vartheta_p}{1-m}}}(T-t)^{\frac{p}{1-m}}
    =: \frac{\kb_1[u_0](T-t)^{\frac{p}{1-m}}}{t^{\frac{2p(p-m+1)\vartheta_p}{1-m}}}\,.
\]
Integrating  \eqref{lem.lower.L1.01x} on $[t,T]$ we obtain
\[\begin{split}
\int_\Omega u(t)\Phi_1\dx
&\ge \frac{\kb_2[u_0] }{\frac{2p(p-m+1)\vartheta_p}{1-m}+1}\left(T^{\frac{2p(p-m+1)\vartheta_p}{1-m}+1}-t^{\frac{2p(p-m+1)\vartheta_p}{1-m}+1}\right)
\ge  \kb_2[u_0] t^{\frac{2p(p-m+1)\vartheta_p}{1-m}}\left(T -t \right)\\
\end{split}
\]
where we have used the convexity inequality $a^{\alpha+1} -b^{\alpha+1}\ge (\alpha+1)b^\alpha (a-b)$. Notice that since $\frac{1}{1-m}>1$ we have that $\frac{T-t}{T}\ge (\frac{T-t}{T})^{\frac{1}{1+m}+\frac{2\vartheta_p}{1-m}(p-2m+1)\,[p-(1+m)]_+} $ hence the above inequality implies
\[
\int_\Omega u(t)\Phi_1\dx \ge  \kb_3[u_0,T]
\left(  \frac{t}{T}  \right)^{\frac{2p(p-m+1)\vartheta_p}{1-m}}(T-t)^{\frac{1}{1-m}+\frac{2\vartheta_p}{1-m}(p-2m+1)\,[p-(1+m)]_+}
\]
which is exactly \eqref{lem.lower.L1.ineq} when $t\le \frac{2}{3}T$. Note that $\kb_3[u_0,T]$ depends on $u_0, T, m,N,p,\Omega$ and has the expression
\[
\kb_3[u_0,T]=\kb_2[u_0 ] T^{1+\frac{2p(p-m+1)\vartheta_p}{1-m}}
=\frac{\lambda_1  \,T^{1+\frac{2p(p-m+1)\vartheta_p}{1-m}}  }{ \c_p ^{\frac{p}{1-m}}2^{\frac{p-m+1}{1-m}} \ka^{p-m+1}  \|u_0\|_{\LL^p(\Omega)}^{2p(p-m+1)\vartheta_p}}
\]

\noindent$\bullet~$\textsc{Step 2. }\textit{The case $t\in[ \frac{2}{3}T, T]$. }We combine the inequality \eqref{lem.lower.L1.01} with the upper estimates \eqref{decay.Linfty.p} on the $\LL^\infty$ norm, and we split some cases.

\noindent \textit{First we consider $m\in (m_s,1)$, }so that we can take $p=m+1$ and the inequalities  \eqref{decay.Linfty.p}  and \eqref{thm.GHP.up.T.ineq.const} become
\begin{equation}\label{lem.lower.L1.02}\begin{split}
\|u(t)\|_{\LL^\infty(\Omega)}
    &\le \overline{\c}   \Q^*[u_0]^{\frac{4m}{1-m}\vartheta_{1+m}} (T-t)^{\frac{1}{1-m}}\\
\left\|\frac{u^m(t)}{\Phi_1}\right\|_{\LL^\infty(\Omega)}
    &\le c_{m,p,N,\Omega} \Q^*[u_0]^{\frac{4m}{1-m}\vartheta_{1+m}}
    (T-t)^{\frac{m}{1-m}}\\
\end{split}\end{equation}
Combining \eqref{lem.lower.L1.01} and \eqref{lem.lower.L1.02} we obtain
\[
\frac{\rd}{\dt}\int_\Omega u\Phi_1\dx
\le  -\frac{\lambda_1 (T-t)^{\frac{1+m}{1-m}}}{ \c_{1+m}^{\frac{1+m}{1-m}}  \|u\|_{\LL^\infty(\Omega)}^{1-m}\left\|\frac{u^m}{\Phi_1}\right\|_{\LL^\infty(\Omega)}}
\le  -\frac{\lambda_1 (T-t)^{\frac{m}{1-m}}}{\c_{1+m}^{\frac{1+m}{1-m}} \overline{\c}^{1-m}c_{m,p,N,\Omega}  \Q^*[u_0]^{4m\frac{2-m}{1-m}\vartheta_{1+m}}}\,,
\]
which integrated over $[t,T]$ gives
\begin{equation}\label{lem.lower.L1.03b}
\int_\Omega u(t)\Phi_1\dx  \ge \frac{\lambda_1(1-m) (T-t)^{\frac{1}{1-m}}}{\c_{1+m}^{\frac{1+m}{1-m}} \overline{\c}^{1-m}c_{m,p,N,\Omega}  \Q^*[u_0]^{4m\frac{2-m}{1-m}\vartheta_{1+m}}}:=\underline{\c}(u_0)(T-t)^{\frac{1}{1-m}}\,.
\end{equation}
This proves \eqref{lem.lower.L1.ineq} when $m>m_s$.

\noindent\textit{When $m>m_{c,1}$ }we have the $\LL^1_{\Phi_1}-\LL^\infty$ smoothing effects, i.e., we are allowed to use Proposition \ref{prop.smoothing.L1phi}, and Theorem  \eqref{thm.GHP.up.T.1}
so that for all $t\in [\tfrac{T}{2},T]$ we have
\begin{equation}\label{lem.lower.L1.03c}
\|u(t)\|_{\LL^\infty(\Omega)}\le \overline{\c}
(T-t)^{\frac{1}{1-m}}\qquad\mbox{and}\qquad \left\|\frac{u^m(t)}{\Phi_1}\right\|_{\LL^\infty(\Omega)}
\le \ka \, (T-t)^{\frac{m}{1-m}} \,,
\end{equation}
where $\overline{\c}>0$ only depends on $m,N,\Omega$. Combining \eqref{lem.lower.L1.01} and \eqref{lem.lower.L1.03c} we obtain
\[
\frac{\rd}{\dt}\int_\Omega u\Phi_1\dx
\le  -\frac{\lambda_1 (T-t)^{\frac{1+m}{1-m}}}{\c_p^{\frac{1+m}{1-m}}\|u\|_{\LL^\infty(\Omega)}^{1-m}\left\|\frac{u^m}{\Phi_1}\right\|_{\LL^\infty(\Omega)}}
\le  -\frac{\lambda_1 (T-t)^{\frac{m}{1-m}}}{\c_{1+m}^{\frac{1+m}{1-m}} \overline{\c}^{1-m}   \ka}\,,
\]
which integrated over $[t,T]$ gives
\begin{equation}\label{lem.lower.L1.03b}
\int_\Omega u(t)\Phi_1\dx  \ge \frac{\lambda_1(1-m) (T-t)^{\frac{1}{1-m}}}{\c_{1+m}^{\frac{1+m}{1-m}} \overline{\c} }:=\underline{\c}(u_0)(T-t)^{\frac{1}{1-m}}\,.
\end{equation}

\noindent\textit{In the remaining case,  when $m\le m_s$, }we have always that $p>1+m$ and $\vartheta_{m+1}<0$, so that the inequality  \eqref{decay.Linfty.p} becomes (notice that here $p-2m>1-m>0$)
\begin{equation}\label{lem.lower.L1.04}\begin{split}
\|u(t)\|_{\LL^\infty(\Omega)}&\le \overline{\c} \Q^*[u_0]^{\frac{4m}{1-m}\vartheta_{p}}
\left(\frac{\|u_0\|_{\LL^p(\Omega)}^{2p\vartheta_p}}{T^{N\vartheta_p}}\right)^{2\vartheta_p(p-(1+m))_+}
(T-t)^{\frac{1}{1-m} -\frac{2\vartheta_p}{1-m}[p-(1+m)]}\\
&:=\overline{\c}_p(T-t)^{\frac{1}{1-m} -\frac{2\vartheta_p}{1-m}[p-(1+m)]}\,,
\end{split}\end{equation}
and Theorem  \ref{thm.GHP.up.T} gives
\begin{equation}\label{lem.lower.L1.04b}\begin{split}
\left\|\frac{u^m(t)}{\Phi_1}\right\|_{\LL^\infty(\Omega)}
&\le \ka[u_0,t] \,
(T-t)^{\frac{m}{1-m}-\frac{2\vartheta_p }{1-m}[p-(1+m)]} \,,
\end{split}
\end{equation}
where $\ka[u_0,t] $ depends  $m,p,N,\Omega$ and $u_0$ (and possibly on $t$) and is given in \eqref{thm.GHP.up.T.ineq.const}.
Note that when $t\ge T/2$ the constant $\ka[u_0,t]$ does not depend on $t$:
\[
\ka[u_0,t]=\ka[u_0]= c_{m,p,N,\Omega} \Q^*[u_0]^{\frac{4m}{1-m}\vartheta_p}\left(\frac{\|u_0\|_{\LL^p(\Omega)}^{2p\vartheta_p}}{T^{N\vartheta_p}}\right)^{2\vartheta_p(p-(1+m))_+}
\]

Combining \eqref{lem.lower.L1.01}, \eqref{lem.lower.L1.04}, and \eqref{lem.lower.L1.04b}  we obtain
\[
\frac{\rd}{\dt}\int_\Omega u\Phi_1\dx
\le -\frac{\lambda_1(T-t)^{\frac{p}{1-m}}}{\c_p^{\frac{p}{1-m}}\|u\|_{\LL^\infty(\Omega)}^{p-2m}\left\|\frac{u^m}{\Phi_1}\right\|_{\LL^\infty(\Omega)}}
\le  -\frac{\lambda_1 (T-t)^{\frac{m}{1-m}+\frac{2\vartheta_p }{1-m}(p-2m+1)[p-(1+m)]} }{\c_p^{\frac{p}{1-m}}\overline{\c}_p^{p-2m}\ka[u_0]}\,,
\]
which integrated over $[t,T]$ gives, letting $\alpha_p=\tfrac{2\vartheta_p }{1-m}(p-2m+1)[p-(1+m)]  $,
\[
\int_\Omega u(t)\Phi_1\dx  \ge \frac{\lambda_1}{\c_p^{\frac{p}{1-m}}\overline{\c}_p^{p-2m}\ka[u_0]}
\frac{(T-t)^{\frac{1}{1-m}+\alpha_p} }{\frac{1}{1-m}+\alpha_p}
:=\underline{\c}(u_0)(T-t)^{\frac{p}{1-m}-\frac{p-m}{1-m}\frac{\vartheta_{p}}{\vartheta_{1+m}}+1}
\]
which is \eqref{lem.lower.L1.ineq} when   $m \le m_s$.   Notice that
\[
\overline{\c}_p^{p-2m}\ka[u_0]=  \overline{\c}^{p-2m}c_{m,p,N,\Omega} \left[\Q^*[u_0]^{\frac{4m}{1-m}\vartheta_p}\left(\frac{\|u_0\|_{\LL^p(\Omega)}^{2p\vartheta_p}}{T^{N\vartheta_p}}\right)^{2\vartheta_p(p-(1+m))_+}\right]^{p-2m+1}\,.
\]
This concludes the proof. \qed

\begin{thm}[Global lower bounds II, the case $m>m_s$]\label{thm.GHP.low.II}
Let $u$ be a   WDS to \eqref{CDP}    corresponding to $0\le u_0\in \LL^{m+1}(\Omega)$ and let $T=T(u_0)>0$ be its extinction time.  Then we have
\[
\frac{u^m(t,x)}{\Phi_1(x) }\ge  \kb[u_0,t]  \left(\frac{T-t}{T}\right)^{\frac{m}{1-m}}\qquad\mbox{for all $0\le t\le T$ and all $x\in \overline{\Omega}$,}
\]
where   $\kb[u_0,t]=c_0 \underline{\c}[u_0,t]$,   $c_0$ is the constant in the lower bound of the Green function \eqref{G4}, and $\underline{\c}(u_0)$ is given in \eqref{c0.low}.
\end{thm}

\noindent {\bf Proof.~}Inequality \eqref{thm.GHP.low.I.ineq} and the convexity of $a\mapsto a^{\frac{1}{1-m}}$ imply
\[
\frac{u^m(t,x)}{\Phi_1(x) } 
\ge   c_0\left(\frac{t}{T}\right)^{\frac{m}{1-m}}\frac{\|u(t)\|_{\LL^1_{\Phi_1}(\Omega)}}{T -t }
\ge c_0 \left(\frac{t}{T}\right)^{\frac{m}{1-m}} \underline{\c}[u_0,t] \left(\frac{T-t}{T}\right)^{\frac{m}{1-m}}\,,
\]
where we have used the lower bound for the $L^1_{\Phi_1}$-norm \eqref{lem.lower.L1.ineq} with $p=1+m$, i.e., inequality \eqref{lem.lower.L1.03b}.\qed

Note that in the very good FDE range, $m\in (m_{c,1},1)$, the constant $\underline{\c}(u_0)$ does not depend on $u_0$, see the explicit form  \eqref{c0.low}, and we obtain lower bounds independent on the initial datum, namely:
\begin{cor}[Global lower bounds II, the case $m>m_{c,1}$]\label{thm.GHP.low.II.b}
Let $m\in (m_{c,1},1)$ and $u$ be a   WDS to \eqref{CDP}    corresponding to $0\le u_0\in \LL^{m+1}(\Omega)$, and let $T=T(u_0)>0$ be its extinction time.  Then  for all   $\frac{2}{3}T\le t\le T$   and all $x\in \overline{\Omega}$
\begin{equation}\label{thm.GHP.low.II.b.ineq}
\frac{u^m(t,x)}{\Phi_1(x) }\ge \kb \left(\frac{T-t}{T}\right)^{\frac{m}{1-m}}\,,
\end{equation}
where $\kb =c_0 \underline{\c}$, $c_0$ is the constant in the lower bound of the Green function \eqref{G4}, and $\underline{\c}$ is given in \eqref{c0.low} and does not depend on $u_0$.
\end{cor}
In the very fast diffusion regime, the decay exponent becomes more involved:
\begin{thm}[Global lower bounds III, the case $0<m\le m_s$]\label{thm.GHP.low.III}
Let $m\in (0,1)$ and $u$ be a   WDS to \eqref{CDP}    corresponding to  $0\le u_0\in \LL^p(\Omega)$ with $p\ge 1$ if $m\in (m_c,1)$ and $p>p_c$ if $m\in (0,m_c]$, where $m_c=\frac{N-2}{N}$ and $p_c=\tfrac{N(1-m)}{2}$.
Let $T=T(u_0)>0$ be the extinction time. Then  for all $0\le t\le T$ and all $x\in \overline{\Omega}$
\begin{equation}\label{thm.GHP.low.III.ineq}
\frac{u^m(t,x)}{\Phi_1(x) }\ge   \kb[u_0,t]   (T-t)^{\frac{  m}{1-m}+\frac{2\vartheta_p}{1-m}(p-2m+1)\,[p-(1+m)]_+},
\end{equation}
  where
\[
\kb[u_0,t]=c_0 \underline{\c}(u_0) \left(\frac{t}{T}\wedge \frac{2}{3}\right)^{\frac{m}{1-m}}
\]
and   $c_0$ is the constant in the lower bound of the Green function \eqref{G4}, and $\underline{\c}(u_0)$ is given in \eqref{c0.low}.
\end{thm}
\noindent\textbf{Remark. }The above Theorem indeed holds for all $m\in (0,1)$, but it relevant only when $m\le m_s$, and reduces to Theorem \ref{thm.GHP.low.II} when $m>m_s$, with simplified constant and optimal time decay.

\noindent {\bf Proof.~}We combine inequality \eqref{thm.GHP.low.I.ineq} with the lower bound for the $L^1_{\Phi_1}$-norm \eqref{lem.lower.L1.ineq}
\[
\begin{split}
\frac{u^m(t,x)}{\Phi_1(x) }&\ge \frac{c_0 t^{\frac{m}{1-m}}\|u(t)\|_{\LL^1_{\Phi_1}(\Omega)}}{(1-m)\left(T^{\frac{1}{1-m}}-t^{\frac{1}{1-m}}\right)}
\ge c_0  \left(\frac{t}{T}\right)^{\frac{m}{1-m}} \frac{\|u(t)\|_{\LL^1_{\Phi_1}(\Omega)}}{T -t }\\
&\ge c_0\underline{\c}[u_0,t]\left(\frac{t}{T}\wedge \frac{2}{3}\right)^{\frac{m}{1-m}}   (T-t)^{\frac{  m}{1-m}+\frac{2\vartheta_p}{1-m}(p-2m+1)\,[p-(1+m)]_+}\mbox{\,.\qed}
\end{split}
\]
 
\subsection{Global Harnack Principle for all times}\label{sec.GHP}
Once we have at our disposal global in space and time upper and lower estimates we can combine them into forms of Global Harnack Principle:
\begin{thm}[Global Harnack Principle I]\label{thm.GHP.I}
Let $m\in (0,1)$ and $u$ be a   WDS to \eqref{CDP}    corresponding to  $0\le u_0\in \LL^p(\Omega)$ with $p\ge 1$ if $m\in (m_c,1)$ and $p>p_c$ if $m\in (0,m_c]$, where $m_c=\frac{N-2}{N}$ and $p_c=\tfrac{N(1-m)}{2}$.
Let $T=T(u_0)>0$ be the extinction time.  Then  for all $0\le t\le T$ and all $x\in \overline{\Omega}$
\begin{equation}\label{thm.GHP.I.ineq}
 \kb[u_0,t]   (T-t)^{\frac{m}{1-m} + \frac{p-2m+1}{1-m}2\vartheta_p[p-(1+m)]_+}
\le \frac{u^m(t,x)}{\Phi_1(x) }\le
\ka[u_0,t] \,
(T-t)^{\frac{m}{1-m}-  \frac{2\vartheta_p }{1-m}[p-(1+m)]_+},
\end{equation}
  Here, $\ka[u_0,t],\kb[u_0,t] $ depend on $m,p,N,\Omega$ and $u_0$ and have the form
\begin{equation}\label{thm.GHP.I.ineq.KA}
 \ka[u_0,t]:= c_{m,p,N,\Omega} \left\{\begin{array}{lll}
   \ka   \|u_0\|_{L^p(\Omega)}^{2p\vartheta_p} \, t^{-\frac{2p\vartheta_p}{1-m}} & \quad\mbox{when }t\in\left(0,\frac{2}{3}T\right)\,,\\
   \Q^*[u_0]^{\frac{4m}{1-m}\vartheta_p}
    \left(\frac{\|u_0\|_{\LL^p(\Omega)}^{2p\vartheta_p}}{T^{N\vartheta_p}}\right)^{2\vartheta_p(p-(1+m))_+} &  \quad\mbox{when }t\in\left[\frac{2}{3}T, T\right]\,,
 \end{array}\right.
\end{equation}
where $c_{m,p,N,\Omega},  \ka  >0$ depend only on $m,p,N,\Omega$, and
\begin{equation}\label{thm.GHP.I.ineq.KB}
\kb[u_0,t]=c_0 \underline{\c}[u_0,t] \left(\frac{t}{T}\wedge \frac{2}{3}\right)^{\frac{m}{1-m}}
 =  \left( \frac{t}{T}\wedge \frac{2}{3} \right)^{\frac{m+2p(p-m+1)\vartheta_p}{1-m}}
\begin{cases}
\underline{\c}_1[u_0,T] & \mbox{when }t\in (0,\frac{2}{3}T)\\
\underline{\c}_2[u_0,T] &\mbox{when } t\in\left[\frac{2}{3}T, T\right]\,,
\end{cases}
\end{equation}
where $\underline{\c}_1[u_0,T]$ and $\underline{\c}_1[u_0,T]$ depend on only on $m,p,N,\Omega$  and $u_0,T$ and have an explicit expression given in \eqref{c0.low-1} and \eqref{c0.low-2} respectively. When $m>m_{c,1}$ we have that for $t\ge \frac{2}{3}T$, the constants $\ka[u_0,t],\kb[u_0,t] $ do not depend on $u_0$ nor on $T$, only on $m,p,N,\Omega$.

\end{thm}
\noindent {\bf Proof.~}Combine the upper bounds \eqref{thm.GHP.up.T.ineq} with the lower bounds \eqref{thm.GHP.low.III.ineq}.   The case when $t\ge \frac{2}{3}T$, follows by the upper bound \eqref{thm.GHP.up.T.ineq.1} combined with the lower bound \eqref{thm.GHP.low.III.ineq}, recalling that $\underline{\c}_2[u_0,T]$ does not depend on $u_0$ nor on $T$, see \eqref{c0.low-2}.\qed
\begin{rem}[Constants and exponents in different regimes]\label{Rem.GHP}
\rm  We consider here large times, namely $t\ge \frac{2}{3}T$,  so that $\ka[u_0,t]$ and $\kb[u_0,t]$ do not depend on $t$. \vspace{-3mm}
\begin{enumerate}[leftmargin=16pt, label=(\roman*)]\itemsep2pt \parskip2pt \parsep0pt
\item\textit{Novelty, sharpness of the exponents and of the data. } The above result is new in the subcritical case $m\in (0,m_s]$, to the best of our knowledge. The GHP \eqref{thm.GHP.I.ineq} gives the sharp spatial behaviour for all $m\in (0,1)$, namely that
    \[
    u(t,x)\asymp\Phi_1^{1/m}(x)\asymp\dist(x,\partial\Omega)^{\frac{1}{m}}\qquad\mbox{for all }(t,x)\in (0,T)\times\overline{\Omega}\,.
    \]
    It also provides the sharp time decay in the supercritical range $m>m_s$, see point (iii) below.  With respect to the GHP of DiBenedetto-Kwong-Vespri \cite{DKV1991}, valid only when $m>m_s$, we obtain an improvement in terms of the initial data, see point (iii).
    When $m>m_{c,1}$, we also eliminate the dependence on $u_0$ in the constants, and admit nonnegative $\LL^1_{\Phi_1}$ data, see (iv) below and Theorem \ref{thm.GHP.II}. On the other hand, in the subcritical range $m\in (0,m_s]$, we have that the upper and lower powers of time do not match, but this seems to be optimal. The GHP in the form \eqref{thm.GHP.I.ineq} is the only result known so far in the subcritical case, and in general it can not be improved (as far as we know): see point (ii) and Theorem \ref{thm.NoGHP} below.

\item  In the critical case $m=m_s$, i.e., the case of the Yamabe exponent, Sire-Wei-Zheng \cite{SWZ2022} show that there are data for which the solution extinguish at the following rate
\begin{equation}\label{GK+SWZ.result.0}
u(t,x)\sim  (T-t)^{\frac{1}{1-m}} \left|\log\left( T-t \right)\right|^{\frac{1}{2m}}.
\end{equation}
The same behaviour was pioneered by Galaktionov-King \cite{GK2002} on the unit ball. We refer to Section \ref{Ssec.Yamabe} for more details.

In this case, we cannot take $p=m+1=p_c$ in \eqref{thm.GHP.I.ineq}, but we can take any $p$ arbitrarily close to $m+1$, and our result implies that for all $\varepsilon>0$ (think $\varepsilon\sim p-(m+1)$)
\[
\kb[u_0](T-t)^{\frac{m}{1-m} + \varepsilon}  \le \frac{u^m(t,x)}{\Phi_1(x) }\le
\ka[u_0] \,
(T-t)^{\frac{m}{1-m}-\varepsilon} \,,
\]
  and the constants $\kb[u_0], \ka[u_0]$ may degenerate or blow-up respectively, when   $\varepsilon\to 0^+$. This is essentially sharp due to the possible logarithmic behaviour described above. Also, we are not aware of a proof that shows that all solutions must decay as in \eqref{GK+SWZ.result.0}.

\item When $m>m_s$ the time exponents and the constants considerably simplify and  \eqref{thm.GHP.I.ineq} becomes
\[
\kb[u_0](T-t)^{\frac{m}{1-m}} \le \frac{u^m(t,x)}{\Phi_1(x) }\le
\ka[u_0] \,
(T-t_0)^{\frac{m}{1-m}} \,,
\]

and we recover the GHP of \cite{DKV1991}, with the improvement that we have computable constants $\ka[u_0,t]$ and $\kb[u_0]$  that do not depend on $\|\nabla u_0^m\|_{\LL^2(\Omega)}<\infty$ as in \cite{DKV1991}, but only $u_0\in\LL^{1+m}(\Omega)$. More precisely,  $\ka[u_0,t]$ and $\kb[u_0]$  have explicit expressions given in the proofs and  depend on $m,N,p,\Omega$ and on $\Q^*[u_0]$. The latter is finite when  $u_0\in\LL^{1+m}(\Omega)$,  thanks to the  HLS inequality \eqref{HLS}.

\item When $m>m_{c,1}$ the constants do not even depend on $u_0$ as in the following:\vspace{-6mm}
\end{enumerate}
\end{rem}
\begin{thm}[Global Harnack Principle II]\label{thm.GHP.II}
Let $m\in (m_{c,1},1)$, let $u$ be the   WDS to \eqref{CDP}    corresponding to $0\le u_0\in \LL^1_{\Phi_1}(\Omega)$, and let $T=T(u_0)$ be its extinction time. Then there   constants $\ka,\kb>0$ depending only on $m,N,\Omega,\varepsilon$ such that  for all $\varepsilon\in (0,1)$ and all $\varepsilon T\le t\le T$ and all $x\in \overline{\Omega}$
\[
\kb\,(T-t)^{\frac{m}{1-m}} \le \frac{u^m(t,x)}{\Phi_1(x) }\le
\ka\,(T-t)^{\frac{m}{1-m}} \,.
\]
\end{thm}
\noindent {\bf Proof.~}Combine the upper bounds \eqref{thm.GHP.up.T.ineq.1} with the lower bounds \eqref{thm.GHP.low.II.b.ineq}.   It is clear from the previous proofs, that the constants $\ka[u_0,t]$ and $\kb[u_0,t]$ can be replaced by $\ka[u_0,\varepsilon T]$ and $\kb[u_0,\varepsilon T]$, a small change in the proofs that avoids the small times case. Notice that in the limit $\varepsilon\to 0^+$, the constants degenerate: $\ka\to\infty$ and $\kb\to 0^+$.\qed

We now rely on results from the next section to show that in general, when $m\le m_s$, it is impossible to have a GHP with matching powers up to the extinction time. As a consequence, the GHP of Theorem \ref{thm.GHP.I} is optimal in arbitrary domains. \vspace{-2mm}
\begin{thm}[Obstruction to a GHP in the subcritical case]\label{thm.NoGHP}Let $m\in (0,m_s]$ and $\Omega\subset\RR^N$ be star-shaped\,\footnote{Indeed this theorem holds true for any domain in which the existence of positive stationary solutions fails. Here, we consider star-shaped domains because we use the celebrated Pohozaev identity \cite{P1965}.}. Consider the class $\mathcal{B}$ of WDS solutions $u$ starting from  $0\le u_0\in \LL^p(\Omega)$ with $p>p_c\ge 1+m$, where $p_c=\tfrac{N(1-m)}{2}$, and let $T=T(u_0)>0$ be their extinction time.
Then it is impossible to obtain a GHP of the form
\begin{equation}\label{thm.NoGHP.ineq}
\kb\,(T-t)^{\frac{m}{1-m}} \le \frac{u^m(t,x)}{\Phi_1(x) }\le \ka\,(T-t)^{\frac{m}{1-m}}
\end{equation}
that holds true when $t\to T^-$ for all $x\in \Omega$ and for all solutions in the class $\mathcal{B}$, with constants $\ka\ge \kb>0$ that are allowed to depend on $u_0, m, N,p,\Omega$.
\end{thm}
\noindent {\bf Proof.~}When $m=m_s$ a counterexample to the GHP \eqref{thm.NoGHP.ineq}  has been recently proven by Sire-Wei-Zheng \cite{SWZ2022}, as discussed in Remark \ref{Rem.GHP}(ii), see also Section \ref{Ssec.Yamabe}.
We provide here a proof by contradiction valid in the whole range $m\in (0,m_s]$.

Assume by contradiction the validity of \eqref{thm.NoGHP.ineq} for a solution $u(\tau, x)\in \mathcal{B}$. On the one hand, the corresponding rescaled solution  $w(t,x)$  (according to \eqref{rescaling.FDE.1}) is a solution to \eqref{RCDP} and satisfies
\[
0<\kb \le \|w^m(t)\|_{\LL^q(\Omega)}\le \ka<\infty \qquad\mbox{for all $t>0$ and all $q\ge1$.}
\]
On the other hand, we know by Lemma \ref{Lem.Omega.limits} that $w^m(t)$ converges strongly in $\LL^q$   (along a subsequence $t_j\to\infty$, for any $q\in [1,2^*)$)   to a function $S$, with $S^m\in H^1_0(\Omega)$ (that may different depending on the subsequence), which is an energy solution to the associated elliptic problem \eqref{SDP} $-\Delta S^m =\c S$. Since $\|S^m\|_{\LL^q(\Omega)}\ge \kb>0$, this implies in particular the existence of a nontrivial nonnegative solution to  $-\Delta S^m =\c S$. Also, by the strong maximum principle, $S>0$ inside $\Omega$. This gives a contradiction, since we know by the results of Pohozaev \cite{P1965} that, if $\Omega$ is star-shaped, there are no positive solutions to this equation. \qed

\subsection{Boundary Harnack inequalities for all times}\label{Ssec.BHI-T}
The GHP of the previous section can be rewritten in different forms that have the flavour of Boundary Harnack Inequalities:
\begin{thm}[Boundary Harnack inequalities]\label{thm.BHI.1}
Let $m\in (0,1)$ and $u$ be a   WDS to \eqref{CDP}    corresponding to  $0\le u_0\in \LL^p(\Omega)$ with $p\ge 1$ if $m\in (m_c,1)$ and $p>p_c$ if $m\in (0,m_c]$, where $m_c=\frac{N-2}{N}$ and $p_c=\tfrac{N(1-m)}{2}$, and let $T=T(u_0)>0$ be its extinction time.  Then  for all $0\le t\le T$
\[
\sup_{x\in \Omega}\frac{u^m(t,x)}{\Phi_1(x) }
\le \frac{\overline{\mathsf h}[u_0,t]}{(T-t)^{\frac{2\vartheta_p }{1-m}[p+2(1-m)][p-(1+m)]_+}}
\inf_{x\in \Omega}\frac{u^m(t,x)}{\Phi_1(x) } \,,
\]
where   $\overline{\mathsf h}[u_0,t]=\ka[u_0,t]\,\kb[u_0,t]^{-1}$, $\ka[u_0,t]$ is given in \eqref{thm.GHP.I.ineq.KA} and $\kb[u_0,t]^{-1}$ is given in \eqref{thm.GHP.I.ineq.KB}, and depend on $m,N,p,\Omega$ and $u_0,T$.

When $t\ge  \frac{2}{3}T$,  $\overline{\mathsf h}[u_0]$ does not depend on $t$, and when $m>m_{c,1}$ it becomes independent on $u_0,T$ and only depends on $m,N,\Omega$.
\end{thm}
\noindent {\bf Proof.~}Follows from inequality \eqref{thm.GHP.I.ineq} and by the expression of the constants.\qed


\section{Asymptotic behaviour}\label{Sec.Asypt.beh}
Once we know that solutions exists globally in time (up to the extinction time), and that they are positive in the interior and smooth up to the boundary, the natural question is how they will behave ``for large times'', i.e., as they approach the extinction time:
\begin{center}
  \textit{What is the behaviour of nonnegative solutions to the \eqref{CDP} close to the extinction time?}
\end{center}
The answer to this simple question is quite involved and strongly depends on the value of $m$. We shall distinguish 3 ranges: the supercritical case $m\in (m_s,1)$, the critical or Yamabe case $m=m_s=\frac{N-2}{N+2}$, and the subcritical case $m\in (0,m_s)$.

\subsection{The supercritical case $m>m_s$}
As we shall discuss below, in this range there are bounded stationary solutions, and sharp stabilization results have been established for the first time by the authors in \cite{BF2021}.  After that, a number of improvements appeared, providing alternative proofs and extending the convergence to different norms by means of global in time regularity estimates \cite{Ak2021, JX2019, JX2021, CMCS2022}. The purpose of this subsection is to put these results in context. Also, we shall give a sketch of the proof of the nonlinear entropy method developed in \cite{BF2021}, together with some improvements. We also report on nowadays classical and recent results about uniqueness and non-degeneracy of solutions to the Lane-Emden-Fowler equation, that for us represents the  stationary equation.

\subsubsection{Convergence to stationary profiles.} This question has been addressed for the first time by Berryman-Holland in their pioneering paper \cite{BH}: under some regularity assumptions (that nowadays we know to be true in view of the previous discussion) they were able to confirm what they saw in the experiments that motivated the model, i.e., that after a short time, solutions behave like a separate variable solution until they extinguish:
\begin{equation}\label{Sep.var.soln.FDE}
u(t,x)\sim  \mathcal{U}(t,x):= \left(\frac{T-t}{T}\right)^{\frac{1}{1-m}}S(x)\qquad\mbox{as $t\to T^-$,}
\end{equation}
where $T=T(u_0)$ is the extinction time of $u(t,\,\cdot\,)$ and $S$ is a solution to \textit{the stationary problem}
\begin{equation}\label{SDP}\tag{SDP}
		\left\{\begin{array}{lll}
			-\Delta S^m= \c S& \qquad\mbox{in }\Omega\,,\\
			S=0 & \qquad\mbox{on }\partial\Omega\,,\\
			S>0  & \qquad \mbox{in } \Omega\,,
		\end{array}\right.
\end{equation}
where $\c=1/[(1-m)T]>0$. Setting $p=1/m$ and $S^m=V$, we recognize the celebrated semilinear elliptic equation of Lane-Emden-Fowler \cite{F1931,L1869}
\begin{equation}\label{LEF}\tag{LEF}
		\left\{\begin{array}{lll}
			-\Delta V= \c V^p& \qquad\mbox{in }\Omega\,,\\
			V=0 & \qquad\mbox{on }\partial\Omega\,,\\
			V>0  & \qquad \mbox{in } \Omega\,,
		\end{array}\right.
\end{equation}
We make a short digression to analyze the properties of the stationary solutions relevant to the study of the asymptotic behaviour of the \eqref{CDP}.
\subsubsection{The Lane-Emden-Fowler (LEF) problem }\label{Ssec.LEF}

The semilinear elliptic equation $-\Delta V=\c V^p$ was introduced by Lane \cite{L1869} and  Emden-Fowler \cite{F1931}, and  is probably the most studied nonlinear elliptic PDE. We focus on the more relevant contribution related to the study of the asymptotic behaviour of the fast diffusion flow.

\noindent\textbf{Existence. }Existence of stationary solutions $S$, or $V$, can be guaranteed by nowadays standard methods of Calculus of Variations, see for instance \cite{S2008}. Notice that nontrivial solutions may fail to exists for some values of $m$, or $p$, namely when $m\in (0,m_s]$, or $p\geq p_s$: this heavily depends on the geometry of the domain. Indeed, Pohozaev \cite{P1965} showed non-existence in star-shaped domains, while Bahri-Coron \cite{BC1988} proved existence when the geometry is nontrivial.

\textit{The existence of positive (nontrivial) stationary states is guaranteed in the exponent range to $m\in (m_s,1)$, or equivalently $1<p<p_s=\tfrac{N+2}{N-2}$. }

It is hopeless to account for a complete bibliography about Lane-Emden-Fowler equations  \eqref{LEF}: besides the above mentioned papers, we refer to the monographs by Quittner-Souplet \cite{QS2019}, Cazenave-Haraux \cite{CH1988} and Struwe \cite{S2008}; see also \cite{BH,BGV2012-MilanJ, BGV2013-ContMath, BGV2012-JMPA,BF2021,JX2021,Ak2021} and references therein.   We mention the celebrated papers of Brezis-Turner \cite{BT}, Gidas-Spruck \cite{GS81}, DeFigueredo-Lions-Nussbaum \cite{dFLN},   for absolute upper bounds (i.e., $0\le u\le C$ in $\Omega$, with  $C>0$ is independent of the solution) in different ranges of parameters. Constructive proofs of local and global Harnack inequalities and absolute bounds have been obtained by Grillo, V\'azquez, and the first author \cite{BGV2012-MilanJ, BGV2013-ContMath}. Global inequalities of the form
\[
V(x)\asymp \dist(x,\partial\Omega)\asymp \Phi_1(x)\qquad\mbox{for all }x\in \overline{\Omega\,, }
\]
hold true with explicit (upper and lower) constants, see \cite{BGV2012-MilanJ, BGV2013-ContMath}. Note that also the GHP of Theorems \ref{thm.GHP.I} and \ref{thm.GHP.II}  provides explicit and simpler expressions of the constants, when applied to separation of variables solutions \eqref{Sep.var.soln.FDE}.  These quantitative results were essential in the asymptotic analysis by Grillo, V\'azquez, and the first author \cite{BGV2012-JMPA}, in order to obtain convergence in relative error and the first rates of convergence to equilibrium.

\noindent\textbf{Boundary Regularity. }Since nonnegative solutions turn out to be bounded and strictly positive in $\Omega$, they are smooth and continuous up to the boundary, i.e., $C^\infty(\Omega)\cap C^{\alpha}(\overline{\Omega})$ by standard elliptic regularity,   for some
$\alpha\in [0,1]$ that may depend on the regularity of the boundary of $\Omega$.    The boundary regularity depends on the regularity of $\Omega$ and can be improved to $C^{2+\frac{1}{m}}(\overline{\Omega})$ or even $C^\infty(\overline{\Omega})$. We have discussed above the (optimal) boundary regularity of the solutions to the parabolic problem: the spatial regularity is the same, since we can apply the parabolic result to solutions of separation of variables of the form \eqref{Sep.var.soln.FDE}, see Section \ref{Ssec.Regularity}.

\noindent\textbf{Uniqueness of solutions to the Lane-Emden-Fowler equation: an open problem. }This is a longstanding and mostly unsolved open question: \textit{whether or not \eqref{LEF} has a unique solution when $\Omega$ is convex}, see for instance \cite[Remark 6.9(ii)]{QS2019}. In 1988, Dancer \cite{Dancer88} conjectured that the answer is affirmative.

To the best of our knowledge, the conjecture has been verified only for $p$ close to $1$ for $N\ge 2$ by Lin \cite{L1994}, and if
$p$ is close to $+\infty$ for $N = 2$ by De Marchis-Grossi-Ianni-Pacella \cite{dMGIP2019} and by Grossi-Ianni-Luo-Yan \cite{GILY2022}.
Some other results about uniqueness have been proven when the domain $\Omega$ possesses some kind of symmetries.
For example, when is a ball, the famous symmetry result by Gidas-Ni-Nirenberg \cite{GNN} shows that any solution to \eqref{LEF} is radial, so that uniqueness follows by ODE techniques, see also Ni-Nussbaum \cite{NN}. When the domain is both symmetric and convex with respect to $N$ orthogonal directions  Damascelli-Grossi-Pacella \cite{DGP1999} show that uniqueness hold in dimension $N = 2$, and Grossi \cite{G2000} in dimension $N\ge 3$ when $p$ is sufficiently close to $p_s$; in these papers, the symmetry and special ``convexity'' play an essential role.

In a recent preprint \cite{LWZ2022}, Li-Wei-Zu have shown that solutions of \eqref{LEF} are unique on a smooth bounded convex domain $\Omega$, provided the Robin function is a Morse function (i.e., its critical points are non-degenerate) when $p$ is sufficiently close to $p_s$. They also provide some (more technical) conditions for non-convex domains, and prove uniqueness of solutions that blow up at one or many points.

\noindent\textbf{Nondegeneracy of stationary solutions. }When we analyze the stabilization properties of the nonlinear flow, an important aspect is the non-degeneracy of the solution to \eqref{LEF}, which is related to its linearization. Again, this depends essentially on the geometry of the domain, and Saut-Temam showed in \cite{ST1979} that this happens generically. This plays an essential role in the asymptotic behaviour of solutions that will be discussed thoroughly in Sections \ref{Ssec.Isolated} and \ref{Ssec.Generic}.

\subsubsection{Logarithmic time rescaling of the nonlinear flow} In order to better understand the asymptotic results, it is convenient to consider rescaled solutions, so that the solution by separation of variables becomes stationary: let $T=T(u_0)>0$ be the extinction time, and define \begin{equation}\label{rescaling.FDE.1}
u(\tau,x) = \left(\frac{T-\tau}{T}\right)^{\frac{1}{1-m}}w(t,x), \qquad
t= T\log\left(\frac{T}{T-\tau}\right).
\end{equation}
In this way, the time interval $0<\tau<T$ becomes $0<t<\infty$, and the Problem \eqref{CDP} is mapped to
\begin{equation}\label{RCDP}\tag{RCDP}
\begin{split}
\left\{\begin{array}{lll}
w_{t}=\Delta (w^m)+\dfrac{w}{(1-m)T} &\mbox{for all }(t,x)\in(0,\infty)\times\Omega\,,\\
w(0,x)=u_0(x) & \mbox{for all }x\in \Omega\,,\\
w(t,x)=0 & {\rm for}~  t >0 ~{\rm and}~ x\in\partial\Omega.
\end{array}\right.
\end{split}
\end{equation}
The  transformation \eqref{rescaling.FDE.1} can also be expressed as
\[
w(t,x)=\ee^{\frac{t}{(1-m)T}}\, u\left(T-T\ee^{-t/T},x\right),
\]
so that \textit{the behaviour near extinction (i.e., as $\tau\to T^-$) for the original flow corresponds now to the behaviour as $t\to\infty$ in the rescaled flow. }It is quite common in the literature to use also the variable $v=w^m$, so that letting $p=1/m$ we get
\begin{equation}\label{RCDP.v}\tag{RCDP-V}
\begin{split}
\left\{\begin{array}{lll}
\partial_t v^p = \Delta v +  \c v^p\; &\mbox{for all }(t,x)\in(0,\infty)\times\Omega\,,\\
v(0,x)=u^m_0(x) & \mbox{for all }x\in \Omega\,,\\
v(t,x)=0 & {\rm for}~  t >0 ~{\rm and}~ x\in\partial\Omega\,,
\end{array}\right.
\end{split}
\end{equation}
where $  \c=\frac{1}{(1-m)T}$.
The main result of Berryman-Holland \cite{BH} reads: \it for every $u_0^m\in H^1_0(\Omega)$\,, there exists a sequence $t_n\to \infty$ such that
\[
\|v(t_n)-V\|_{H^1_0(\Omega)}\xrightarrow[n\to\infty]{} 0\,.
\]
\rm
where $V$ solves \eqref{LEF}.
However, this first result left many questions open. On the one hand, if the stationary state is unique\footnote{The stationary state of the Porous Medium case $m>1$ is always unique, see \cite{Ar-Pe, V2,BFV2016,BSV2015}, and this makes the analysis much simpler. Sharp results can be proven by comparison as in Aronson-Peletier \cite{Ar-Pe} or  by entropy methods as in V\'azquez \cite{V2}. A finer entropy method has been developed by Grillo, V\'azquez, and the first authors in \cite{BGV2012-JMPA}. A direct proof (merely using the GHP) of the sharp convergence in relative error has been done by Sire, V\'azquez, and the authors in \cite{BSV2015,BFV2016}, in the case of nonlocal porous medium equations. We refer to the introduction of our paper \cite{BF2021} for more details. } the result is quite satisfactory: the convergence holds for all times $t\to \infty$, hence the behaviour of $v$ for long times is well described by $V$ in the $H^1$ topology, and also in stronger topologies thanks to the regularity results of Jin-Xiong \cite{JX2019,JX2021}. The next question in this case would be whether or not there are (sharp) convergence rates.

On the other hand, we have seen that the problem of characterizing for which domains there hold uniqueness of solutions to the \eqref{LEF} remains open nowadays. For instance, it can happen (when $\Omega$ is a suitable annulus) that there are non-isolated stationary states. In this case, the result of Berryman-Holland \cite{BH} only guarantees that solutions are approaching the $\omega$-limit, the set of (infinitely many) positive solutions to \eqref{LEF}. A priori we can not even guarantee that given an initial datum, the solution will converge  to a unique stationary state: we cannot exclude that along different subsequences, the solution converges to different stationary states, i.e., it asymptotically oscillates between different equilibria, see also \cite{Ak2016,Ak-Ka2014,FS2000, BV2007, BGV2012-JMPA, BF2021} for further discussions.

\subsubsection{Stabilization towards a unique profile} Given the possibility of multiple different stationary states, a natural delicate question left open in \cite{BH} was to understand whether or not the solution $v$ converges to a unique stationary profile. The answer was given by Feireisl-Simondon in \cite{FS2000}. By means of a  Lojasiewicz-type inequality\footnote{Inspired and adapted from the celebrated work of Simon \cite{S1983} on stabilization of gradient flows under suitable analyticity conditions, see also Jendoubi \cite{J1998} and Akagi \cite{Ak2021}.} they proved that a nonnegative bounded weak solution to \eqref{RCDP} converges uniformly towards a unique stationary profile $S$. More precisely, we rephrase here Theorem 3.1 of \cite{FS2000}, adapted to our setting: \it Any nonnegative weak solution $w\in\LL^{\infty}\big((0,\infty)\times\Omega\big)$ of \eqref{RCDP}  is continuous for all $t>0$, and there exists a classical solution $S$ to \eqref{SDP}, depending on the initial datum, such that $w(t)\to S$ as $t\to \infty$ in the strong $C^0(\overline{\Omega})$ topology.

\rm The arguments of \cite{FS2000} heavily rely on compactness, hence no rates of convergence were provided. There has been recently some refinements and relevant improvements as we shall discuss in Section \ref{Ssec.SharpRates}, in particular see Remark \ref{rem.other.results}.

An intriguing open question is  \textit{how the initial datum selects the stationary solution. }For instance, in the case of the Cauchy problem on the whole space, the particular equilibrium (Barenblatt profile, in selfsimilar variables) is determined by the conservation of mass (at least in the good range, otherwise by conservation of relative mass) and by the   conservation of the   center of mass (when finite), see for instance \cite{BBDGV,BDGV,BDNS2022}.  In the whole space, the asymptotic profiles (Barenblatt solutions) are explicit and this is a clear advantage. On bounded domains we are not aware of explicit asymptotic profiles, nor of preserved quantities that would eventually allow to understand this selection process, still obscure nowadays.

\subsubsection{Uniform Convergence in Relative Error}
In 2012, the first author together with Grillo and V\'azquez established convergence in relative error to a stationary profile. Theorem 2.1 of \cite{BGV2012-JMPA} adapted to our current setting reads:   \it

Let $m \in (m_s,1)$, let $w$ be a bounded solution to  \eqref{RCDP} corresponding to the initial datum  $0\le u_0\in \LL^p(\Omega)$ with $p\ge 1$ if $m\in (m_c,1)$ and $p>p_c$ if $m\in (0,m_c]$, where $m_c=\frac{N-2}{N}$ and $p_c=\tfrac{N(1-m)}{2}$.    Let $S$ be the positive classical solution to the elliptic problem \eqref{SDP} such that $\|w(t)-S\|_{\LL^\infty(\Omega)}\to 0$ as $t\to\infty$.  Then
\begin{equation}\label{Rel.Err.Conv.RCDP}
\lim_{t\to\infty}\left\|\frac{w(t,\cdot)}{S(\cdot)}-1\right\|_{\LL^{\infty}(\Omega)}=0\,.
\end{equation}
\rm
The proof is based on the GHP of DiBenedetto-Kwong-Vespri \cite{DKV1991}, the uniform convergence of Feiresl-Simondon \cite{FS2000}, and a delicate barrier argument close to the boundary.  The GHP of \cite{DKV1991} forces the initial data to be such that $0\le u_0^m\in H^1(\Omega)$, but with the new results of this paper the $H^1$ assumption can be removed according to the weaker assumptions of   Theorem \ref{thm.GHP.I}.

The limit \eqref{Rel.Err.Conv.RCDP} can be equivalently stated as follows: there exists a positive function $\delta(t)\to 0$ as $t\to \infty$, such that
\begin{equation}\label{Rel.Err.Conv.RCDP.explicit}
[1-\delta(t)] S(x)\le  w(t,x)  \le [1+\delta(t)] S(x)\qquad\mbox{for all $x\in \Omega$ and all $t\ge t_0\ge 0$.}
\end{equation}
  Here $t_0>0$ is such that $\delta(t)\le 1$ for all $t\ge t_0$.
We can rephrase \eqref{Rel.Err.Conv.RCDP} and \eqref{Rel.Err.Conv.RCDP.explicit} in original variables: \\
\it there exists  a positive function $\delta(\tau)\to 0$ as $\tau\to T^-$, such that
\[
\begin{split}
\lim_{\tau\to T^-}\left\|\frac{u(\tau,\cdot)}{\mathcal{U}(\tau,\cdot)}-1\right\|_{\LL^{\infty}(\Omega)}& =0\quad\mbox{or,   there exists $\tau_0\in [0,T)$ such that} \quad\\
|u(\tau,x)-\mathcal{U}(\tau,x)|&\le \delta(\tau)S(x)\left(\frac{T-\tau}{T}\right)^{\frac{1}{1-m}}\quad\mbox{for all $x\in \Omega$ and $\tau\in [\tau_0,T]$,}
\end{split}
\]
where $\mathcal{U}(\tau,x)=S(x)\left(\frac{T-\tau}{T}\right)^{\frac{1}{1-m}}$ is the solution by separation of variables \eqref{Sep.var.soln.FDE}. \rm

\noindent\textsl{An improved Global Harnack Principle. }The latter inequality  can be rephrased as: there exists a positive function $\delta(\tau)\xrightarrow[]{\tau\to T^-}0$ such that for $\tau\in [\tau_0,T]$
\begin{equation}\label{UCRE.vs.GHP}
(1-\delta(\tau)) S(x)\left(\frac{T-\tau}{T}\right)^{\frac{1}{1-m}}\le u(\tau,x)\le (1+\delta(\tau)) S(x)\left(\frac{T-\tau}{T}\right)^{\frac{1}{1-m}}\,.
\end{equation}
This is an asymptotic improvement of the GHP: on the one hand, we get the sharp boundary behaviour, as for the GHP, recalling that $S^m\asymp \dist(\,\cdot\,,\partial\Omega)$. On the other hand,   the ``constants'' $1\pm\delta(\tau)$   become asymptotically sharp, i.e., they both converge to $1$ as $\tau\to T^-$.

In rescaled variables, we see that solutions to \eqref{RCDP.v} satisfy, for all $t \gg1$,
\[
(1-\delta(t)) V(x) \le v(t,x)\le (1+\delta(t)) V(x)\,.
\]
and this will play a crucial role in the sharp asymptotic analysis.
\subsubsection{The asymptotic regime and the linearized problem}The estimates \eqref{UCRE.vs.GHP}  in general do not hold for small times: it may take some time for the solution to become uniformly small in relative error. After that time, that may depend on the initial datum, all solution enter the ``asymptotic regime'' in which the behaviour is dictated by a suitable linearized problem, and where it is possible to measure quantitatively the convergence rates towards the equilibrium (uniquely chosen by the initial datum).

In order to have quantitative rates of convergence, it is fundamental to perform a thorough analysis of the linearized problem: we follow the $p$-notation here. According to \cite{BF2021} the linearization around a stationary solution $V$ (i.e., a solution to \eqref{LEF}) is the following:
\begin{equation}\label{linearized.V-FDE}
pV^{p-1}\partial_t f =\Delta f + \c pV^{p-1} f\,.
\end{equation}
We shall briefly analyze the fine asymptotic behaviour of the homogeneous Cauchy-Dirichlet problem for linearized equation \eqref{linearized.V-FDE}, since it provides a key-tool for the nonlinear entropy method, i.e., the validity of \textit{improved Poincar\'e inequalities} under appropriate orthogonality conditions. We briefly resume the main results and ideas, referring to \cite[Section 2]{BF2021} for a more detailed discussion and complete proofs.   Solutions to the linearized problem are regular up to the boundary, see \cite{JX2019} and references therein.

Let us begin by noticing a trivial yet important fact: $V$ \textit{is not a stationary solution }to equation \eqref{linearized.V-FDE}:
\[
-\Delta V = \c V^p\neq \c pV^p\qquad\mbox{since $p>1$.}
\]
Next we observe that stationary solutions $\varphi$ must satisfy the homogeneous Dirichlet problem associated to the linear elliptic equation
\begin{equation}\label{linearized.V-FDE.elliptic}
 -\Delta \varphi = \c pV^{p-1} \varphi\,.
\end{equation}
Whether  or not the above linear elliptic equation (an elliptic Schr\"odinger equation with potential $V^{p-1}$) admits nontrivial solutions will be essential for the understanding of the asymptotic behaviour of the linear flow \eqref{linearized.V-FDE}. This is related to non-degeneracy of solutions to the elliptic Dirichlet problem for \eqref{LEF}, and to the best of our knowledge, still an open question, that we shall briefly discuss below.

\subsubsection{Fast VS Slow equilibria}\label{Ssec.Isolated}
When uniqueness of stationary solutions fails, it can happen that there are infinitely many solutions  of \eqref{LEF}, as for the case of some annuli \cite{Ak2016,Ak-Ka2014,Ak-Ka2013,Ak2013,DGP1999, Dancer90}, and that some of these solutions may be degenerate (in other words, one can find solutions $V$ of \eqref{LEF} for which \eqref{linearized.V-FDE.elliptic} admits a nontrivial solution $\varphi$).
As we shall see, this complicates the panorama when investigating the convergence to equilibrium of FDE flows.

On the one hand, we know that $u_0$ selects only one stationary point where the solution convergences, but the selection mechanism is not explicit and this makes the convergence analysis rather involved, see for instance Akagi \cite{Ak2016, Ak2021} and Akagi-Kagikiya \cite{Ak-Ka2014}. On the other hand, when a solution to \eqref{LEF} is non-degenerate (i.e., \eqref{linearized.V-FDE.elliptic} admits no non-trivial solutions),  we have shown in \cite{BF2021} that, in rescaled variables, solutions to the \eqref{RCDP.v} converge exponentially fast to those equilibria, see \eqref{Thm.Main.Rel.Err}.
Because of this fact, we shall call non-degenerate stationary solutions {\it ``the fast'' equilibria of the FDE flow}.
As we shall explain later, this is always the case in generic domains.

In a recent work by Jin-Xiong \cite{JX2021}, it is shown that the convergence to non-isolated critical points can be slower, i.e., polynomial  decay in rescaled variables. Also, a recent contribution by Choi-McCann-Seis \cite{CMCS2022} shows that, when critical points are degenerate, the convergence must be at best polynomial. We shall call such stationary solutions {\it ``the slow'' equilibria of the FDE flow. }We will discuss in more detail these convergence results in Section  \ref{Ssec.SharpRates} and in Remark \ref{rem.other.results}.

\subsubsection{The kernel of the linearized operator is generically trivial}\label{Ssec.Generic}A natural question arises: can we characterize the domains which only have fast equilibria? Or, equivalently
\begin{center}
  \textit{For which domains $\Omega\subset\RR^N$ all the solutions of the Lane-Emden-Fowler are non-degenerate?}
\end{center}
Indeed, a characterization of such domains is still missing nowadays. There are many intriguing partial results that we discuss below, but we wish to start from a result by Saut-Temam \cite{ST1979}: the answer is affirmative for generic domains. Let us be more precise. Fix $\alpha \in (0,1)$, and define the set
\begin{equation}
\label{eq:C2a sets}
\mathcal{O}:=\left\{\Omega\subset\RR^N\,:\,\Omega\mbox{ is open}, \,\mbox{$\overline{\Omega}$ is compact, and $\partial\Omega\in C^{2,\alpha}$}\right\}.
\end{equation}
The topology on $\mathcal{O}$ can be defined through a family of neighborhoods as follows:   let $\varepsilon>0$ and
\[
\mathcal{N}_\varepsilon(\Omega):=\left\{\Omega'\in \mathcal{O}\,:\,\exists\;\Phi\in C^{2,\alpha}(\RR^N;\RR^N)\mbox{ with $\|\Phi-{\rm Id}\|_{C^{2,\alpha}}< \varepsilon$ s.t.   $\Omega'=\Phi(\Omega)$} \right\}.
\]
The question above can be rewritten as follows:
\begin{center}
\noindent $(H_\Omega)$ \it For any $V$ solution to \eqref{LEF}, there is no nontrivial solution (i.e., $\varphi\not\equiv 0$) to the equation
\[
-\Delta \varphi = \c p V^{p-1}\varphi\quad\mbox{in }\Omega\,,\qquad \varphi=0\quad\mbox{on }\partial\Omega\,,
\]
\end{center}
\rm Defining the weighted space $\LL^2_V(\Omega)$ via the norm $\|f\|_{\LL^2_V(\Omega)}^2:=\int_\Omega f^2V^{p-1}\dx$, this can also be equivalently restated:
\begin{center}
$(H_\Omega)$\qquad \it For any $V$ solution to \eqref{LEF}, $\c p$ is not an eigenvalue for the \\
\,\,\,\,\,\,\,\,Dirichlet Laplacian on $\LL^2_V$, i.e., $\c p\not\in {\rm Spec}_{\LL^2_V(\Omega)}(-\Delta)$.
\end{center}\rm
This fact is not easy to check in general, and it depends on the geometry of the domain. However, Saut-Temam \cite{ST1979} showed that this result is generically true. More precisely, let $\mathcal O$ be defined as in \eqref{eq:C2a sets} endowed with the $C^{2,\alpha}$ topology.
We define the family of (good) sets for which $(H_\Omega)$ holds:
\[
\mathcal{G}:=\{\Omega\in \mathcal{O}\;:\; \c p\not\in {\rm Spec}_{\LL^2_V(\Omega)}(-\Delta) \}\,.
\]
The result of Saut-Temam \cite[Theorem~1.2]{ST1979}, reads:
\it The set $\mathcal{G}\subset \mathcal{O}$ is open and dense.\rm

\noindent\textsl{Some examples and related results. }\hfill(See also \cite[Section 2.4]{BF2021})\vspace{-2mm} 
\begin{enumerate}[label=(\roman*),leftmargin=16pt]\itemsep1pt \parskip1pt \parsep0pt

\item\textit{Positive answer for domains with symmetries. } By the results of Damascelli-Grossi-Pacella \cite[Theorem 4.2]{DGP1999}, we know that $(H_\Omega)$ is true on balls, namely that $B_r(x_0)\in \mathcal{G}$ for all $x_0\in \RR^N$ and $r>0$, for any $N\ge 2$.
    Also, by the results of Zou \cite{Zou} we know that $(H_\Omega)$ is stable under $C^1$ perturbation of the balls.

    In dimension $N=2$, Damascelli-Grossi-Pacella \cite{DGP1999} show that $(H_\Omega)$ holds for domains which are convex in the directions $e_i$ and symmetric with respect to the hyperplanes $\{x_i=0\}$, $i=1,2$, and Grossi \cite{G2000} in dimension $N\ge 3$, but always when $p$ is sufficiently close to $p_s$.
\item\textit{Positive answer for convex domains and $p$ close to $1$ or $p_s$. }When $p$ is close to $1$ for $N\ge 2$, condition $(H_\Omega)$ has been proven to hold by Lin \cite{L1994}; if $p$ is close to $+\infty$ for $N = 2$, $(H_\Omega)$ was shown by De Marchis-Grossi-Ianni-Pacella \cite{dMGIP2019} and by Grossi-Ianni-Luo-Yan \cite{GILY2022}.
    In a recent preprint \cite{LWZ2022}, Li-Wei-Zou have proved that solutions of \eqref{LEF} are non-degenerate on a smooth bounded convex domain $\Omega$, provided the Robin function is a Morse function, i.e., its critical points are non-degenerate.

\item\textit{Negative answer on some annuli. } We know that $(H_\Omega)$ is not true for some annuli, see for instance \cite{Ak2016,Ak-Ka2014,Ak-Ka2013,Ak2013,DGP1999, Dancer90}. However, the result of Saut-Temam \cite{ST1979} implies that if we perturb a bit the annulus in the $C^{2,\alpha}$ topology, then for most perturbations $(H_\Omega)$ holds true. A similar phenomenon happens for a dumb-bell shaped domain, \cite{Dancer88,Dancer90}. More in general, if the geometry is nontrivial (as for the annulus), the results of Bahri-Coron \cite{BC1988} indicate that degenerate solutions can exist.

\end{enumerate}
\subsubsection{Sharp rates. Non-degeneracy VS weighted Poincar\'e inequalities}\label{Ssec.Seeking.Rates}
When $V$ is non-degenerate, i.e., when condition $(H_\Omega)$ is fulfilled, we can study the spectral properties of the linearized elliptic operator: this is nothing but the Dirichlet Laplacian $-\Delta$ as a linear unbounded selfadjoint operator on $\LL^2_{V^{p}}(\Omega)$.

We summarize here the results of Lemma 2.1 of \cite{BF2021}, about the spectrum of $-\Delta$ on $\LL^2_{V^{p}}(\Omega)$: \it
%
\vspace{-2mm}\begin{enumerate}[leftmargin=*]\itemsep1pt \parskip1pt \parsep0pt
\item[$(i)$]The inverse operator $(-\Delta)^{-1}:\LL^2_V\to \LL^2_V$ is a compact operator   with eigenvalues $\{\mu_{V,k}\}_{k\in \NN}$ such that $0<\mu_{V,k}\to 0^+$  as $k\to \infty$. We denote by $\V_k\subset \LL^2_V$ the finite dimensional spaces generated by the eigenfunctions associated to the $k^{th}$ eigenvalue, and  by $\pi_{\V_k}: \LL^2_V\to \V_k$ the projection on the eigenspace $\V_k$. We also denote by $N_k={\rm dim}(\V_k)$ and by $\pkj$ with $j=1,\dots, N_k$ the elements of a basis of $\V_k$ made of normalized eigenfunctions, $\|\pkj\|_{\LL^2_V}=1$.
\item[$(ii)$]The operator $ -\Delta $ is a linear unbounded selfadjoint operator on $\LL^2_V$, which is the Friedrichs extension associated to the Dirichlet form $Q(f)=\int_\Omega |\nabla f|^2\dx$. $-\Delta$ has a discrete spectrum on $\LL^2_V$, with the same eigenfunctions (and consequently the same eigenspaces $\V_k$) as $(-\Delta)^{-1}$ and eigenvalues $\lambda_{V,k}=\mu_{V_k}^{-1}$\,, so that
    \[
    0<\lambda_{V,1}<\lambda_{V,2}<\dots<\lambda_{V,k}<\lambda_{V,k+1}\to \infty
    \]
\item[$(iii)$]The smallest eigenvalue $\lambda_{V,1}=\c>0$ is simple, namely the corresponding eigenspace $\V_1$ is 1-dimensional, i.e., $N_1=1$.
Also the first positive eigenfunction is $\phi_{1,1} = V/\|V\|_{\LL^2_V}=V/\|V\|_{\LL^{p+1}}^{(p+1)/2}$.
\item[$(iv)$] All the eigenfunctions are of class $C^{2,\alpha}(\Omega)\cap C^{\alpha}(\overline{\Omega})$ for some $\alpha\in (0,1)$, and have a similar boundary behaviour: for all $x\in \Omega$ there are constants $c_{j,k,\Omega}>0$ such that
\[
\begin{split}
c_{1,1,\Omega}^{-1}\,\dist(x,\partial\Omega)\le \phi_{1,1} &\le c_{1,1,\Omega}\,\dist(x,\partial\Omega)
\qquad\mbox{and}\qquad |\pkj(x)|  \le c_{j,k,\Omega}\,\phi_{1,1}(x).
\end{split}
\]
\end{enumerate}\rm
\noindent\textbf{Seeking the sharp rates. } The assumption $(H_\Omega)$ (i.e., $\c p\not\in {\rm Spec}_{\LL^2_V(\Omega)}(-\Delta)$) guarantees that we can always define the integer $k_p>1$ as the largest integer $k$ for which $p\c>\lambda_{V,k}$, so that
\[
0<\lambda_{V,1}=\c<\dots<\lambda_{V, k_p}<p\c<\lambda_{V, k_p+1}.
\]
and show the validity of the following \textit{Improved Poincar\'e Inequality}, which is Corollary 2.2 of \cite{BF2021}: \it
Under assumption $(H_\Omega)$, let $\varphi\in \LL^2_V$ be such that
\begin{equation}\label{cor.Poincare2.hyp.OG}
\varphi_k=\pi_{\V_k}(\varphi)=0\qquad\mbox{for all }k\le k_p \,.
\end{equation}
Then the following inequality holds true:
\begin{equation}\label{cor.Poincare1.Ineq}
0< \lambda_{V, k_p+1} \int_\Omega \varphi^2 V^{p-1}\dx \le \int_\Omega |\nabla \varphi|^2\dx\,.
\end{equation}\rm

\noindent\textbf{The sharp rates. }We can rewrite the above inequality in a form which is more useful for the Entropy method. To this end, we define
\begin{equation}\label{lambda.p.def}
\lambda_p:=\lambda_{V,k_p+1}-\c p>0\,,
\end{equation}
so that \eqref{cor.Poincare1.Ineq}  becomes (under the same assumptions \eqref{cor.Poincare2.hyp.OG})
\begin{equation}\label{cor.Poincare2.Ineq}
\lambda_p\int_\Omega \varphi^2 V^{p-1}\dx  \le \int_\Omega |\nabla \varphi|^2\dx-\c p\int_\Omega \varphi^2 V^{p-1}\dx\,.
\end{equation}\rm
It turns out that $\lambda_p$ not only provides the sharp rates of convergence for the linearized flow, but it also provides the sharp rates for the nonlinear one.

\begin{rem}\label{Rem.Poincare2}\rm
Consider \textit{the limit as $p\to 1^+$ in the above Poincar\'e inequality }\eqref{cor.Poincare2.Ineq}: it has been shown by Grillo-V\'azquez and the first author in \cite{BGV2012-JMPA,BGV2013-ContMath}, that $V\to \Phi_1$ and $\lambda_{V,1}=\c\to \lambda_1$ as  $p\to 1^+$, where $(\lambda_1, \Phi_1)$ are the first eigen-elements of the classical Dirichlet Laplacian on $\Omega$.  As a consequence, the above Poincar\'e inequality \eqref{cor.Poincare2.Ineq} becomes the ``second Poincar\'e inequality'', namely $\lambda_2\|\varphi \|_{\LL^2}^2\le \|\nabla\varphi \|_{\LL^2}^2$\,, and holds for functions orthogonal to $\Phi_1$, that is $\int_\Omega\varphi\Phi_1\dx=0$.
It follows that $\lambda_p\to\lambda_2-\lambda_1$ as $p\to 1^+$. In particular, (H2) holds true for any smooth bounded domain for $p$ sufficiently close to $1$ (he closeness depending on the domain $\Omega$).
\end{rem}

\subsubsection{The linear entropy method: orthogonality and improved Poincar\'e inequalities}The simplest method to prove the asymptotic behaviour of the linear flow \eqref{linearized.V-FDE} is by means of Fourier analysis, but this cannot be extended to the nonlinear case. Hence we briefly show how to prove the asymptotic behaviour by means of an entropy method, whose main ingredients will appear also in the nonlinear case.  Let us define the linear Entropy functional $\EL[f]$ and the linear Fisher information or linear Entropy-Production functional $\IL[f]$ as follows:
\begin{equation}\label{linear.entropy}
\EL[f]=\int_\Omega f^2 V^{p-1}\dx\qquad\mbox{and}\qquad\IL[f]=\int_\Omega |\nabla f|^2\dx - p\c\int_\Omega f^2 V^{p-1}\dx\,.
\end{equation}
It is not difficult to check that $\IL[f]$ is (minus) the time derivative along the flow of the entropy $\EL[f]$:
\begin{align}
\frac{\rd}{\dt}\EL[f(t)]&= 2\int_\Omega f(t,x) f_t(t,x) V^{p-1}(x)\dx
= \frac{2}{p}\int_\Omega f(t,x) \left[\Delta f(t,x) + \c pV^{p-1}(x) f(t,x)\right]\dx \label{linear.fisher}\\
&=-\frac{2}{p}\left(\int_\Omega |\nabla f(t,x)|^2\dx - p\c\int_\Omega f^2(t,x) V^{p-1}(x)\dx\right)
=-\frac{2}{p}\,\IL[f(t)]\nonumber
\end{align}
\noindent\textbf{The Improved Poincar\'e inequality: orthogonality and convergence rates. }The first observation is that we already need an improved Poincar\'e inequality to guarantee that the Fisher information $\IL$ is nonnegative: the first Poincar\'e inequality \eqref{cor.Poincare1.Ineq} (i.e., with constant $\lambda_{V,1}=\c$) is not sufficient to guarantee the nonnegativity of $\IL$ since $p>1$. Hence, we need the improved Poincar\'e inequality \eqref{cor.Poincare2.Ineq}:
\begin{equation}\label{cor.Poincare2.Ineq-b}
\lambda_p\EL[f]=\lambda_p\int_\Omega f^2 V^{p-1}\dx  \le \int_\Omega |\nabla f|^2\dx-\c p\int_\Omega \varphi^2 V^{p-1}\dx = \IL[f]\,,
\end{equation}
but the price to pay is the orthogonality condition \eqref{cor.Poincare2.hyp.OG} that must be preserved along the linear flow, as we shall see below. Once the validity of improved Poincar\'e inequality \eqref{cor.Poincare2.Ineq-b} along the flow is ensured we can combine it with the Entropy-Entropy Production \eqref{linear.fisher} to obtain
\[
\frac{\rd}{\dt}\EL[f(t)]=-\frac{2}{p}\IL[f(t)]\le  -\frac{2\lambda_p}{p}\EL[f(t)]\,,
\]
which finally implies the exponential decay of the Entropy:
\[
\EL[f(t)] \le \ee^{-\frac{2\lambda_p}{p}t}\EL[f_0]\,,\qquad\mbox{where }\lambda_p=\lambda_{V,k_p+1}-p\c>0\,.
\]
Hence $f(t)$ converges exponentially fast to $0$ in $\LL_V^2$.
\begin{rem}\label{Rem.Asumpt.p=1}\rm
In the limit $p\to 1^+$ the above exponential decay becomes (cf. also Remark \ref{Rem.Poincare2})
\[
\int_\Omega |f(t,x)|^2 \dx \le \ee^{-2(\lambda_2-\lambda_1)t}\int_\Omega |f_0(x)|^2 \dx\,.
\]
and holds for initial data $f_0$ orthogonal to the first eigenfunction $\Phi_1$ in the $\LL^2$ sense. This is the optimal result for the classical heat equation on bounded domains with Dirichlet boundary conditions, more specifically for the equation $f_t=\Delta f +\lambda_1f$.
\end{rem}

\noindent\textbf{The orthogonality condition is preserved along the linear flow. }In order to apply the  Poincar\'e inequality \eqref{cor.Poincare2.Ineq-b} in the above entropy method, we have to make sure that the orthogonality conditions are preserved along the linear evolution: more precisely, we want to show that
\begin{equation}\label{linear.ortogonality.flow}\begin{split}
&\mbox{If $\pi_{\V_k}(f(t_0))=0$ for all $k=1,\dots,k_p$, then $\pi_{\V_k}(f(t))=0$ for all $t\ge t_0$ and all $k=1,\dots,k_p$}\,.
\end{split}\end{equation}
Indeed, given $\psi_k\in \V_k$, we know that $-\Delta \psi_k =\lambda_{V,k}V^{p-1}\psi_k$, so we can  compute
\[
\begin{split}
\frac{\rd}{\dt} \int_\Omega f(t,x)\, \psi_k (x) \,V^{p-1}(x)\dx&
= \int_\Omega f_t(t,x)  \psi_k (x) \,V^{p-1}(x)\dx\\
& =  \frac{1}{p}\int_\Omega f(t,x) \Delta \psi_k (x) \dx  + \c \int_\Omega f(t,x)\psi_k (x) V^{p-1}(x)\dx \\
& =  \frac{p\c-\lambda_{V,k}}{p} \int_\Omega f(t,x)\psi_k (x) V^{p-1}(x)\dx\,.
\end{split}\]
As a consequence, for all $\psi_k\in \V_k$
\[
\int_\Omega f(t,x)\, \psi_k (x) \,V^{p-1}(x)\dx = \ee^{\frac{p\c-\lambda_{V,k}}{p} (t-t_0)}\int_\Omega f(t_0,x)\, \psi_k (x) \,V^{p-1}(x)\dx\,,
\]
which clearly implies \eqref{linear.ortogonality.flow}. It is important to observe that if we do not impose the orthogonality condition at the initial time, the projections of the solution eventually blow up in infinite time and with an exponential rate, namely $\left|\int_\Omega f(t,x)\, \psi_k (x) \,V^{p-1}(x)\dx\right| \to \infty$ as $t\to \infty$,  for all $\psi_k\in \V_k$, $k\in \{1,\ldots,k_p\}$.

\subsubsection{First rates of convergence. }By means of quantitative continuity arguments, exponential decay to equilibrium (i.e., inequality \eqref{Thm.Main.Rel.Err} with some $\lambda_m>0$) was shown by Grillo, V\'azquez, and the first author in \cite{BGV2012-JMPA} when $m\in (m_{\sharp}, 1)$, where $m_\sharp$ has an explicit yet involved expression, that contains the constants of Harnack inequalities for solutions of the \eqref{LEF}, computed explicitly in \cite{BGV2012-MilanJ, BGV2013-ContMath}.
For all $m\in(m_{\sharp}, 1) $ a weighted Poincar\'e inequality is shown to hold as a consequence of quantitative and constructive global Harnack estimates combined with a quantitative convergence result of nonlinear eigenvalues to the linear one. More precisely, it is shown in  \cite{BGV2012-JMPA,BGV2012-MilanJ, BGV2013-ContMath} how the nonlinear eigen-pair $(\lambda_p, V_p)$ converges to the linear one $(\lambda_1, \Phi_1)$ as $p\to 1^+$. 
For more details, see Sections 4 and 5 of \cite{BGV2012-JMPA} or also the last example at the end of Subsection 2.4 of \cite{BF2021}.  When dealing with strictly positive Dirichlet data, an entropy method similar   to the one presented in   \cite{BGV2012-JMPA} has been developed in \cite{BLMV}.

\subsubsection{Sharp rates of convergence. }\label{Ssec.SharpRates}In \cite{BF2021} we have proven the sharp convergence rates in all the range $(m_s,1)$,  under the assumption $(H_\Omega)$ of non-degeneracy of the stationary solutions, which we know to be generically true by the above discussion. More precisely, Theorem 1.2 of \cite{BF2021} reads:\\ \it
There exists an open and dense\footnote{ The topology and the set $\mathcal{O}$ has been defined in Section \eqref{Ssec.Generic}.} set $\mathcal{G}\subset\mathcal{O}$, such that for any domain $\Omega\in \mathcal{G}$ the following holds. Let $p\in (1,p_s)$, and let $v$ be a solution to Problem \eqref{RCDP.v} on $[0,\infty)\times \Omega$ corresponding to the initial datum $0\le u_0\in \LL^1_{\Phi_1}(\Omega)$ when $m\in (m_{c,1},1)$ and $0\le u_0\in \LL^{1+m}(\Omega)$ when $m\in (m_s, m_{c,1})$. 
Let $V$ be the positive classical solution to \eqref{LEF} such that $\|v(t)-V\|_{\LL^\infty(\Omega)}\to 0$ as $t\to\infty$.  Also, let  $\lambda_p>0$ be defined as in \eqref{lambda.p.def}. Then there exists $\kappa>0$ such that, for all $t>0$ large,
\begin{equation}\label{Thm.Main.Rel.Err}
\int_\Omega \left|\frac{v(t,x)}{V(x)}-1\right|^2 \, V^{p+1}(x)\dx
\le \kappa\, \ee^{-2\lambda_p\,t}
\end{equation}
and the decay rate $\lambda_p$ is sharp. Also, for all $t>0$ large we have
\begin{equation}\label{Thm.Main.Rel.Err2}
\left\|\frac{v(t,\cdot)}{V(\cdot)}-1\right\|_{\LL^{\infty}(\Omega)}\le \kappa\, \ee^{-\lambda_p t}.
\end{equation}\rm
\begin{rem}\label{rem.main.results}\rm
\begin{enumerate}[label=(\roman*),leftmargin=16pt]\itemsep1pt \parskip1pt \parsep0pt
\item We have introduced some improvements with respect to the original statement of Theorem 1.2 of \cite{BF2021}. First, we can allow for more general data, i.e., $0\le u_0\in \LL^1_{\Phi_1}(\Omega)$ when $m\in (m_{c,1},1)$ and $0\le u_0\in \LL^{1+m}(\Omega)$ when $m\in (m_s, m_{c,1})$. These are the minimal assumptions under which the GHP \eqref{thm.GHP.I} holds, which is sufficient to prove the convergence in relative error of \cite[Theorem 2.1]{BGV2012-JMPA} (this was previously based on the GHP of \cite{DKV1991}, which required $u_0^m\in H^1(\Omega)$).
    The second improvement concerns the rates in the $\LL^\infty$ norms, which originally had a rate $\lambda_p/4N$. Using the regularity estimates recently obtained by Jin-Xiong \cite{JX2019}, it is possible to eliminate the annoying factor $1/4N$. Actually much more can be obtained, as we shall discuss in Remark \ref{rem.other.results}(i).

\item In original variables, the estimates of the above Theorem can be stated as follows: there exists $T_0 \in [0,T)$ such that
\[
\left\|\frac{u^m(\tau,\cdot)}{\mathcal{U}^m(\tau,\cdot)}-1\right\|_{\LL^{\infty}(\Omega)}\le \kappa' \left(\frac{T-\tau}{T}\right)^{\frac{\lambda_p}{T}}\qquad\mbox{for all $\tau\in [T_0,T]$.}
\]
where $\mathcal{U}$ is the separate variable solution defined in \eqref{Sep.var.soln.FDE}. Also,
\[
\int_\Omega \left|\frac{u^m(\tau,x)}{\mathcal{U}^m(\tau,x)}-1\right|^2 S^{1+m}(x)\dx\le \kappa' \left(\frac{T-\tau}{T}\right)^{\frac{2}{T}\lambda_p}\qquad\mbox{for all $\tau\in [T_0,T]$.}
\]
\item\textit{ About the sharpness of $\lambda_p$. } As we have explained in Section \ref{Ssec.Seeking.Rates} the rate $\lambda_p$ turns out to be the same as in the linear case, hence no better rate shall be expected in this degree of generality.
\item As $p\to1^+$ it holds $\lambda_p\to \lambda_2-\lambda_1$, and the rate is the same as  for the linear Heat equation, see  Section 4 of \cite{BGV2012-JMPA} for further details.
\end{enumerate}
\end{rem}

\begin{rem}[More recent results and alternative proofs]\label{rem.other.results} \rm After our results \cite{BF2021} appeared in arXiv in 2019, some questions were still open. A number of alternative proofs and complementary results have recently appeared, by Jin-Xiong \cite{JX2019,JX2021}, Akagi \cite{Ak2021} and  Choi-McCann-Seis \cite{CMCS2022}. We shall comment these recent results and complete the picture.  \rm\vspace{-2mm}
\begin{enumerate}[label=(\roman*),leftmargin=16pt]\itemsep1pt \parskip1pt \parsep0pt
\item\textit{Sharp rates in stronger norms. } A first related result is due to Jin-Xiong \cite{JX2019} where, as a consequence of their remarkable regularity estimates, they were able to extend the sharp rates \eqref{Thm.Main.Rel.Err2} to $\LL^\infty$ norms  and even to stronger norms: Corollary 1.4 of \cite{JX2019} shows that, under the same assumptions as our Theorem 1.2 of \cite{BF2021}, there exists $\kappa>0$ such that, for all   $p\in (1,p_s)$   and  $t\gg1$
    \[
    \left\|\frac{v(t,\cdot)}{V(\cdot)}-1\right\|_{C^{1+p}(\overline{\Omega})}\le \kappa\, \ee^{- \lambda_p t}\,.
    \]
    Also, when $p$ is integer, for all $k\ge 0$ it holds
        \[
    \left\|\frac{v(t,\cdot)}{V(\cdot)}-1\right\|_{C^{k}(\overline{\Omega})}\le \kappa\, \ee^{- \lambda_p t},
    \]
    where now $\kappa$ depends on $k$.

\item\textit{Power like rates for all domains. }The above results show that, for generic domains, there is (sharp) exponential convergence to non-degenerate profiles. A natural question is what happens when degenerate profiles exist. The answer has been given by Jin-Xiong \cite{JX2021}, where they shows that when the condition $(H_\Omega)$ is not satisfied, there is still a convergence rate but it is no more exponential. Theorem 1.2 of \cite{JX2021} reads:
\it Let $p\in (1,p_s)$ and let $v$ be a classical solution to \eqref{RCDP.v}, with extinction time $T$. Then there exists a positive stationary solution $V$ to \eqref{LEF} and two constants $\kappa, \sigma>0$ such that
    \begin{equation}\label{Conv.Rates.power}
    \left\|\frac{v(t,\cdot)}{V(\cdot)}-1\right\|_{C^{2}(\overline{\Omega})}\le \frac{\kappa}{t^\sigma}\qquad\mbox{for all }t\ge 1\,.
    \end{equation}
    If moreover condition $(H_\Omega)$ is satisfied, there exist $\kappa, \lambda>0$ such that
    \[
    \left\|\frac{v(t,\cdot)}{V(\cdot)}-1\right\|_{C^2(\overline{\Omega})}\le \kappa\, \ee^{- \lambda  t}\qquad\mbox{for all }t\ge 1\,.
    \]
    \rm
The proof relies on the careful analysis of the evolution of the curvature term $v_t/v$ (as already commented in Sections \ref{SSec.new.proof.BC} and \ref{Ssec.Regularity}) and on a Simon-Lojasiewicz type \cite{S1983} inequality (see also Jendoubi \cite{J1998}), in the spirit of Del Pino-Saez \cite{dPS2001} and Feiresl-Simondon \cite{FS2000}.

\item\textit{An alternative proof of sharp convergence rates. }Recently, Akagi \cite{Ak2021} provided a proof of the sharp convergence rates for the energy
    \[
    J[v]=\frac{1}{2}\int_\Omega|\nabla v|^2\dx - \frac{\c}{p+1}\int_\Omega |v|^{p+1}\dx\,.
    \]
    More precisely Theorem 1.3 of \cite{Ak2021} shows that,  under the same assumptions as our Theorem 1.2 of \cite{BF2021},  there exists $\kappa>0$ such that, for all $t>0$ large,
    \[
    0\le J[v(t)]-J[V]\le \kappa e^{-\lambda_p\, t}\,.
    \]
    \rm As a consequence, convergence in weighted $\LL^2$ norm \eqref{Thm.Main.Rel.Err} is obtained with a different proof (and then convergence in higher norms follows, as in (i) above). The proofs are essentially functional analytic, and rely on the higher differentiability of the functional $J$ on $H^1_0$, and are different from the ones of \cite{BF2021} and \cite{JX2019,JX2021}\,. The key ingredients are fine $H^{-1}$ energy estimates, a gradient-type inequality for the functional $J$, together with refined ``spectral estimates'' for the linearized operator. Also, the proofs have the advantage to be true also for signed solutions.

\item\textit{An interesting dichotomy. }When the profiles are degenerate, we know by the results of Jin-Xiong \cite{JX2021} that there is at least polynomial decay \eqref{Conv.Rates.power}. The question is: can this convergence be improved? The negative answer is provided in a recent preprint by Choi-McCann-Seis \cite{CMCS2022}, whose main result, Theorem 3.1 reads: \it Let $p\in (1,p_s)$, and let $v$ be a solution to \eqref{RCDP.v} in $[0,\infty)\times \Omega$ converging in relative error to $V$, a positive classical solution to \eqref{LEF}, i.e., $\|v(t)/V-1\|_{\LL^\infty(\Omega)}\to 0$ as $t\to\infty$. Then exactly one of the two alternatives holds:\\
    (i) The relative error decays algebraically or slower:
    \[
    \left\|\frac{v(t)}{V} -1\right\|_{\LL^\infty(\Omega)}\gtrsim  \left\|\frac{v(t)}{V} -1\right\|_{\LL^2_V(\Omega)}\gtrsim\frac{1}{t}
    \]
    (ii) The relative error decays exponentially fast:
    \[
     \left\|\frac{v(t)}{V} -1\right\|_{\LL^2_V(\Omega)}\lesssim\left\|\frac{v(t)}{V} -1\right\|_{\LL^\infty(\Omega)}\lesssim e^{-\lambda t}
    \]
    where $\lambda\ge \lambda_K$, and $\lambda_K$ is the first positive eigenvalue of a suitable linearized operator.
    \rm

    Note that when condition $(H_\Omega)$ holds, $\lambda=\lambda_p$ and they recover the sharp results of \cite{BF2021}. The novelty here is represented by the fact that the authors analyze the case in which the kernel of the linearized operator is not trivial, and found that in some cases (for instance, when the kernel of the linearized operator is integrable) exponential convergence still holds. We refer to \cite{CMCS2022} for more details.

\item\textit{Conclusion. }On generic domains, i.e., when condition $(H_\Omega)$ holds, there is the sharp exponential convergence to the isolated profiles \eqref{Thm.Main.Rel.Err}, proven in \cite{BF2021} and then in \cite{JX2021,Ak2021,CMCS2022}. When condition $(H_\Omega)$ is not satisfied, only polynomial decay rates can be guaranteed and the following dichotomy holds: either there is polynomial convergence with rate at most $1/t$, or the convergence is exponential.
\rm
\end{enumerate}
\end{rem}

\rm

\subsubsection{The nonlinear entropy method. Poincar\'e inequalities and almost orthogonality. }
Let us explain with some details the nonlinear entropy method developed by us in \cite{BF2021}, since it introduces a number of new features and ideas that can potentially be extended to different contexts. Special attention will be devoted to  the ``almost orthogonality condition'', which allows us to export sharp linear spectral results to the nonlinear settings.

Let us define the Entropy functional
\[
\E[v]=\int_\Omega \left[\left(v^{p+1}-V^{p+1}\right) -\frac{p+1}{p}(v^p-V^p)V \right]\dx\,,
\]
which will be the nonlinear analogue of the linear entropy functional $\EL[f]$ defined in \eqref{linear.entropy}.

For our entropy method to work, we will need to be in the asymptotic regime, namely we need the relative error to be small, and this is guaranteed in the range $p\in (1,p_s)$ by the results of \cite{BGV2012-JMPA}, as we have already discussed, cf. \eqref{Rel.Err.Conv.RCDP}. Hence, without loss of generality we can assume that:\\ \it for any $\delta\in(0,1)$ there exists a time $t_0>0$ such that
\begin{equation}\tag*{(H1${\rm '}$)$_\delta$}
|f(t,x)|\le \delta \,V(x)\qquad\mbox{for a.e. }(t,x)\in [t_0,\infty)\times\Omega
\end{equation}\rm
\textit{The idea of this nonlinear entropy method }is to mimic the linear case as much as possible: we differentiate the entropy along the flow, obtaining (up to errors that we have to control carefully) a nonlinear Fisher information which we can compare with the linear one, so that we can apply the (improved) Poincar\'e inequalities and conclude the (sharp) exponential decay. Although the strategy is simple, there are two difficulties.\\
First, in order to use the improved Poincar\'e inequalities we need orthogonality conditions, which are not preserved along the nonlinear flow. We solve this first difficulty with improved Poincar\'e valid under ``Almost Orthogonality'' (AO) conditions, and we show that AO is preserved along the flow.
This is the most technical and delicate part of the proof, see Sections 3.3--3.6 of \cite{BF2021} for more details.\\
The second important step is to quantitatively (and constructively) compare linear and nonlinear quantities. This requires the weighted smoothing effects of \cite[Section 4]{BF2021}, or the quantitative regularity estimates of  \cite{JX2019}.

\noindent\textbf{Comparing linear and nonlinear Entropy and Entropy production. }To begin with, we need to quantitatively compare linear and nonlinear quantities: Lemma 3.2 and Proposition 3.3 of \cite{BF2021} read:\\ \it
Let $w=v^p$  be a solution to the (RCDP), let $f=v-V$, and assume $\mathrm{(H1{\rm '})_\delta}$ with $0<\delta<1/2p$. Then, for all $t\ge t_0$ we have
\[
\frac{p+1}{2(1+\overline{c}_p\delta)^2}\,\EL[f]\,\le\, \E[v]\,\le\, \frac{p+1}{2}(1+\overline{c}_p\delta)^2\,\EL[f]\,.
\]
and also
\[
\frac{\rd}{\dt}\E[v(t)]= -\frac{p+1}{p}\,\IL[f(t)]\,+\,\R_p[f(t)]
\]
where
\begin{equation}\label{Prop.Entropy.Lin.Nonlin.2}
\big|\R_p[f]\big|  \le   \kappa_p  \int_\Omega |f|^3 V^{p-2}\dx\,.
\end{equation}\rm
At this point we need some improved Poincar\'e inequalities to conclude.

\noindent\textbf{Almost-orthogonality conditions and improved Poincar\'e inequalities. }
Recall that  we have defined $\lambda_p:=\lambda_{V,k_p+1}-\c p>0\,,$ where $k_p$ was the largest $k$ such that $p\c>\lambda_{V,k}$\,. It is convenient to express the orthogonality conditions \eqref{cor.Poincare2.hyp.OG} in an equivalent way, by means of Rayleigh-type quotients:
\[
\QL_{k,j}[\psi]:=\frac{\left|\int_\Omega \psi\,\pkj  \,V^{p-1}\dx\right|}{\left(\int_\Omega \psi^2 \,V^{p-1}\dx\right)^{\frac{1}{2}}}
=\frac{\big|\langle \psi, \pkj\rangle_{\LL^2_{V}}\big|}{\|\psi\|_{\LL^2_{V}}}=0,
\]
for all $k=1,\dots,k_p$ and $j=1,\dots,N_k$.
As we have already explained, the above orthogonality conditions are preserved along the linear flow, but not along the nonlinear flow. Therefore, we have introduced a new concept of \textit{almost-orthogonality, }which plays an analogous role for the nonlinear flow and allows us to use improved Poincar\'e inequalities along the nonlinear flow. More precisely, we say that a function $f\in\LL^2_V$ satisfies the $\varepsilon$-almost-orthogonality condition for the linear functional, (AOL)$_\varepsilon$ for short,  if the Rayleigh quotients $\QL_{k,j}$ are small: namely, \it given $\varepsilon\in (0,1)$,
\begin{equation}\label{AOL}\tag*{(AOL)$_\varepsilon$}
\QL_{k,j}[f]\le \varepsilon \qquad\mbox{for all $k=1,\dots,k_p$ and all $j=1,\dots,N_k$}\,.
\end{equation}\rm
As mentioned above we would like to show that, given $\varepsilon>0$ small,
the condition (AOL)$_\varepsilon$ holds after some time along the nonlinear flow.
This is very difficult to show, and in fact we do not know how to prove it directly. Instead, we can control  a nonlinear version of \ref{AOL}, the \textit{nonlinear Rayleigh quotients} defined  below, that we can prove to remain uniformly small along the nonlinear flow and asymptotically converge to zero:
\[
\Q_{k,j}[v]:=\frac{\left|\int_\Omega \big(v^p-V^p\big)\,\pkj  \,\dx\right|}
    {\left(\int_\Omega \left[\left(v^{p+1}-V^{p+1}\right) -\frac{p+1}{p}(v^p-V^p)V \right]\dx\right)^{\frac{1}{2}}}:=\frac{\A_{k,j}[v]}{\E[v]^{\frac{1}{2}}}\,.
\]
As we shall discuss below, the nonlinear Rayleigh quotients $\Q_{k,j}$ are quantitatively comparable to the linear ones $\QL_{k,j}$ and, as a consequence, the \ref{AOL} condition stated in terms of $\QL_{k,j}$ is essentially equivalent to the one stated in terms of $\Q_{k,j}$, namely
\begin{equation}\label{AON}\tag*{(AON)$_\varepsilon$}
\Q_{k,j}[v]\le \varepsilon \qquad\mbox{for all $k=1,\dots,k_p$ and all $j=1,\dots,N_k$}\,.
\end{equation}
The equivalence between (AOL)$_\varepsilon$ and (AON)$_\varepsilon$ follows by Lemma 3.4 of \cite{BF2021} and can be summarized in the following way: in the asymptotic regime, i.e., when the relative error is smaller than $\delta$, there exists $\kappa_{p}>1$ such that, taking $\delta\ll \varepsilon$, we have
\[
(AOL)_\varepsilon\qquad \Longrightarrow \qquad (AON)_{\kappa_{p}\varepsilon}\qquad \Longrightarrow\qquad(AOL)_{\kappa_{p}^2\varepsilon}.
\]
See Remark 3.5 of \cite{BF2021} for more details.

\noindent\textbf{Improved Poincar\'e inequality for almost-orthogonal functions }holds true once we impose the AO conditions, as in Lemma 3.6 of \cite{BF2021} that reads: \\ \it
Assume $(H_\Omega)$, and let $\varphi\in \LL^2_V$ satisfy (AOL)$_\varepsilon$.
Then, the following improved Poincar\'e inequality holds:
\[
( p\c + \lambda_p -\gamma_p\,\varepsilon^2) \int_\Omega \varphi^2 V^{p-1}\dx\le \int_\Omega|\nabla\varphi|^2\dx\,,
\]
where $\gamma_p:=(\lambda_{V,k_p+1}-\lambda_{V,1})k_pN_{k_p}$. We can also rewrite it as follows
\[
\begin{split}
(\lambda_p -\gamma_p\,\varepsilon^2)\,\EL[\varphi]&=(\lambda_p -\gamma_p\,\varepsilon^2) \int_\Omega \varphi^2 V^{p-1}\dx
 \le \int_\Omega|\nabla\varphi|^2\dx - \c p\int_\Omega \varphi^2 V^{p-1}\dx = \IL[\varphi]\,.
\end{split}\]

\rm

\noindent\textbf{Nonlinear Entropy-Entropy Production inequalities for almost orthogonal functions. }Assuming temporarily that AO is preserved along the nonlinear flow, we can combine the above improved Poincar\'e inequalities with the entropy production \eqref{Prop.Entropy.Lin.Nonlin.2}. \it Let  $w=v^p$ be a solution to the (RCDP) and let   $v=f+V$   satisfy $\mathrm{(H1{\rm '})}_\delta$ with $0<\delta<1/2p$. We have the following: \vspace{-2mm}
\begin{enumerate}[label=(\roman*),leftmargin=16pt]\itemsep1pt \parskip1pt \parsep0pt
\item Assume   $(H_\Omega)$   and that $f(t)$ satisfies (AOL)$_\varepsilon$ for $t\ge t_0$.
    Then, choosing $\delta, \varepsilon\ll 1$ so  that $\kappa_{p}\varepsilon^2+\delta<2\lambda_p/(p\tilde\gamma_p)$ for some (small  and explicit) $\tilde \gamma_p>0$, we have that
    \begin{equation}\label{Entropy.decay.quasi.OG.ineq.nonlin}
        \frac{\rd}{\dt}\E[v(t)]\le -\left(\frac{2\lambda_p}{p}-\tilde\gamma_p(\kappa_{p}\varepsilon^2+\delta)\right)\,\E[v(t)]<0\,.
    \end{equation}
\item  Assume   $(H_\Omega)$   and that, for some $\eta>0$, we have
    \begin{equation}\label{Entropy.decay.quasi.OG2.hyp}
        \left\|\frac{v(t)-V}{V}\right\|_{\LL^\infty(\Omega)}\le \ka \, \E[v(t-1)]^{\eta} \quad\mbox{and}\quad
        \QL_{k,j}[f(t)]\le \overline{c}_{p,k,j} \, \E[v(t-1)]^{\frac{\eta}{2}}\,,
    \end{equation}
    for all $t\ge t_0\ge 1$ and all $k=1,\dots,k_p$, $j=1,\dots,N_k$. Then, for all $t\ge t_0\ge 1$ we obtain
    \begin{equation}\label{Entropy.decay.quasi.OG2.ineq}
        \frac{\rd}{\dt}\E[v(t)]\le
            - \frac{2\lambda_p}{p}\,\E[v(t)] + \kappa_p\, \E[v(t-1)]^\eta \,\E[v(t)]\,.\vspace{-1mm}
    \end{equation}
\end{enumerate}\rm
Notice that \eqref{Entropy.decay.quasi.OG.ineq.nonlin} implies the (almost sharp) exponential decay of the entropy, by a simple integration. This is the result that follows from a qualitative AO condition. On the other hand, inequality \eqref{Entropy.decay.quasi.OG2.ineq} is a ordinary differential inequality with delay which allows how to show that, if we control the AO condition with a power of the Entropy, then we obtain \textit{the sharp decay of the entropy} as in Proposition 3.6 of \cite{BF2021}: \it Let (ii) above hold true. Then, there exists a $T_0\ge t_0\ge 0$ such that for all $t\ge T_0$ we have the following sharp decay rates of the entropy
\[
 \E[v(t)]\le \ka_0 \ee^{-\frac{2\lambda_p}{p} t}\,,
\]
where $\ka_0>0$ depends on $p,N, \eta, T_0, \E[v(T_0)]$.\rm

\begin{rem}\label{rem.entropy.exp}\vspace{-2mm}\rm
\begin{enumerate}[label=(\roman*),leftmargin=16pt]\itemsep1pt \parskip1pt \parsep0pt
\item On the one hand, in order to get almost sharp exponential decay rates we only need to ensure the validity of AO conditions along the nonlinear flow for large times.
\item On the other hand,  to obtain sharp decay rates we need to ensure the validity of hypothesis \eqref{Entropy.decay.quasi.OG2.hyp}, i.e., we need a quantitative control of the AO condition. This can be obtained through a weighted $\LL^2-\LL^\infty$ smoothing effect \cite[Theorem 4.1]{BF2021}. 
     In particular, Corollary 4.2 of \cite{BF2021} reads: \\ \it
    Assume $\mathrm{(H1{\rm '})}_\delta$ with $0<\delta<1/2p$ and that $t_0$ is large enough so that $\E[v(t_0)]\le 1$ and $\frac{\rd}{\dt}\E[v(t)]<0$ for all $t\ge t_0-1$. Then the following estimates hold true for any $t\ge t_0$:
    \[
        \left\| \frac{v(t)}{V}-1\right\|_{\LL^\infty(\Omega)}\le  \ka_{\infty}\E[v(t-1)]^{\frac{1}{4N}}\,,
    \]
    where $\ka_\infty>0$ depends on $N,p,\c,\Omega, \|V\|_{\LL^\infty(\Omega)}, \|V\|_{\LL^{p+1}(\Omega)}$\,.  \rm
(Thanks to the recent optimal regularity results of Jin-Xiong \cite{JX2019,JX2021,JX2022}, the exponent  $1/4N$ can be replaced by $1/2$.)
\end{enumerate}
\end{rem}\rm
Hence, in order to ensure the exponential decay of the entropy, it only remains to check that the AO conditions are true along the nonlinear flow.

\noindent\textbf{Almost orthogonality improves along the nonlinear flow. }We shall show that the AO conditions are not only preserved, but they improve along the nonlinear flow, as in Propositions 3.13 and 3.14 of \cite{BF2021} that read:
 \it Assume $(H_\Omega)$. Let  $w=v^p$ be a solution to the (RCDP) and let   $v=f+V$   satisfy $\mathrm{(H1{\rm '})}_\delta$ with $0<\delta<1/2p$. We have that the following holds: \vspace{-2mm}
\begin{enumerate}[label=(\roman*),leftmargin=16pt]\itemsep1pt \parskip1pt \parsep0pt
\item \textsl{Qualitative almost orthogonality along the nonlinear flow. }For every $\varepsilon>0$ there exists $t_\varepsilon\ge t_0\ge0$ such that if $\mathrm{(H1{\rm '})}_\delta$ holds for some $\delta<\kb_0 \varepsilon$, then
    \[
        \Q_{k,j}[v(t)]\le \varepsilon\quad\mbox{for all $t\ge t_\varepsilon$ and for   all $k=1,\dots,k_p$ and $j=1,\dots,N_k$\,.}\vspace{-2mm}
    \]
\item \textsl{Quantitative almost orthogonality along the nonlinear flow. }Assume that for some $\gamma>0$ inequality \eqref{Entropy.decay.quasi.OG2.hyp} holds.
     Then, there exists a time $T_0\ge t_0\ge 0$ such that
     \[
        \Q_{k,j}(v(t))\le  \E[v(t-1)]^{\frac{\gamma}{2}}\,,\qquad\mbox{for all $t\ge T_0$ and all $1\le k\le k_p$.}
     \]
\end{enumerate}\rm
\begin{rem}\label{rem.entropy.exp}\vspace{-2mm}\rm
The proof of these results is quite involved and technical. This is actually the core of the results of \cite{BF2021} and it contains the most delicate proofs. Roughly speaking, we are showing that the nonlinear is more stable than the linear.
    The proof is based on a quantitative contradiction argument (constructive), that exploits the fact that the nonlinear projections (i.e., the numerator of $\QL_{k,j}$)\vspace{-1mm}
    \[
        \A_{k,j}[v(t)]:=\left|\int_\Omega\left(v^p(t,x)-V^p(x)\right)\pkj (x)  \dx\right|\vspace{-1mm}
    \]
    may explode (exponentially) if the AO conditions are violated at some time, see  \cite[Lemma 3.11]{BF2021} for more details.
\end{rem}\rm

This discussion concludes our analysis of the convergence for the fast diffusion equation in the supercritical case $m>m_s$.

\subsection{The critical case $m=m_s$. The Yamabe flow}\label{Ssec.Yamabe}
When $m=m_s=\frac{N-2}{N+2}$ and $N\ge 3$, the FDE $u_t=\Delta u^m$ corresponds to the scalar curvature of the metrics along the Yamabe flow, introduced by Hamilton \cite{H1989}. The asymptotic behaviour on smooth compact Riemannian manifolds has been proved by Ye \cite{Y1994}, Schwetlick-Struwe \cite{SS2003} and Brendle \cite{B2005, B2007}. To the best of our knowledge, sharp asymptotic results are still missing for the Dirichlet problem. A pioneering result is due to Galaktionov-King \cite{GK2002} and is valid for radial solutions on a ball; see also \cite[Section 5]{K2010} for a formal discussion for general domains. More recently, a generalization of the previous result is given by Sire-Wei-Zheng in \cite{SWZ2022}, where they showed the existence of some initial data such that the solution to \eqref{CDP} extinguish with the rate (the same of \cite{GK2002} on $B_1(0)$)\vspace{-1mm}
\begin{equation}\label{GK+SWZ.result}
u(t,x)\sim  (T-t)^{\frac{1}{1-m}} \left|\log\left( T-t \right)\right|^{\frac{1}{2m}}\vspace{-1mm}\,.
\end{equation}
The proof uses a variant of the gluing method for parabolic equations, in the spirit of Cort\'azar-del Pino-Musso \cite{CDM2020} and D\'avila-del Pino-Wei \cite{DDW2020}.

\noindent\textit{A new Global Harnack Principle. }Theorem \ref{thm.GHP.I} shows that, in this case, all bounded solutions satisfy for all $t\in [\frac{2}{3}T,T]$,   $x\in \overline{\Omega}$, and $\varepsilon>0$,\vspace{-1mm}
\[
  \kb[u_0]   (T-t)^{\frac{m}{1-m} +  (p-2m+1) \varepsilon} \left(\frac{t}{T}\right)^{\frac{m}{1-m}} \le \frac{u^m(t,x)}{\Phi_1(x) }\le
\ka[u_0] \,
(T-t)^{\frac{m}{1-m}-\varepsilon} \,,\vspace{-1mm}
\]
This result is almost sharp, in view of  the logarithmic rate \eqref{GK+SWZ.result}. On the other hand we do not know if all solutions must decay as \eqref{GK+SWZ.result} (  cf.   Remark \ref{Rem.GHP}).

\noindent\textit{Blow up as $t\to \infty$ in rescaled variables. }If we reinterpret the above result \eqref{GK+SWZ.result} in terms of the rescaled solution $w(t)$, this correspond to the case of initial data that produce solutions that blow up (in the $\LL^\infty$ norm) when $t\to\infty$ at a power-like rate, namely $\|v(t)\|_{\LL^\infty}\sim t^a$ as $t\to \infty$ for some $a>0$.

\noindent\textit{Bubble towers. }Daskalopoulos-del Pino-Sesum \cite{DDS2018} have constructed a class of type II ancient solutions to the Yamabe flow, which are rotationally symmetric and converge to a tower of spheres as $t\to -\infty$. Bubble tower solutions for the energy critical heat equation were constructed in del Pino-Musso-Wei \cite{DMW2022}. Sire-Wei-Zheng in \cite{SWZ2022} conjecture the existence of bubble tower solutions also for the Dirichlet problem under study.

\noindent\textbf{Asymptotic results for small perturbation in the critical case. }It we perturb the FDE by a ``Brezis-Nirenberg term'', $u_t=\Delta u^m +b u$ with $b>0$, then stationary states exist, as shown in the celebrated paper of Brezis-Nirenberg \cite{BN1983}. Jin-Xiong \cite{JX2021} recently proved that for any $b>0$ the solution in this case stabilizes to a stationary state, as it happens in the case $m>m_s$. The proof relies on a delicate adaptation the blow up analysis of Struwe \cite{S1984}, Bahri-Coron \cite{BC1988}, Schwetlick-Struwe \cite{SS2003} and Brendle \cite{B2005}. The asymptotic results fail when $b=0$.  The entropy method of \cite{BF2021} can be adapted to this case as well, and provides explicit convergence rates in terms of the linearized problem.

\subsection{New results in the subcritical case $m<m_s$. }\label{Ssec.Asymp.Subcrit.New}
When $m\in (0,m_s)$, an asymptotic analysis performed at a formal level by King \cite[Section 5]{K2010}  found that different kinds of selfsimilar solutions seem to provide the correct asymptotic behaviour; see also Galaktionov-King
\cite[Section 1.2]{GK2002} for a brief discussion on this subcritical range. The situation here is unclear but there are some partial results that indicate that convergence rates shall be worse in general, as in the case $m=m_s$. In particular, in rescaled variables they shall imply blow up when $t\to\infty$, as in the critical case.  Our contribution in this case is given by the GHP \eqref{thm.GHP.I.ineq} close to the extinction time, namely we can show that bounded solutions (that have to start with data in $\LL^q$ with $q>p_c>1+m$ in this case) have a decay that can be estimated from above and below as follows,   for all $\tau\in [\frac{2}{3}T,T]$:
\[
\kb[u_0]  (T-\tau)^{\frac{m}{1-m} + \frac{q-2m+1}{1-m}2\vartheta_q[q-(1+m)]_+}
\le \frac{u^m(\tau,x)}{\Phi_1(x) }\le
\ka[u_0] \,
(T-\tau)^{\frac{m}{1-m}-  \frac{2\vartheta_q }{1-m}[q-(1+m)]_+},
\]
cf. also Theorem \ref{thm.GHP.I} and Remark \ref{Rem.GHP} for more precise statements. In rescaled variables, i.e., for solutions $w$ of the \eqref{RCDP}, the result reads:
\begin{equation}\label{GHP.subcrit.resc.var}
\kb[u_0]\, e^{-\frac{q-2m+1}{1-m}2\vartheta_q[q-(1+m)]_+ \frac{t}{T}} \le \frac{w^m(t,x)}{\Phi_1(x)}\le
\ka[u_0] \,
e^{\frac{2\vartheta_q }{1-m}[q-(1+m)]_+\frac{t}{T}} \,,
\end{equation}
This result is in clear contrast with the supercritical case, in which  solutions were stabilizing towards a positive profile in relative error, which means that the  quotient $w^m(t)/\Phi_1$ is bounded and bounded away from zero uniformly for large  times.

We shall see that, in the subcritical range, the GHP \eqref{GHP.subcrit.resc.var} describes quite well the behaviour of the solutions, that converge to zero in $\LL^q$ for $q$ small, while they blow up in high $\LL^q$ norms.

Our main result in the subcritical range is the following:

\begin{thm}[Asymptotic behaviour in the subcritical case, star-shaped domains]\label{thm.asympt.subcrit}
Let $m\in (0,m_s]$, $\Omega$ be a smooth star-shaped domain, and let $w$ be a solution to \eqref{RCDP} corresponding to   $u_0\in L^{1+m}(\Omega)\cap H^{-1}(\Omega)$.   Then $w^m(t)$ converges to zero strongly in $\LL^q(\Omega)$ for all $q<2^\ast$, that is:
\[
\lim_{t\to\infty}\|w(t)\|_{\LL^{qm}(\Omega)}=0\qquad\mbox{for all }0<q<2^\ast.
\]
Moreover, the $\|w(t)\|_{\LL^{1+m}(\Omega)}$ norm is uniformly bounded for all times:\footnote{Notice that we always have $qm< 1+m$ since $m\le m_s$ and $q<2^*$.   Notice also that when $q\ge p_c$ then $u_0\in L^q(\Omega)$ authomatically implies $u_0\in H^{-1}(\Omega)$ by the HLS inequality \eqref{HLS}. Recall that the dual quotient is defined as $\Q^*[u_0]= \|u_0\|_{\LL^{1+m}(\Omega)}^{1+m}\|u_0\|_{H^{-1}(\Omega)}^{-(1+m)}$.}
\begin{equation}\label{thm.asympt.subcrit.energy.1+m.bounds}
\sup_{t>0}\|w(t)\|_{\LL^{1+m}(\Omega)} \le \c_{m,\Omega}^{\frac{1}{1+m}} \Q^*[u_0]^{\frac{2m}{1-m^2}}<+\infty\,.
\end{equation}
When  $u_0\in L^q(\Omega)$ with $q\ge  p_c$, the solution blows up when $t\to\infty$ in $\LL^q$,
\[
\lim_{t\to\infty}\|w(t)\|_{\LL^{q}(\Omega)}=+\infty\,.
\]
Also, when $q>p_c$, have upper bounds on the explosion rate given by \eqref{GHP.subcrit.resc.var}.
\end{thm}
Before beginning the proof, we shall introduce some concepts and preliminary results.
\subsubsection{$H^1$ Lyapunov functional, $\omega$-limits, and proof of Theorem \ref{thm.asympt.subcrit}}
Consider the Lyapunov functional $\F: \dom(\F)\subseteq H^1_0(\Omega)\to \RR$ defined on its domain $\dom(\F)= \{w\in \LL^{1+m}(\Omega) \;|\; w^m\in H^1_0(\Omega)\} $
\[
\F[w]:= \frac{1}{2} \int_\Omega \left|\nabla w^m\right|^2\dx - \frac{\c m}{1+m}\int_\Omega |w|^{1+m}\dx\,.
\]
We first observe that this functional is bounded below when $m\ge m_s$, since $w^m\in H^1_0(\Omega)$ implies $w\in L^{1+m}(\Omega)$ by the standard Sobolev inequality in $H^1_0(\Omega)$. On the other hand, when $m<m_s$ the $\LL^{1+m}$ norm cannot be controlled by the $H^1_0$ norm of $w^m$, hence a priori the functional can be unbounded from below when we only ask $w^m\in H^1_0(\Omega)$.

In any case, the time derivative of $\F$ along the flow is given by
\[
\frac{\rd}{\dt}\F[w(t)]=-m\int_\Omega w^{m-1}(w_t)^2 \dx = - \J[w(t)] \le 0
\]
  It follows by Proposition \eqref{Prop.Q*-estimates}, estimate \eqref{ext.1+m}, that,   for all $m,\varepsilon \in (0,1)$ and $\tau\in (\varepsilon T,T]$,
\[
\|u(\tau)\|_{\LL^{1+m}(\Omega)} \le \c_{\varepsilon, m,\Omega}^{\frac{1}{1+m}} \Q^*[u_0]^{\frac{2m}{1-m^2}}
(T-\tau)^{\frac{1}{1-m}}\,.
\]
Hence, the rescaled solution $w$ satisfies (whenever $u_0\in \LL^{1+m}(\Omega)$ we have that $\Q^*[u_0]<+\infty$)
\[
\sup_{t>T_\varepsilon}\|w(t)\|_{\LL^{q}(\Omega)} \le \c_{\varepsilon,m,q,\Omega}^{\frac{1}{1+m}} \Q^*[u_0]^{\frac{2m}{1-m^2}}=: \ka_\varepsilon[u_0]<+\infty\qquad\mbox{for all $0<q\le 1+m$}
\]
where $T_\varepsilon=T|\log(1-\varepsilon)|$, and we have used H\"older inequality with $q\le 1+m$.

Recall that the gradient term is bounded for all positive times: indeed, by Proposition \ref{Prop.Q*-estimates} we have that for all $\varepsilon\in (0,1)$ and all $\tau\in [\varepsilon T,T]$
\begin{equation}\label{ext.grad}
\|\nabla u^m(\tau)\|_{\LL^{2}(\Omega)}^2 \le \ka_{0,\varepsilon}[u_0](T-\tau)^{\frac{2m}{1-m}}\,,\quad\mbox{that implies}\quad
\sup_{t>T_\varepsilon}\|\nabla w^m(t)\|_{\LL^{2}(\Omega)}^2 \le \ka'_\varepsilon[u_0]\,,
\end{equation}
according to \eqref{rescaling.FDE.1}, where $T_\varepsilon$ is the same as above. Note that $\ka_\varepsilon[u_0]$ and $\ka'_\varepsilon[u_0]$  depend on $N,m,q,\Omega,\Q^*[u_0]$ and $\varepsilon$ (and may blow up as $\varepsilon\to 0$), but do not depend on $t$.

Hence the Entropy is bounded for all times $t\ge   T_\varepsilon   >0$, namely $|\F[w(t)]|\lesssim \ka_\varepsilon[u_0]$.
Being a bounded non-increasing function along the flow, the Entropy has a limit as $t\to\infty$:
\[
-  \ka_\varepsilon[u_0]\lesssim   \F_\infty=\lim_{t\to\infty}\F[w(t)] \lesssim \ka'_\varepsilon[u_0]\,.
\]
\begin{lem}\label{Lem.Lq.I.estimates}For all $\varepsilon\in(0,1)$,   all $q\in (\tfrac{1+m}{2},1+m]$,   all $h\ge 0,$ and all $t\ge T_\varepsilon$, there exists a constant   $\ka_{\varepsilon,q}[u_0]>0$   such that
\[
\int_\Omega \left|w^q(t+h)-w^q(t)\right|\dx \le \ka_{\varepsilon,q}[u_0]\, h^{\frac{1}{2}}
\left(\int_t^{\infty}\J[w(\tau)]\rd\tau\right)^{\frac{1}{2}}\,,
\]
where
\[
\ka_{\varepsilon,q} =q  \, |\Omega|^{1-\frac{q}{1+m}}
\, \ka_{\varepsilon}[u_0]^{ q- \frac{1+m}{2}}\,.
\]
\end{lem}
\noindent {\bf Proof.~}Since $\tfrac{1+m}{2}<q$ and $t\ge T_\varepsilon$, we have
\[
\begin{split}
\int_\Omega \left|w^q(t+h)-w^q(t)\right|&\dx
= \int_\Omega \left|\int_t^{t+h}\partial_t(w^q) \rd\tau\right|\dx
\le q\int_\Omega \int_t^{t+h}w^{q-1}|w_t|  \rd\tau \dx\\
&\le q\left(\int_t^{t+h}\int_\Omega  w^{m-1}(w_t)^2 \dx\rd\tau\right)^{\frac{1}{2}}\left(\int_t^{t+h}\int_\Omega w^{2q-(1+m)} \dx\rd\tau\right)^{\frac{1}{2}}\\
&\le  q h^{\frac{1}{2}} |\Omega|^{1-\frac{q}{1+m}} \left(\sup_{t>T_\varepsilon}\|w(t)\|_{\LL^{2q-(1+m)}(\Omega)}\right)^{ q- \frac{1+m}{2}}
\left(\int_t^{\infty}\J[w(\tau)]\rd\tau\right)^{\frac{1}{2}}\,.
\end{split}\]
This concludes the proof. \qed

\begin{lem}[Convergence and characterization of the $\omega$-limit]\label{Lem.Omega.limits}
Let $m\in (0,1)$, and $w$ be a solution to \eqref{RCDP} corresponding to $u_0\in L^{1+m}(\Omega)  \cap H^{-1}(\Omega) $. Then the semiorbit of $w^m$, given by
\[
\gamma(w^m,t_1)=\{w^m(t)\in \LL^q(\Omega)\;|\; t\ge t_1\}
\]
is precompact in $\LL^q(\Omega)$ for all $q<2^\ast$ and all $t_1\ge 0$. Moreover,  $\omega$-limit set of $w$ respect to the $\LL^q(\Omega)$ strong topology defined by
\[
\mathcal{M}:=\bigcap_{t_0>0}\overline{\bigcup_{t_1\ge t_0}\gamma(w^m;t_1)}=\{S^m\in \LL^q(\Omega)\;|\;\mbox{there exist $t_k\to\infty$ such that }\|w^m(t_k)-S^m\|_{\LL^q(\Omega)}\to 0\,\}
\]
is non-void and is the set of nonnegative solutions $S$ to the stationary problem \eqref{SDP}.
\end{lem}
\begin{rem}The $\omega$-limit in the different ranges. \rm We have seen in Section \ref{Ssec.LEF} that when $m>m_s$ the $\omega$-limit set $\mathcal{M}$ can have infinitely many elements, for instance when $\Omega$ is a certain annulus. On the other hand, when $m\le m_s$ we have that $\mathcal{M}=\{0\}$ at least when the domain is star-shaped, as a consequence of Pohozaev identity \cite{P1965}. In the subcritical case the asymptotic results are more delicate, both because of the possible absence of nontrivial steady states, and since solutions starting from $u_0\in \LL^{1+m}(\Omega)$ may blow-up (in finite or infinite time).
\end{rem}
\noindent {\bf Proof.~}The proof is quite standard when $m>m_s$, see \cite{BH,FS2000,V2}, but it requires some nontrivial changes when $m\le m_s$. We provide here a proof that holds for every $m\in (0,1)$. One of the main issues when $m\le m_s$ is that the gradient part of the entropy does not control the lower part, i.e., the Sobolev inequality does not help anymore. However, as shown before, the boundedness of the functional $\F$ is ensured along the flow. Also, the entropy has a limit when $t\to \infty$, that we called $\F_\infty$, and satisfies
\begin{equation}\label{Lem.Omega.limits.01}
0\le \F[w(t)]-\F_\infty=\int_t^\infty \J[w(\tau)]\rd\tau \qquad\mbox{hence}\qquad \lim_{t\to\infty}\int_t^\infty \J[w(\tau)]\rd\tau =0\,.
\end{equation}

\noindent$\bullet~$\textsc{Step 1. }\textit{Compactness in $\LL^p$. }We fix $\varepsilon\in (0,1)$. As a consequence of Proposition \ref{Prop.Q*-estimates}, as discussed above, we know that there exists a constant $\ka'_\varepsilon[u_0]>0$ such that for all $t>T_\varepsilon$
\begin{equation}\label{Lem.Omega.limits.11}
0\le \|\nabla w^m(t)\|_{\LL^2(\Omega)}^2 + \|w(t)\|_{\LL^{1+m}(\Omega)}^2\le \ka'_\varepsilon[u_0]\,.
\end{equation}
By Kondrachov Theorem, the semiorbit of $w^m$ is precompact in all $\LL^q(\Omega)$ with $q<2^\ast$. Hence, being $\mathcal{M}$ the limit set of the semiorbits, it contains the limits as $t\to\infty$ along subsequences.

Another consequence is that the Entropy is bounded below, which ensures the validity of \eqref{Lem.Omega.limits.01}.

\noindent$\bullet~$\textsc{Step 2. }\textit{The $\omega$-limit is invariant under time shift. }For any $q\le 1+m$ and any $\tau\ge 0$ we have that
\[
\lim_{t\to \infty}\big|\|w(t+\tau)\|_{\LL^q(\Omega)}^q-\|w(t)\|_{\LL^q(\Omega)}^q\big|
\le \lim_{t\to \infty}\int_\Omega |w^q(t+\tau)-w^q(t)|\dx=0
\]
as a consequence of Lemma \ref{Lem.Lq.I.estimates}. In particular, this implies strong convergence in $\LL^q$ by   the results of Bresiz-Lieb, that can be found in \cite[Theorem 1.9]{LiebLoss}.   Hence, given a sequence $t_j\to\infty$ such that $w_j\to S$ strongly in $\LL^q$, we have that for all $\tau\ge 0$ the limit is the same for a.e. $x\in \Omega$,
\begin{equation}\label{Lem.Omega.limits.22}
\lim_{t_j\to \infty} \|w(t_j+\tau)-w(t_j)\|_{\LL^q(\Omega)}=0
\end{equation}
Define the shifted function $\tau \mapsto w_j(\tau)=w(t_j+\tau)$. On the one hand, by the above discussion, we know that any fixed $\tau\ge 0$ the sequence $w_j(\tau)$ is convergent to a limit that we call $S(\tau)$ (in the strong $\LL^q$ topology and up to subsequences that we do not relabel). On the other hand, we have by \eqref{Lem.Omega.limits.22} that the limit is the same in $\LL^q$, hence almost everywhere in $\Omega$, namely $S(\tau,x)=S(0,x):=S(x)$. This implies that the limit does not depend on $\tau$.

\noindent$\bullet~$\textsc{Step 3. }\textit{Characterization of the $\omega$-limit. }Fix $\tau> 0$, $0<t_j\to \infty$ and $w_j(\tau)=w(t_j+\tau)$, and $S(\tau)$ as above. It is clear that $w_j(\tau,x)$ satisfies the same equation as $w(\tau, x)$, hence by the very weak formulation of the equation \eqref{RCDP} it is easy to deduce that for all $\varphi\in C_c^\infty(\Omega)$ we have
$$
\int_\Omega w_j(\tau,x)\varphi(x)\dx-\int_\Omega w_j(0,x)\varphi(x)\dx
=\int_{0}^{\tau}\int_\Omega\left[w_j^m(\eta,x)\Delta \varphi(x)-\c w_j(\eta,x)\varphi(x)\right] \dx\rd\eta
$$
Letting $j\to \infty$, thanks to the $L^q$ convergence to $S$ we get
$$
\int_\Omega S(\tau,x)\varphi(x)\dx-\int_\Omega S(0,x)\varphi(x)\dx
=\int_{0}^{\tau}\int_\Omega\left[S^m(\eta,x)\Delta \varphi(x)-\c S(\eta,x)\varphi(x)\right] \dx\rd\eta.
$$
Since $S(\tau)=S(0)=S(\eta)=S$, we get
$$
\tau \int_\Omega\left[S^m(x)\Delta \varphi(x)-\c S(x)\varphi(x)\right] \dx=0\,.
$$
Hence $S$ is a weak solution to the stationary problem \ref{SDP}, since it is a very weak solution (according to the above definition), and moreover we know that $S^m\in H^1_0(\Omega)$ and $S\in \LL^{1+m}(\Omega)$ by \eqref{Lem.Omega.limits.11}.\qed

We show next that when the $\LL^q$ norm is uniformly small, rescaled solutions must extinguish in finite. Note that this is impossible, since rescaled solutions are precisely defined (via a suitable time rescaling) so that they are nontrivial on the whole time interval $[0,\infty)$.

\begin{lem}\label{Lem.Subcr}
Let $m\in (0,m_s]$, $\Omega\subset\RR^N$ be a smooth domain, and let $w$ be a solution to \eqref{RCDP} corresponding to $u_0\in L^q(\Omega)$ with $q\ge p_c$. Assume that $w$ has uniformly small  $L^q$ norm at some time: there exists $t_0\geq 0$ such that
\begin{equation}\label{Lem.Subcr.hyp}
\|w(t_0)\|_{\LL^q(\Omega)}\le\delta\,,
\end{equation}
with
\[
0\le   \delta^{1-m}   \le \frac{c_{m,q}}{2q\c}, \qquad\mbox{where}\qquad c_{m,q}=\frac{4mq(q-1)}{(q+m-1)^2}\mathcal{S}_2^{-2}
\]
and $\mathcal{S}_2$ is the constant in the Sobolev-Poincar\'e inequality \eqref{Sob.Poinc} with $p=q$.  Then $w$ must extinguish at a finite time.
\end{lem}
\noindent {\bf Proof.~}Let us estimate the time derivative of the $\LL^q$ norm as follows:
\begin{equation}
\label{eq:ODEq}
\begin{split}
\frac{\rd}{\dt}\int_\Omega w(t,x)^q\dx &=-\frac{4mq(q-1)}{(q+m-1)^2}\int_\Omega \left|\nabla w^{\frac{q+m-1}{2}}\right|^2\dx
    +q\c \int_\Omega w^q\dx\\
&\le -\frac{4mq(q-1)}{(q+m-1)^2}\mathcal{S}_2^{-2} \left(\int_\Omega w^q\dx\right)^{1-\frac{1-m}{q}}+q\c \int_\Omega w^q\dx\\
&=  -c_{m,q} \left(\int_\Omega w^q\dx\right)^{1-\frac{1-m}{q}}\left[1- \frac{q\c}{c_{m,q}} \left(\int_\Omega w^q\dx\right)^{\frac{1-m}{q}}\right]\,,
\end{split}
\end{equation}
where $\mathcal{S}_2^2$ is the constant in the Sobolev-Poincaré inequality \eqref{Sob.Poinc} with $p=q$,
which we can use since $q\ge p_c$.

Assume now that \eqref{Lem.Subcr.hyp} holds. Then we claim that $\|w(t)\|_{\LL^q(\Omega)}\le\delta$ for all $t\geq t_0$. Indeed, set
$$
t_*:=\sup\{t \in [t_0,\infty)\,: \,\|w(t)\|_{\LL^q(\Omega)}\le\delta\}
$$
and assume by contradiction that $t_*<\infty$. Since $\|w(t_*)\|_{\LL^q(\Omega)}\le\delta$, by continuity there exists $\tau>0$ such that $\|w(t)\|_{\LL^q(\Omega)}\le 2\delta$ for $t \in [t_*,t_*+\tau].$
But then, it follows by \eqref{eq:ODEq} and the definition of $\delta$ that
$$
\frac{\rd}{\dt}\int_\Omega w(t,x)^q\dx\leq 0 \qquad \text{on }[t_*,t_*+\tau],
$$
therefore,   for $t>t_*$ we have    $\|w(t)\|_{\LL^q(\Omega)}\le\|w(t_*)\|_{\LL^q(\Omega)}\le \delta$, a contradiction to the maximality of $t_*$.

Now, since we know that $t_*=+\infty$, it follows by \eqref{eq:ODEq} and the definition of $\delta$ that, for $t \geq t_0$,
$$
\frac{\rd}{\dt}\int_\Omega w(t,x)^q\dx \leq  -c_{m,q} \left[1- \frac{q\c}{c_{m,q}}   \delta^{1-m}  \right]\left(\int_\Omega w^q\dx\right)^{1-\frac{1-m}{q}}\leq -\frac{c_{m,q}}{2} \left(\int_\Omega w^q\dx\right)^{1-\frac{1-m}{q}}.
$$
Integrating the differential inequality on $[t_0,\infty)$ implies that $\|w(t)\|_{\LL^q(\Omega)}$ vanishes in finite time.\qed

We are now in the position of giving the proof of our main result of this part.

\noindent {\bf Proof of Theorem \ref{thm.asympt.subcrit}.~}We split three cases.

\noindent\textit{Decay to zero when $q<2^*$. }We know by Lemma \ref{Lem.Omega.limits} that $w(t_k)\to S$ as $t_k\to \infty$, where $S$ is a nonnegative solution to \eqref{LEF}. By the result of Pohozaev \cite{P1965} we know that $S\equiv 0$\,. Since the $\omega$-limit consists of one element, the convergence does not hold only along subsequences, but also for all $t\to\infty$. Hence we have that $w^m(t)\to 0$ as $t\to \infty$ strongly in $\LL^q$, for all $q<2^\ast$.

\noindent\textit{Boundedness when $q=1+m$. }\eqref{thm.asympt.subcrit.energy.1+m.bounds} follows by Proposition \ref{Prop.Q*-estimates} and from the change of variables \eqref{rescaling.FDE.1}.

\noindent\textit{Blow up when $q\ge p_c$. }We claim that, for $q\ge 2^\ast$, the $\LL^q$ norm explodes when $t\to\infty$. Indeed, assume by contradiction that this is false, i.e., that for some $q>2^\ast$ we have that $\liminf\limits_{t\to\infty}\|w^m(t)\|_{\LL^q(\Omega)} <\infty$. We also know that the $H^1_0$ norm is uniformly bounded by \eqref{ext.grad}, cf. also Proposition \ref{Prop.Q*-estimates}. Hence
\[
\liminf_{t\to\infty}\|w^m(t)\|_{\LL^q(\Omega)} + \|\nabla w^m(t)\|_{\LL^2(\Omega)}<\infty.
\]
By Kondrachov Theorem, this implies that there exists a sequence $t_k\to \infty$ such that $w^m(t_k)$ is precompact in $\LL^q(\Omega)$. Also, thanks to Lemma \ref{Lem.Omega.limits}, up to a subsequence we have that $w^m(t_k)\to S^m\in \mathcal{M}$ in $\LL^q(\Omega)$, and
\[
\|S^m\|_{\LL^q(\Omega)}= \lim_{t_k\to \infty}\|w^m(t_k)\|_{\LL^q(\Omega)} =  c_q\ge 0
\]
Now there appear two different cases, $c_q>0$ and $c_q=0$, that we will analyze separately.

In the case when $c_q>0$, we immediately achieve a contradiction: indeed we would have proven the existence of a nontrivial nonnegative solution $S$ (hence positive) of the stationary problem \eqref{SDP}. This contradicts the results of Pohozaev \cite{P1965}: on star-shaped domains there are no positive solution to the Dirichlet problem \eqref{SDP}.

In the case $c_q=0$, since $\|w^m(t_k)\|_{\LL^q(\Omega)}$ tends to zero, given $\delta>0$ as in Lemma~\ref{Lem.Subcr} there exists $k_q\in\NN$ such that
\[
 \|w^m(t_k)\|_{\LL^q(\Omega)} \le \delta.
\]
Lemma \ref{Lem.Subcr} implies that in this case $w(t)$ must extinguish at a finite time, which is impossible by the construction of the rescaled solution.\qed

%
%
\bigskip
\hrule
\bigskip
\addcontentsline{toc}{section}{~~~Acknowledgments}

\noindent\textbf{Acknowledgments. } M.B. was partially supported by the Projects MTM2017-85757-P and PID2020-113596GB-I00 (Ministry of Science and Innovation, Spain).  M.B. acknowledges financial support from the Spanish Ministry of Science and Innovation, through the ``Severo Ochoa Programme for Centres of Excellence in R\&D'' (CEX2019-000904-S) and by the E.U. H2020 MSCA programme, grant agreement 777822. {A.F. has received funding  from the European Research Council under the Grant Agreement No 721675. }
Essential parts of this work were done while M.B. was visiting A.F. at ETH Z\"urich (CH) in the year 2022. M.B. would like to thank the FIM (Institute for Mathematical Research) at ETH Z\"urich for the kind hospitality and for the financial support. We thank the anonymous referees for their useful suggestions on a preliminary version of our paper.
%
%
%
%


\end{document}